\documentclass[11pt,reqno]{amsart}
\usepackage{geometry}
\usepackage{amssymb,amsmath,amsthm}
\usepackage[utf8]{inputenc}
\usepackage[T1]{fontenc}
\usepackage{enumerate}
\usepackage[all]{xy}
\usepackage{fullpage}
\usepackage{comment}
\usepackage{array}
\usepackage{longtable}
\usepackage{stmaryrd}
\usepackage{mathrsfs} 
\usepackage{xcolor}
\usepackage[dvipsnames]{xcolor} 
\usepackage{mathdots}
\usepackage{mathtools}
\usepackage{tikz}
\usepackage{float}
\usepackage{siunitx}
\usepackage{subcaption}

\def\Z{\mathbb{Z}}
\def\R{\mathbb{R}}
\def\N{\mathbb{N}}
\def\Q{\mathbb{Q}}
\def\K{\mathbb{K}}
\def\F{\mathbb{F}}
\def\A{\mathcal{A}}

\def\deg{\mathrm{deg}}

\def\int{\mathrm{Int}}

\allowdisplaybreaks

\usepackage{hyperref}
\hypersetup{
    colorlinks=true,
    linkcolor=blue,
    filecolor=red,  
    citecolor=green,
    urlcolor=cyan,
    pdftitle={Escape of Mass of the p Cantor},
    pdfpagemode=FullScreen,
    }

\usepackage{cleveref}
\crefformat{section}{\S#2#1#3}
\crefformat{subsection}{\S#2#1#3}
\crefformat{subsubsection}{\S#2#1#3}
\usepackage{enumitem}
\usepackage{tikz}

\newtheorem{theorem}{Theorem}[section]
\newtheorem{lemma}[theorem]{Lemma}
\newtheorem{corollary}[theorem]{Corollary}

\newtheorem{example}[theorem]{Example}
\newtheorem{conjecture}[theorem]{Conjecture}

\theoremstyle{definition}
\newtheorem{definition}[theorem]{Definition}

\theoremstyle{remark}
\newtheorem{remark}[theorem]{Remark}

\makeatletter
\newcommand{\subalign}[1]{%
  \vcenter{%
    \Let@ \restore@math@cr \default@tag
    \baselineskip\fontdimen10 \scriptfont\tw@
    \advance\baselineskip\fontdimen12 \scriptfont\tw@
    \lineskip\thr@@\fontdimen8 \scriptfont\thr@@
    \lineskiplimit\lineskip
    \ialign{\hfil$\m@th\scriptstyle##$&$\m@th\scriptstyle{}##$\hfil\crcr
      #1\crcr
    }%
  }%
}
\makeatother

\numberwithin{equation}{section}

\newcommand{\nth}{\textup{th}}

\title{Escape of Mass of the $p$-Cantor Sequence}
\author{Noy Soffer Aranov}
\email{noy.sofferaranov@tugraz.at}
\address{Graz University of Technology, Institute of Analysis and Number Theory, 8010 Graz, Austria}

\author{Steven Robertson}
\email{
steven.robertson@postgrad.manchester.ac.uk
}
\address{Department of Mathematics, University of Manchester, Manchester, United Kingdom}
\begin{document}
\begin{abstract}
    Let $p$ be a prime. In 2017, Kemarsky, Paulin, and Shapira (KPS) conjectured that any Laurent series over $\F_p$ exhibits \textit{full escape of mass} with respect to any irreducible polynomial $P(t)\in\F_p[t]$. In 2025, this was shown to be false in the case $p=2$ and $P(t)=t$ by Nesharim, Shapira and the first named author. This work shows that for \textit{any} odd prime $p$ and \textit{any} irreducible polynomial $P(t)\in\F_p[t]$, the so-called $p$-Cantor sequence provides a counterexample to the aforementioned conjecture over $\F_p$. \\

    \noindent Furthermore, the concepts of \textit{maximal} escape of mass and \textit{generic} escape of mass are introduced. These lead to two natural variations of the KPS conjecture, both of which are shown to hold for all previous counterexamples.
\end{abstract}
\maketitle
\tableofcontents
\newpage
\section{Introduction}
\noindent Let $\alpha\in\R$ be a quadratic irrational. A classical theorem of Lagrange states that the regular con\-tinued fraction expansion of $\alpha$ is eventually periodic. For any prime $p\in \mathbb{N}$ and any $k\in\N$, it is clear that $p^k\alpha$ is also a quadratic irrational, and hence, its continued fraction expansion is also eventually periodic. That is, for every $k\in \mathbb{N}$, there exist $\ell_k,m_k\in \mathbb{N}$ and $b_{0}(p^k\alpha),\dots,b_{m_k}(p^k\alpha),a_{1}(p^k\alpha),\dots,$ $a_{\ell_k}(p^k\alpha)\in \mathbb{N}$ such that the continued fraction expansion of $\alpha$ is 
$$p^k\alpha=\left[b_{0}\left(p^k\alpha\right);\dots,b_{m_k}\left(p^k\alpha\right),\overline{a_{1}\left(p^k\alpha\right),\dots,a_{\ell_k}\left(p^k\alpha\right)}\right].$$
\noindent Above, the bar notation denotes the periodic part of the continued fraction expansion. It is natural to investigate how the continued fraction expansion of $p^k\alpha$ evolves as $k$ grows. There are two key components of this evolution - the length of the period $\ell_k$ and the magnitude of the partial quotients comprising the periodic part of the continued fraction expansion. The former line of enquiry was resolved by Aka and Shapira \cite{AS} when they proved that the period of the continued fraction expansion of $p^k\alpha$ grows exponentially.
\begin{theorem}{\cite[Theorem 1.2]{AS}}
\label{thm:ASPer}
    Let $p$ be a prime and let $\alpha$ be a quadratic irrational. Then, there exists a constant $c$ depending on $\alpha$ and on $p$, such that for every $k\in \mathbb{N}$, $$\ell_k=cp^k+o\left(p^{\frac{15}{16}k}\right).$$
\end{theorem}
\noindent Furthermore, they proved that no individual partial quotient in the periodic part of the continued fraction expansion of $p^k\alpha$ makes up a positive logarithmic proportion of the total mass as $k$ tends towards infinity (see \cite[Theorem 1.2]{AS} or \cite[page 2]{PS}):
\begin{theorem}{\cite[Theorem 1.2]{AS}}
\label{thm:ASEsc}
    For every quadratic irrational $\alpha\in \mathbb{R}$, and for every prime number $p$, 
    \begin{equation}
    \label{eqn:noEscR}
    \lim_{k\rightarrow\infty}\frac{\max_{i=1,\dots,\ell_k}\log\left(a_i(p^k\alpha)\right)}{\sum_{i=1}^{\ell_k}\log(a_i(p^k\alpha))}=0.
\end{equation}
\end{theorem}
\noindent If (\ref{eqn:noEscR}) holds, then the sequence $\{p^k\alpha\}_{k\in \mathbb{N}}$ is said to exhibit \textbf{no escape of mass}. \\

\noindent In order to prove Theorems \ref{thm:ASPer} and \ref{thm:ASEsc}, the aforementioned authors utilized the well established connection between the continued fraction expansion and the geodesic flow on $X=\operatorname{PGL}_2(\mathbb{R})/\operatorname{PGL}_2(\mathbb{Z})$ \cite{EW,AS}, which is identified with the space of lattices up to homothety. Moreover, they investigated the dynamics of the diagonal group $A\subseteq \operatorname{PGL}_2(\mathbb{R})$ on $X$ to obtain results about the evolution of continued fraction expansions.  \\

\noindent Given a compact $A$ orbit, $Ax\subseteq X$, it is understood that $Ax$ supports a unique $A$ invariant probability measure $\mu_{x}$. By Dirichlet's units theorem \cite[Chapter 3]{CF}, periodic $A$ orbits correspond to quadratic irrationals. For a quadratic irrational $\alpha$, let $x=x_{\alpha}\in X$ be the associated lattice. In \cite{AS}, the authors analysed the evolution of measures $\mu_{p^k\alpha}:=\mu_{x_{p^k\alpha}}$, where $p$ is prime and $\alpha$ is a quadratic irrational to obtain Theorems \ref{thm:ASPer} and \ref{thm:ASEsc}. In order to do so, they examined weak star limits of measures corresponding to nested lattices satisfying specific arithmetic relationships.
\begin{theorem}{\cite[Theorem 4.8(A)]{AS}}
\label{thm:ASHeckeGen}
     Let $\Lambda_0\in X$ be a lattice with a compact $A$ orbit. Let $\{\Lambda_n\}_{n\in \mathbb{N}\cup\{0\}}\subseteq X$ be a sequence of lattices, such that for every $n\geq 0$, one has $\Lambda_{n+1}\subseteq \Lambda_n$ and $\Lambda_n/\Lambda_{n+1}\cong \mathbb{Z}/p\mathbb{Z}$ as $\mathbb{Z}$ modules. Then, besides possibly two exceptional sequences, the measures $\mu_{\Lambda_n}$ converge to $m_X$ in the weak star topology, where $m_X$ is the Haar measure supported on $X$ satisfying $m_X(X)=1$.
\end{theorem}
\begin{remark}
In fact, \cite[Theorem 4.8(A)]{AS} is more general than Theorem \ref{thm:ASHeckeEsc}. However, its general form is not relevant to this paper and is therefore omitted. 
\end{remark}
\noindent In particular, when $\Lambda_k$ corresponds to the lattice $\Lambda_{p^k\alpha}$, one obtains the following theorem, which implies Theorem \ref{thm:ASEsc}. 
\begin{theorem}{\cite[a specific case of Theorem 4.8(A)]{AS}}
\label{thm:ASHeckeEsc}
    Let $p$ be prime and let $\alpha$ be a quadratic irrational. Then, the measures $\left\{\mu_{p^k\alpha}\right\}_{k\in \mathbb{N}}$ converge to $m_X$ in the weak star topology.
\end{theorem}
\noindent Kemarsky, Paulin, and Shapira \cite{KPS,PS} investigated function field analogues of Theorems \ref{thm:ASPer}, \ref{thm:ASEsc}, \ref{thm:ASHeckeGen}, and \ref{thm:ASHeckeEsc}. In order to state these results, the function field setting is introduced. 
\subsection{Continued Fraction Expansions in $\mathbb{F}_q(\!(t^{-1})\!)$}\label{Sect: EoM_intro}\hfill\\
\noindent Let $q$ be a positive power of a prime, and let $\F_q$ be the unique finite field of cardinality $q$. Define the field of formal Laurent series over $\mathbb{F}_q$ by
$$\mathbb{F}_q(\!(t^{-1})\!)=\left\{\sum_{n=-h}^{\infty}a_nt^{-n}:a_n\in \mathbb{F}_q,h\in \mathbb{Z},a_{-h}\neq 0\right\}.$$
\noindent Similarly, define the ring of polynomials over $\F_q$ as 
$$\mathbb{F}_q[t]=\left\{\sum_{n=0}^ha_nt^n:a_n\in \mathbb{F}_q,h\in \mathbb{N},a_h\neq 0\right\}.$$
In the function field setting, the analogue of $\R$ ($\Z$, respectively) is given by $\F_q(\!(t^{-1})\!)$ ($\F_q[t]$, respectively). In line with this analogy, the function field version of $\Q$ is given by the field of rational functions: \[\F_q(t)=\left\{\frac{M(t)}{N(t)}:M(t),N(t)\in\F_q[t],~ N(t)\neq0\right\}.\] A Laurent series $\Theta(t)\in\F_q(\!(t^{-1})\!)$ is called \textbf{irrational} if $\Theta(t)\not\in\F_q(t)$. Additionally, the counterpart to the prime numbers is given by the set of irreducible polynomials over $\F_q[t]$. \\

\noindent In analogy with the real numbers, every Laurent series $\Theta(t)\in \mathbb{F}_q(\!(t^{-1})\!)$ has a unique continued fraction expansion of the form
\begin{equation}
    \Theta(t)=A_0^{[\Theta]}(t)+\frac{1}{A_1^{[\Theta]}(t)+\frac{1}{A_2^{[\Theta]}(t)+\frac{1}{\ddots}}}=\left[A_0^{[\Theta]}(t);A_1^{[\Theta]}(t),A_2^{[\Theta]}(t),\dots\right],
\end{equation}
where $A_i^{[\Theta]}(t)\in \mathbb{F}_q[t]$ for $i\geq 0$ and $\deg\left(A_i^{[\Theta]}(t)\right)\geq 1$ for $i\geq 1$. Furthermore, an irrational Laurent series $\Theta(t)\in\F_q(\!(t^{-1})\!)$ is the solution to a quadratic equation with coefficients in $\mathbb{F}_q[t]$ if and only if its continued fraction expansion is eventually periodic \cite[Theorem 3.1]{DAFFsurvey}. That is, $\Theta(t)\in\F_q(\!(t^{-1})\!)$ is a quadratic irrational if and only if then there exist $m_{\Theta},\ell_{\Theta}\geq 1$, and there exist polynomials $B_i^{[\Theta]}(t)\in\F_q[t]$ for $0\le i \le m_{\Theta}\in\N$ and $A_i^{[\Theta]}(t)\in\F_q[t]$ for $0\le i \le \ell_{\Theta}\in\N$ such that
\begin{equation}
    \Theta(t)=\left[B_0^{[\Theta]}(t);B_1^{[\Theta]}(t),\dots,B_{m_{\Theta}}^{[\Theta]}(t),\overline{A_0^{[\Theta]}(t),\dots ,A_{\ell_{\Theta}}^{[\Theta]}(t)}\right].
\end{equation}

\noindent Let $P(t)\in \mathbb{F}_q[t]$ be an irreducible polynomial. Similarly to the real case, it is natural to study how the degree of partial quotients of $P(t)^k\cdot\Theta(t)$ evolve as $k\rightarrow \infty$. This was studied by de Mathan and Teulié in \cite{dMT}, who established (when $P(t)=t$) that the maximal degree of these periodic partial quotients is unbounded. 
\begin{theorem}{\cite[Theorem 4.5]{dMT}}
\label{thm:dMT}
For every quadratic irrational $\Theta(t)\in\F_q(\!(t^{-1})\!)$, one has
\begin{equation}
    \limsup_{k\rightarrow\infty}\max_{i=1,\dots,\ell_{t^k\Theta}}\deg\left(A_i^{[t^k\Theta]}(t)\right)=\infty.
\end{equation}
\end{theorem}

\noindent So far, the analogy between the real and function field case has remained consistent. However, this breaks down in view of Theorem \ref{thm:PS}, which is the function field analogue of Theorem \ref{thm:ASEsc}.

\begin{theorem}{\cite[Theorem 1]{PS}}
    \label{thm:PS}
    Let $P(t)\in \mathbb{F}_q[t]$ be an irreducible polynomial, and let $\Theta(t)\in \mathbb{F}_q(\!(t^{-1})\!)$ be a quadratic irrational. Then, there exists $0<c\leq 1$, depending only on $P(t)$ and $\Theta(t)$, such that 
    \begin{equation}
        \label{eqn:Esc1Excursion}
        \lim_{n\rightarrow\infty}\liminf_{k\rightarrow \infty}\frac{\max_{i=1,\dots,\ell_{\Theta\cdot P^k}}\max\left\{\deg\left(A_i^{[\Theta\cdot P^k]}(t)\right)-n,0\right\}}{\sum_{i=1}^{\ell_{\Theta\cdot P^k}}\deg\left(A_i^{[\Theta\cdot P^k]}(t)\right)}\geq c.
    \end{equation}
\end{theorem}
\noindent In fact, Theorem \ref{thm:PS} is saying more than just that Theorem \ref{thm:ASEsc} does not hold over function fields. The additional `minus $n$' term in the numerator of equation (\ref{eqn:Esc1Excursion}) is saying that Theorem \ref{thm:ASEsc} fails to adapt to function fields even when one is allowed to remove any fixed amount from the size of the maximal partial quotient (that is, from the numerator of (\ref{eqn:noEscR})). This lead to the following conjecture.

\begin{conjecture}{\cite[Conjecture 6]{KPS}}
    \label{conj:KPS}
    For every irreducible polynomial $P(t)\in \mathbb{F}_q[t]$ and every quadratic irrational Laurent series $\Theta(t)\in \mathbb{F}_q(\!(t^{-1})\!)$, one has \begin{equation}\label{eqn:MassEqn}
    \lim_{n\rightarrow \infty}\liminf_{k\rightarrow \infty}\frac{\max_{i=1,\dots,\ell_{\Theta\cdot P^k}}\max\left\{\deg\left(A_{i}^{[\Theta\cdot P^k]}(t)\right)-n,0\right\}}{\sum_{i=1}^{\ell_{\Theta\cdot P^k}}\deg\left(A_{i}^{[\Theta\cdot P^k]}(t)\right)}=1.
    \end{equation}
    
\end{conjecture}

\noindent In 2025, the first named author, Erez Nesharim, and Uri Shapira \cite[Theorem 1.2]{AN} disproved Conjecture \ref{conj:KPS} when $P(t)=t$ and $q=2$.
\begin{theorem}{\cite[Theorem 1.2]{AN}}\label{thm: TM_partial}
    Let $\{\tau_n\}_{n\in \mathbb{N}}\subseteq \mathbb{F}_2$ be the Thue-Morse sequence, and let $\tau=\sum_{n\in \mathbb{N}}\tau_nt^{-n}\in \mathbb{F}_2(\!(t^{-1})\!)$. Then, \begin{equation}\label{eqn:TMMassEqn}
    \lim_{n\rightarrow \infty}\liminf_{k\rightarrow \infty}\frac{\sum_{i=1}^{\ell_{\tau\cdot t^k}}\max\left\{\deg\left(A_{i}^{[\tau\cdot t^k]}(t)\right)-n,0\right\}}{\sum_{i=1}^{\ell_{\tau\cdot t^k}}\deg\left(A_{i}^{[\tau\cdot t^k]}(t)\right)}=\frac{2}{3}\cdotp
    \end{equation}
 \end{theorem}
\noindent Above, the sum in the numerator of equation (\ref{eqn:TMMassEqn}) replaces the maximum in the numerator of equation (\ref{eqn:MassEqn}). Clearly, this implies that Theorem \ref{thm: TM_partial} is a counterexample to Conjecture \ref{conj:KPS}.
\subsection{Escape of Mass in Hecke Trees over $\mathbb{F}_q(\!(t^{-1})\!)$}\label{Sect: Hecke}\hfill\\
\noindent The contents of subsection \ref{Sect: Hecke} provide an alternate description of the results in subsection \ref{Sect: EoM_intro} in the language of so called \textit{Hecke trees}. This information is provided for completeness, but it is not necessary to understand the paper as a whole.\\ 

\noindent By adapting the analysis in \cite{EW,P,KPS}, both \cite{PS} and \cite{AS} used results on distributions of measures in \textit{Hecke trees}\footnote{For more information about Hecke trees see \cite{Ser,P}.} to obtain results regarding the continued fraction expansions of quadratic irrationals. Indeed, many of the theorems in the previous subsection can be rephrased in terms of Hecke trees. \\

\noindent One associates a quadratic irrational $\Theta(t)\in \mathbb{F}_q(\!(t^{-1})\!)$ to a lattice $$\Lambda\in \mathcal{L}_2:=\operatorname{PGL}_2\left(\mathbb{F}_q(\!(t^{-1})\!)\right)/\operatorname{PGL}_2\left(\mathbb{F}_q[t]\right)$$ by using the geometric embedding. Then, by Dirichlet's units theorem, this lattice has a compact orbit under the diagonal group $$A=\Bigg\{\begin{pmatrix}
    a_1&\\
    &a_2
\end{pmatrix}:a_1,a_2\neq 0\Bigg\}.$$ Hence, $\Lambda$ supports a unique $A$ invariant probability measure $\mu_{A\Lambda}$.\\

\noindent One of the fundamental questions in homogeneous dynamics revolves around finding accumulation points of the sequence of measures $\{\mu_{A\Lambda_n}\}_{n\in\mathbb{N}}$, where $\{A\Lambda_n\}_{n\in \mathbb{N}}$ is some sequence of compact $A$ orbits. Due to Dirichlet's units theorem, this question is a way to study the distribution of number fields. This question has also been studied in both the real \cite{ELMV,Sha,AS} and the positive characteristic setting \cite{KPS,AEsc,AN,DPS}.
\begin{definition}\label{def:Hecke_EoM}
    Let $0<c\leq 1$. A sequence of compact $A$ orbits $\{A\Lambda_n\}_{n\in \mathbb{N}}$ exhibits \textbf{$c$ escape of mass} if for every accumulation point $\mu=\lim_{k\rightarrow \infty}\mu_{A\Lambda_{n_k}}$, one has $\mu(\mathcal{L}_2)\leq 1-c$. The escape of mass is said to be \textbf{full} if one can take $c=1$.  
\end{definition}

\noindent In \cite{KPS}, Kemarski, Paulin, and Shapira proved that all compact $A$ orbits arising from a sequence of number fields (which satisfy certain arithmetic relationships) exhibit full escape of mass. They proved this result by studying dynamics of the diagonal group acting on Hecke trees.
\begin{definition}[Hecke Tree]
    Let $P(t)\in \mathbb{F}_q[t]$ be an irreducible polynomial. Lattices $\Lambda_1,\Lambda_2\in \mathcal{L}_2$ are called \textbf{$P(t)$-Hecke friends} if 
    \begin{enumerate}
        \item $\Lambda_2\subseteq \Lambda_1$, and
        \item there exists some $n\geq 1$ such that $\Lambda_1/\Lambda_2\cong \mathbb{F}_q[t]/P(t)^n\mathbb{F}_q[t]$.
    \end{enumerate}
    If $n=1$, $\Lambda_1$ and $\Lambda_2$ are called \textbf{$P(t)$-Hecke neighbors}. 
\end{definition}
\noindent Let $\Lambda_0\in \mathcal{L}_2$. Define $\mathbb{T}_{P}(\Lambda_0)$ as the tree whose nodes are $P(t)$-Hecke friends of $\Lambda_0$, where two lattices are connected by an edge if they are Hecke neighbors. 
\begin{definition}[Branch]
    A branch in $\mathbb{T}_{P}(\Lambda_0)$ is a sequence $\{\Lambda_n\}_{n\in \mathbb{N}}$ such that for every $n\in \mathbb{N}$, $\Lambda_{n+1}$ and $\Lambda_n$ are Hecke neighbors. A branch is identified with a point on the boundary of the tree $\partial \mathbb{T}_{P}(\Lambda_0)\cong \mathbb{P}(\mathbb{F}_q(t)_{P})$, where $\mathbb{F}_q(t)_{P}$ is the completion of $\mathbb{F}_q(t)$ with respect to the $P$-adic norm (see \cite{P,KPS} for details), and $\mathbb{P}(\mathbb{K})$ is the projective line defined over a field $\mathbb{K}$. A branch is called rational if it is identified to a point in $\mathbb{P}(\mathbb{F}_q(t))$.
\end{definition}
\begin{theorem}{\cite[Theorem 1]{KPS}}
    \label{thm:KPS}
    Let $\Lambda_0\in \mathcal{L}_2$, and let be a $P(t)\in \mathbb{F}_q[t]$ be an irreducible polynomial. Then, there exists some $0<c=c_{\Lambda_0,P}\leq 1$, such that for every rational branch in $\mathbb{T}_P(\Lambda_0)$ exhibits $c$ escape of mass.
\end{theorem}
\noindent By the analysis in \cite{P,PS}, Theorem \ref{thm:PS} and Theorem \ref{thm:KPS} are equivalent. Furthermore, Conjecture \ref{conj:KPS} can be rephrased in terms of Hecke Trees: 
\begin{conjecture}{\cite[Conjecture 6]{KPS}}
\label{conj:FullEscMass}
For every choice of rational branch $\{\Lambda_n\}_{n\in \mathbb{N}}$, one has that $\lim_{n\rightarrow \infty}\mu_{A\Lambda_n}=0$. 
\end{conjecture}
\noindent Theorem \ref{thm: TM_partial} is also rephrased in the language of Hecke trees.
\begin{theorem}{\cite[Theorem 1.2]{AN}}\label{thm: TM_partial_Hecke}
    For every rational branch $\{\Lambda_n\}_{n\in \mathbb{N}}$ in $\mathbb{T}_t(\Lambda_{\tau})$, there exists a subsequence $\{n_k\}_{k\in \mathbb{N}}\subseteq \mathbb{N}$, such that the limiting measure $\mu:=\lim_{k\rightarrow \infty}\mu_{A\Lambda_{n_k}}$ satisfies $\mu(\mathcal{L}_2)\geq \frac{1}{3}$.
 \end{theorem}
\subsection{Main Results}\hfill\\
\noindent This paper is intended to be seen as the second in a pair of papers whose collective goal is to disprove Conjecture \ref{conj:FullEscMass}. Indeed, the main results of this paper are attained by closely examining the so-called \textit{number wall of the} $p$-\textit{Cantor sequence} (see sections \ref{sec:NumberWalls} and \ref{sec:pCantor} for definitions), which is a two-dimensional sequence over $\F_p$. The global structure of this sequence was established in \cite{AR} and the theorems contained therein are crucial to this work. \\

\noindent The results in this paper fall into three categories. 
\subsubsection{\textbf{Disproving Conjecture }\ref{conj:KPS}}\hfill\\
 \noindent The primary aim of this paper is to fill the gap in Conjecture \ref{conj:KPS} left by Theorem \ref{thm: TM_partial}. That is, to disprove Conjecture \ref{conj:KPS} in all cases. \\

\noindent To this end, the following definition is made. 
\begin{definition}
\label{def: QuadEscMass}
Let $0<c\leq 1$ and let $P(t)\in \mathbb{F}_q[t]$ be an irreducible polynomial. A quadratic irrational $\Theta(t)\in\F_q(\!(t^{-1})\!)$ exhibits $c$ \textbf{escape of mass} with respect to the sequence $\left\{P(t)^k\right\}_{k\in \mathbb{N}}$ if 
    \begin{equation}\label{eqn:MassEqn}
    \lim_{n\rightarrow \infty}\liminf_{k\rightarrow \infty}\frac{\sum_{i=1}^{\ell_{\Theta\cdot P^k}}\max\left\{\deg\left(A_{i}^{[\Theta\cdot P^k]}(t)\right)-n,0\right\}}{\sum_{i=1}^{\ell_{\Theta\cdot P^k}}\deg\left(A_{i}^{[\Theta\cdot P^k]}(t)\right)}\geq c.
    \end{equation}
    \noindent When $c=1$, $\Theta(t)$ is said to exhibit \textbf{full escape of mass}. If $c$ is maximal, then, $\Theta(t)$ exhibits \textbf{exactly} $c$ \textbf{escape of mass}. 
\end{definition}

\begin{remark}
\label{rem:EquivEscMass}
By the analysis in \cite{P,PS}, Definitions \ref{def: QuadEscMass} and \ref{def:Hecke_EoM} are equivalent. 
\end{remark}

\noindent Just as in Theorem \ref{thm: TM_partial}, if $\Theta(t)\in\F_p(\!(t^{-1})\!)$ exhibits exactly $c<1$ escape of mass, then $\Theta(t)$ is a counterexample to Conjecture \ref{conj:KPS}.\\

 \noindent The disproof of Conjecture \ref{conj:KPS} comes in the form of two theorems. The former such result states that, for \textit{any} odd prime $p$, the so-called $p$-\textit{Cantor sequence} $\left\{c_i^{(p)}\right\}_{i\ge0}$ is an explicit counterexample to Conjecture \ref{conj:FullEscMass} over $\F_p$ in the case $P(t)=t$. The definition of this sequence is postponed to Section 3, as it is not required to state the results of this paper. For convenience, define\footnote{The extra factor of $t^{-1}$ is here as it is often useful to index sequences starting from the zeroth term, but it is also useful to consider Laurent series in $\F_p(\!(t^{-1})\!)$ to have no fractional part. } $\Theta_p(t)=t^{-1}\cdot \sum_{i=0}^\infty c^{(p)}_i t^{-i}\in\F_p(\!(t^{-1})\!)$.  

\begin{theorem}
\label{cor:minEsc}
    For any odd prime $p$, the Laurent series $\Theta_p(t)\in\F_p(\!(t^{-1})\!)$ exhibits exactly $\frac{2}{p}$ escape of mass with respect to the sequence $\{t^k\}_{k\ge0}$. 
\end{theorem}

\noindent In other words, Theorem \ref{cor:minEsc} shows that for every increasing sequence of natural numbers $\{k_i\}_{i\in\N}$, one has that
\begin{equation}
\label{eqn:2/pEsc}
\lim_{n\rightarrow\infty}\lim_{i\rightarrow \infty}\frac{\sum_{i=1}^{\ell_{\Theta_p\cdot t^{k_i}}}\max\left\{\deg\left(A_{i}^{[\Theta_p\cdot t^{k_i}]}(t)\right)-n,0\right\}}{\sum_{i=1}^{\ell_{\Theta_p\cdot t^{k_i}}}\deg\left(A_{i}^{[\Theta_p\cdot t^{k_i}]}(t)\right)}\geq \frac{2}{p}.\end{equation}

\noindent The latter result is a transference principle that allows one to take any Laurent series with $c$-escape of mass with respect to the sequence $\{t^k\}_{k\ge0}$ and return a Laurent series that has $c$-escape of mass with respect to the sequence $\{P(t)^k\}_{k\ge0}$, for \textit{any} choice of irreducible polynomial $P(t)\in\F_p[t]$. To this end, let $\Theta(t)=\sum_{i=-h}^\infty a_i t^{-i}\in\F_p(\!(t^{-1})\!)$ be a Laurent series and $P(t)\in\F_q[t]$ be some polynomial. Then, define the Laurent series \begin{equation}
     \label{eqn: t_to_P(t)} \Theta(P(t)):=\sum_{i=-h}^\infty a_i P(t)^{-i} \in \F_p(\!(t^{-1})\!).
 \end{equation} 

\begin{theorem}\label{thm: t_to_P(t)}
    Let $c\in[0,1]$ and let $\Theta(t)\in\F_p(\!(t^{-1})\!)$ be a Laurent series that exhibits $c$-escape of mass with respect to the sequence $\{t^k\}_{k\ge0}$. Then, for any irreducible polynomial $P(t)\in\F_p[t]$, the Laurent series $\Theta(P(t))\in\F_p(\!(t^{-1})\!)$ exhibits $c$-escape of mass with respect to the sequence $\{P(t)^k\}$.
\end{theorem}
\noindent By applying Theorem \ref{thm: t_to_P(t)} to Theorems \ref{thm: TM_partial} and \ref{cor:minEsc}, the following corollary is immediate.
\begin{corollary}\label{cor: esc_of_mass}
    Recall $\tau(t)$ and $\Theta_p(t)$ from Theorems \ref{thm: TM_partial} and \ref{cor:minEsc}, respectively. For any irreducible polynomial $P(t)\in\F_p[t]$, $\tau(P(t))$ ($\Theta_p(P(t))$, respectively) exhibits exactly $\frac{2}{3}$ $\left(\frac{2}{p}\text{, respectively}\right)$ escape of mass with respect to the sequence $\{P(t)^k\}_{k\ge0}$.
\end{corollary}
\noindent Corollary \ref{cor: esc_of_mass} amounts to a complete disproof of Conjecture \ref{conj:KPS}.\\

\noindent In terms of Hecke trees, Theorems \ref{cor:minEsc} and \ref{thm: t_to_P(t)} imply that Conjecture \ref{conj:FullEscMass} does not hold in any odd characteristic. In particular, due to the analysis of \cite{P,PS}, one obtains the following corollary to Theorem \ref{cor:minEsc}. 
\begin{corollary}
    For every odd $p$, every irreducible polynomial $P(t)\in\F_p[t]$ and for any rational branch $\{\Lambda_n\}_{n\in \mathbb{N}}\subseteq \mathbb{T}_P\left(\Lambda_{\Theta_p}\right)$, there exists a subsequence $\{n_k\}$, such that the limiting measure $\mu:=\lim_{k\rightarrow \infty}\mu_{A\Lambda_{n_k}}$ satisfies $\mu(\mathcal{L}_2)\geq\frac{p-2}{p}$.
\end{corollary}
\subsubsection{\textbf{Maximal Escape of Mass}}\hfill\\
\noindent In view of Corollary \ref{cor: esc_of_mass}, one considers an alternative version of Definition \ref{def: QuadEscMass} where $\limsup$ is used in place of $\liminf$ in equation (\ref{eqn:MassEqn}):
\begin{definition}
    Let $P(t)\in\F_p[t]$ be an irreducible polynomial. A quadratic irrational $\Theta(t)\in \mathbb{F}_q(\!(t^{-1})\!)$ exhibits \textbf{$c$-maximal escape of mass} with respect to the sequence $\{P(t)^k\}_{k\ge0}$ if 
    \begin{equation}
        \lim_{n\rightarrow\infty}\limsup_{k\rightarrow \infty}\frac{\sum_{i=1}^{\ell_{\Theta\cdot P^k}}\max\left\{\deg\left(A_{i}^{[\Theta\cdot P^k]}(t)\right)-n,0\right\}}{\sum_{i=1}^{\ell_{\Theta\cdot P^k}}\deg\left(A_{i}^{[\Theta\cdot P^k]}(t)\right)}\geq c.
    \end{equation}
    When $c=1$, $\Theta(t)$ is said to exhibit \textbf{full maximal escape of mass} with respect to the sequence $\{P(t)^k\}_{k\ge0}$.
\end{definition}
\noindent Similar to the Thue Morse sequence \cite[Corollary 4.8]{AN}, the $p$-Cantor sequence has full maximal escape of mass. 
\begin{theorem}
\label{thm:MaxEsc1}
    For any odd prime $p$ and any irreducible polynomial $P(t)\in\F_p[t]$, the Laurent series $\Phi(t):=\Theta_p(P(t))\in\F_p$ exhibits full maximal escape of mass with respect to the sequence $\{P(t)^k\}_{k\ge0}$. 
\end{theorem}

\noindent Both Theorem \ref{thm:MaxEsc1} and \cite[Corollary 4.8]{AN} provide evidence to the following Conjecture, which refines Conjecture \ref{conj:FullEscMass}.
\begin{conjecture}
    Let $p$ be any prime, $P(t)\in\F_p[t]$ be any irreducible polynomial and $\Theta(t)\in \mathbb{F}_q(\!(t^{-1})\!)$ be any quadratic irrational. Then, $\Theta(t)$ exhibits full maximal escape of mass with respect to the sequence $\{P(t)^k\}_{k\ge0}$. 
\end{conjecture}
\subsubsection{\textbf{Generic Escape of Mass}}\hfill\\
\noindent If $\Theta_p(P(t))$ exhibits exactly $\frac{2}{p}$ escape of mass with respect to the sequence $\{P(t)^k\}_{k\ge0}$, then this means there exists a \textit{specific} increasing sequence of natural numbers $\{k_i\}_{i\in \mathbb{N}}$, for which the inequality in equation \eqref{eqn:2/pEsc} is an equality. Similarly, Theorem \ref{thm:MaxEsc1} states that there exists a different increasing sequence of natural numbers $\{k'_i\}_{i\in\N}$ such that, when $k_i$ replaced by $k'_i$, the left hand side of \eqref{eqn:2/pEsc} is equal to 1. Therefore, it is natural to inquire what one can say about the escape of mass when $\{k_i\}_{i\in\N}$ is replaced by a "\textit{generic}" increasing sequence of integers $k\in\N$. That is, to ask what one can say about the escape of mass from a metric point of view.\\

\noindent Since there is no probability measure supported on all of $\mathbb{N}$, a natural notion of probability on $\mathbb{N}$ is given by the \textit{density}: 
\begin{definition}
    Let $A\subseteq \mathbb{N}$ and let $c\in [0,1]$. Then, $A$ has \textbf{density} $c$ if 
    $$d(A):=\lim_{m\rightarrow \infty}\frac{\vert A\cap \{1,2,\dots,m\}\vert}{m}=c.$$
\end{definition}
\noindent To apply the concept of density to escape of mass problems, the following quantity is defined. 
\begin{definition}\label{def:DiagEscMass}For every $n>1$ and for $k\in \mathbb{N}$, define the $n$-\textbf{escape of mass in the} $k^\textbf{th}$ \textbf{diagonal} of the quadratic irrational $\Theta(t)$ with respect to the sequence $\{P(t)^k\}_{k\ge0}$ as
\begin{equation}
\label{eqn:EscDefWall}
    e_{k,n}:=e_{k,n}(\Theta(t)):=\frac{\sum_{i=1}^{\ell_{\Theta\cdot P^k}}\max\left\{\deg\left(A_i^{[\Theta\cdot P^k]}(t)\right)-n,0\right\}}{\sum_{i=1}^{\ell_{\Theta\cdot P^k}}\deg\left(A_i^{[\Theta\cdot P^k]}(t)\right)}.
\end{equation}
\end{definition}
\begin{remark}
    In this paper, generic escape of mass is only considered in the case $P(t)=t$. Hence, $P(t)$ is not reflected in the above notation.
\end{remark}
\begin{remark}
    The use of the word "diagonal" in Definition \ref{def:DiagEscMass} is explained in Chapter 3. 
\end{remark}
\noindent Clearly, the escape of mass as defined in Definition \ref{def: QuadEscMass} is given by $$\lim_{n\to\infty}\liminf_{k\rightarrow \infty} e_{k,n}(\Theta(t)).$$
\noindent The metric generalisation of escape of mass is now defined.
\begin{definition}
    A quadratic irrational Laurent series $\Theta(t)\in\F_q(\!(t^{-1})\!)$ is said to exhibit \textbf{full generic escape of mass} with respect to the sequence $\{P(t)^k\}_{k\ge0}$ if for every $\varepsilon>0$,
    \begin{equation}
    \label{eqn:genEscFull}
        \lim_{n\rightarrow \infty}d(\{k\in \mathbb{N}:e_{k,n}(\Theta(t))>1-\varepsilon\})=1.
    \end{equation}
\end{definition}
\begin{theorem}
\label{thm:GenEsc=1}
    For every odd prime $p$, the $p$-Cantor sequence $\Theta(t)$ exhibits full generic escape of mass with respect to the sequence $\{t^k\}_{k\ge0}$. 
\end{theorem}
\noindent The same result is established when the $p$-Cantor sequence is replaced with the Thue-Morse sequence.
\begin{theorem}
    \label{thm:TMFullEsc}
    Let $\{\tau_n\}_{n\in \mathbb{N}}$ be the Thue-Morse sequence. Then, $\tau(t)=\sum_{n=1}^{\infty}\tau_nt^{-n}$ exhibits generic full escape of mass with respect to the sequence $\{t^k\}_{k\ge0}$. 
\end{theorem}
\noindent By the analysis of \cite{P,PS}, one obtains the following corollary to Theorem \ref{thm:GenEsc=1} and Theorem \ref{thm:TMFullEsc}.
\begin{corollary}
    For every prime $p$, there exists a quadratic irrational $\nu_p\in \mathbb{F}_p(\!(t^{-1})\!)$, such that every $\varepsilon>0$, and for every rational branch $\{\Lambda_n\}_{n\in \mathbb{N}}\subseteq \mathbb{T}_t\left(\Lambda_{\nu_p}\right)$, one has
    $$\lim_{n\rightarrow \infty}d\left(\{k\in \mathbb{N}:\mu_{A\Lambda_k}(\mathcal{L}_d)>\varepsilon\}\right)=0.$$
\end{corollary}
\noindent Due to Conjecture \ref{conj:FullEscMass}, \cite[Theorem 2]{KPS}, Theorem \ref{thm:GenEsc=1}, Theorem \ref{thm:TMFullEsc}, and \cite[Corollary 4.8]{AN}, the following conjecture is made.
\begin{conjecture}
\label{conj:genEsc=1}
    Let $\Theta(t)\in \mathbb{F}_q(\!(t^{-1})\!)$ be a quadratic irrational and let $P(t)\in\F_q[t]$ be an irreducible polynomial. Then, $\Theta(t)$ exhibits full generic escape of mass with respect to the sequence $\{P(t)^k\}_{k\ge0}$. 
\end{conjecture}
\subsubsection{\textbf{Structure of Paper}}\hfill\\
\noindent The proof of Theorem \ref{thm: t_to_P(t)} is elementary and completed in Section \ref{Sect: 1.5}. Theorems \ref{cor:minEsc}, \ref{thm:MaxEsc1}, \ref{thm:GenEsc=1}, and \ref{thm:TMFullEsc} are achieved by rephrasing Definition \ref{def: QuadEscMass} in terms of the properties of the so-called \textit{number wall} of a sequence, which is defined in \ref{sec:NumberWalls}. Section \ref{Sect: 3} explicitly states the connection between escape of mass and number walls, and Section \ref{sec:pCantor} then introduces the $p$-Cantor sequence and states the construction of its number wall. Finally, \ref{Sect: 7,0} completes the proof of Theorems \ref{cor:minEsc}, \ref{thm:MaxEsc1}, \ref{thm:GenEsc=1}, and \ref{thm:TMFullEsc}.  

\subsection*{Acknowledgments}
The first named author would like to thank Uri Shapira for proposing this question to her. She would also like to thank Erez Nesharim for multiple discussions which led to the ideas in this paper, as well as Omri Solan, Yeor Hafouta, and Jon Chaika for useful intuitions about generic escape of mass. The second named author would like to thank Fred Lunnon for his discussions relating to the number wall of the $p$-Cantor sequence. Both authors would also like to thank Faustin Adiceam for suggesting that they work together on this question. 
\section{Establishing Theorem \ref{thm: t_to_P(t)}}\label{Sect: 1.5}
\noindent Two lemmas are required to prove Theorem \ref{thm: t_to_P(t)}. The first relates the coefficients of $\Theta(t)^{-1}$ to $\Theta(P(t))^{-1}$.
\begin{lemma}\label{lem: t_P(t)_inv}
    Let $\Theta(t)=\sum_{i=-h}^\infty a_it^{-i}\in\F_p(\!(t^{-1})\!)$ and let $\Theta(t)^{-1}=\sum_{i=h}^{\infty} b_it^{-i}$. Then, for any polynomial $P(t)\in\F_p[t]$, $\Theta(P(t))^{-1}=\sum_{i=h}^\infty b_i P(t)^{-i}$.
\end{lemma}

\begin{proof}
    By multiplying $\Theta(t)$ and $\Theta(t)^{-1}$, one obtains \[a_{-h}b_h=1,~~~\text{ and }~~~\sum_{j=0}^i a_{-h+j}b_{h+i-j}=0~~~\text{ for all } i>0.\] To conclude, note that \begin{align*}
        \Theta(P(t))\cdot \sum_{i=h}^\infty b_i P(t)^{-i} &=\sum_{i=-h}^\infty a_i P(t)^{-i}\cdot \sum_{i=h}^\infty b_i P(t)^{-i}\\ &= \sum_{i=0}^\infty P(t)^{-i}\cdot\sum_{j=0}^i a_{-h+j}b_{h+i-j}\\
        &=1 \cdot P(t)^0 +0.
    \end{align*}
\end{proof}
\noindent The second lemma relates the continued fraction coefficients of $\Theta(t)\in\F_p(\!(t^{-1})\!)$ to the partial quotients of $\Theta(P(t))$. 
\begin{lemma}\label{lem: t_P(t)_CF}
    Let $P(t)\in\F_p[t]$ be a polynomial and let $\Theta(t)\in\F_p(\!(t^{-1})\!)$ be a Laurent series whose continued fraction expansion is given by \[\Theta(t)=[A_0(t); A_1(t),A_2(t),\dots].\] Furthermore, let
    \[\Theta(P(t))=[B_0(t); B_1(t),B_2(t),\dots]\]
    be the continued fraction expansion of $\Theta(P(t))$. Then, for every  $i\ge0$, one has $B_i(t)=A_i(P(t))$.
\end{lemma}
\begin{proof}
    \noindent The proof proceeds by induction on $i$. For the base case, note that $A_0(t)$ ($B_0(t)$, respectively) is equal to the polynomial part of $\Theta(t)$ ($\Theta(P(t))$, respectively). Therefore, from the definition of $\Theta(P(t))$ (equation (\ref{eqn: t_to_P(t)})), $B_0(t)=A_0(P(t))$.\\

    \noindent Assume that for all $i'<i$, $B_{i'}(t)=A_{i'}(P(t))$. That is, \[\Theta(P(t))=[A_0(P(t)); A_1(P(t)), \dots, A_{i-1}(P(t)), B_i(t),B_{i+1}(t),\dots].\] For $0\le j\le i$, define $$C_j(t)=[A_j(t); A_{j+1}(t),\dots] ~~~ \text{ and } ~~~ D_j(t)=[A_j(P(t));, A_{j+1}(P(t)),\dots A_{i-1}(P(t)), B_i(t),B_{i+1}(t),\dots].$$ Note that, $C_j(t)$ and $D_j(t)$ satisfy the relations \[C_j(t)=\frac{1}{C_{j-1}(t)-A_{j-1}(t)}~~~\text{ and }~~~ D_j(t)=\frac{1}{D_{j-1}(t)-A_{j-1}(P(t))}\cdotp\] It is now shown by induction that $D_j(t)=C_j(P(t))$ for all $0\le j\le i$. The base case is immediate from the fact that $D_0(t)=\Theta(P(t))$ and $C_0(t)=\Theta(t)$. For the induction step, \begin{align*}
        D_j(t)&=\frac{1}{C_{j-1}(P(t))-A_{j-1}(P(t))}~\substack{\text{Lemma }\ref{lem: t_P(t)_inv}\\=}~ C_j(P(t)).
    \end{align*}
    Finally, note that $B_i(t)$ is the polynomial part of $D_i(t)=C_i(P(t))$, which is equal to $A_i(P(t))$.
\end{proof}

\noindent The proof of Theorem \ref{thm: t_to_P(t)} is now completed.
\begin{proof}[Proof of Theorem \ref{thm: t_to_P(t)}]
    Let $P(t)\in\F_p[t]$ be any irreducible polynomial. For notational ease, define $\Phi(t):=\Theta(P(t))$ and $d_P:=\deg(P(t))$. \\
    
    \noindent With notation from Definition \ref{def: QuadEscMass}, $$A_{i}^{[\Phi\cdot P^k]}(t)~\substack{Lemma~ \ref{lem: t_P(t)_CF}\\=}~A_i^{[\Theta\cdot t^k]}(P(t)),$$ and hence\[\deg\left(A_{i}^{[\Phi\cdot P^k]}(t)\right)=d_P\cdot \deg\left(A_{i}^{[\Theta\cdot t^k]}(t)\right).\]Also by Lemma \ref{lem: t_P(t)_CF}, $\ell_{\Phi\cdot P^k}=\ell_{\Theta\cdot t^k}$. Therefore, \begin{equation*}
        \frac{\sum_{i=1}^{\ell_{\Phi\cdot P^k}}\max\left\{\deg\left(A_{i}^{[\Phi\cdot P^k]}(t)\right)-n,0\right\}}{\sum_{i=1}^{\ell_{\Phi\cdot P^k}}\deg\left(A_{i}^{[\Phi\cdot P^k]}(t)\right)}=\frac{\sum_{i=1}^{\ell_{\Theta\cdot t^k}}\max\left\{\deg\left(A_{i}^{[\Theta\cdot t^k]}(t)\right)-\frac{n}{d_P},0\right\}}{\sum_{i=1}^{\ell_{\Theta\cdot t^k}}\deg\left(A_{i}^{[\Theta\cdot t^k]}(t)\right)}\cdotp
    \end{equation*}
    \noindent Therefore, when $n=d_p\cdot n'$ for any $n'\in\N$, \begin{equation}
        \label{eqn:limP(t)vst}
        \liminf_{k\to\infty}\frac{\sum_{i=1}^{\ell_{\Phi\cdot P^k}}\max\left\{\deg\left(A_{i}^{[\Phi\cdot P^k]}(t)\right)-n,0\right\}}{\sum_{i=1}^{\ell_{\Phi\cdot P^k}}\deg\left(A_{i}^{[\Phi\cdot P^k]}(t)\right)}=\liminf_{k\to\infty}\frac{\sum_{i=1}^{\ell_{\Theta\cdot t^k}}\max\left\{\deg\left(A_{i}^{[\Theta\cdot t^k]}(t)\right)-n',0\right\}}{\sum_{i=1}^{\ell_{\Theta\cdot t^k}}\deg\left(A_{i}^{[\Theta\cdot t^k]}(t)\right)}\cdotp
    \end{equation}
    \noindent To conclude the proof, note that as $n,n'\to\infty$, the above equations both form bounded monotonic decreasing sequences that intersect infinitely often. Hence, the limits as $n'\to\infty$ of the left hand side and the right hand side of \eqref{eqn:limP(t)vst} are equal to one another.
\end{proof}
\section{Number Walls}
\label{sec:NumberWalls}
\noindent All remaining proofs in this paper rely on carefully analyzing the number wall of the $p$-Cantor sequence as constructed in \cite{AR}. In this section, number walls are defined. 
\noindent The following family of matrices provide the bricks with which number walls are built.

\begin{definition} \label{Toe}
A matrix $(s_{i,j})$ for $0\le i\le n, 0\le j \le m$ is called \textbf{Toeplitz} if all the entries on any of its diagonals are equal to one another. Equivalently, $s_{i,j}=s_{i+1,j+1}$ for any $i,j\in\mathbb{N}$ for which these entries are defined. 
\end{definition}
\noindent Given a doubly infinite sequence $\mathbf{S}:= (s_i)_{i\in\mathbb{Z}}$, natural numbers $m$ and $v$, and an integer $n$, define the $(m+1)\times (v+1)$ \textbf{Toeplitz matrix} $T_S(n;m, v):= (s_{i-j+n})_{0\le i \le m, 0\le j \le v}$ as\\
\begin{align*}T_S(n;m,v):=\begin{pmatrix}
s_n & s_{n+1} & \dots & s_{n+v}\\
s_{n-1} & s_{n} & \dots & s_{n+v-1}\\
\vdots &&& \vdots\\
s_{n-m} & s_{n-m+1} & \dots & s_{n-m+v} 
\end{pmatrix}.\end{align*} If $v=m$, this is abbreviated to $T_S(n;m)$. 
\subsection{Definition of a Number Wall}\hfill
\begin{definition}\label{nw}
Let $\mathbf{S}=(s_i)_{i\in\mathbb{Z}}$ be a doubly infinite sequence over a field $\mathbb{K}$. The \textbf{number wall} of the sequence $\mathbf{S}$ is defined as the two dimensional array of numbers $W_\K(\mathbf{S})=(W_{\K}(\mathbf{S})[m,n])_{n,m\in\mathbb{Z}}$ with \begin{equation*}
    W_{\K}(\mathbf{S})[m,n]=\begin{cases}\det(T_S(n;m)) &\textup{ if } m\ge0,\\
    1 & \textup{ if } m=-1,\\
    0 & \textup{ if } m<-1. \end{cases}
\end{equation*}  
When $\K=\F_q$ for some prime power $q\in\N$, $W_\K(\mathbf{S})$ is abbreviated to $W_q(\mathbf{S})$.
\end{definition} 
\noindent To stay consistent with standard matrix notation, $m$ increases as the rows go down the page and $n$ increases from left to right. 

\begin{remark}
     Whilst the first instance of number walls in the literature is due to Conway and Guy \cite[pp. 85--89]{CG} , Toeplitz matrices (and their determinants) have been studied since the 1900s. The theory of number walls was greatly developed by Lunnon \cite{L01} and \cite{L09}, and they have since been applied to problems in Diophantine approximation \cite{ANL,R23,RG}, dynamics \cite{AN}, and discrepancy theory \cite{RobVdC}.

\end{remark}
\subsection{The Profile of a Number Wall}\hfill\\
\noindent To compute escape of mass, it is only important that an entry in the number wall is zero or nonzero (see Corollary \ref{cor:EscWall}), with its exact value in the latter case being irrelevant, whence this definition:
\begin{definition}\label{profile_def}
    Let $W_\K(\textbf{S})$ be the number wall of a sequence $\textbf{S}$ over a field $\K$. The \textbf{profile} of $W_\K(\textbf{S})$, denoted by 
    $$\chi(W_\K(\mathbf{S}))=\chi(W_\K(\mathbf{S}))[m,n]_{m\in\N,~ m\le n \le r-m},$$ is a two dimensional array with the same width and height as $W_\K(\textbf{S})$, is defined by  \[\chi(W_\K(\mathbf{S}))[m,n]=\begin{cases}
        0&\text{ if }W(\mathbf{S})[m,n]=0\\ X&\text{ otherwise.}
    \end{cases}\] Above, $X$ should be interpreted as an arbitrarily chosen symbol. In other words, the profile of a number wall shows only where its zero entries are located.
\begin{figure}[H]
    \centering
    \includegraphics[width=0.75\linewidth]{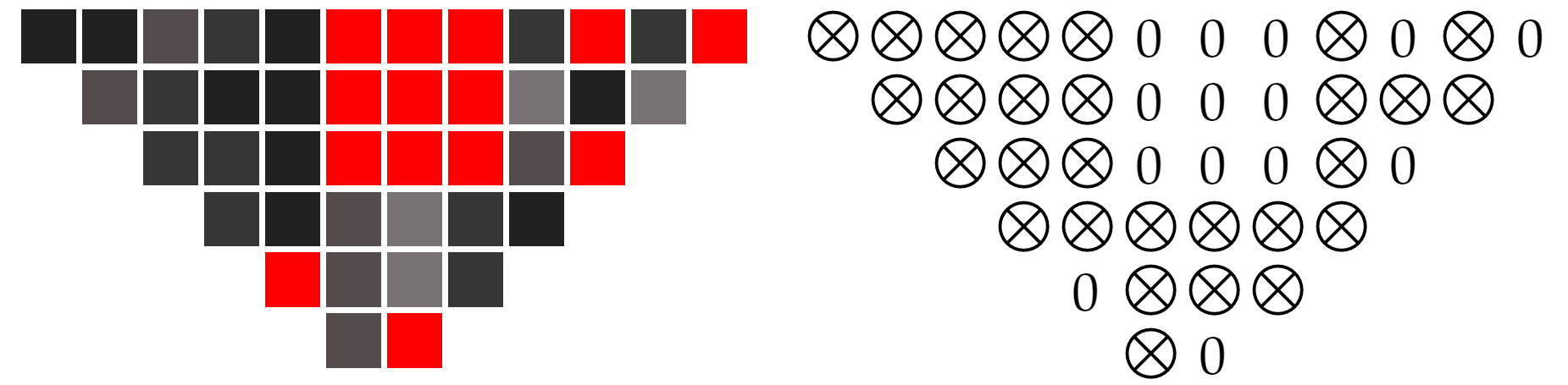}
    \caption{\textbf{Left}: The number wall $W_5(\textbf{S})$ for the sequence $\textbf{S}=(1,1,3,2,1,0,0,0,2,0,2,0)\in\F_5^{12}$, where the zeroes are coloured in red. \textbf{Right}: The profile of $W_\K(\textbf{S})$.}
\end{figure}
\end{definition}

\subsection{Diagonals of Number Walls}\label{subsub:DiagsNW}\hfill\\
\noindent For the purposes of studying the escape of mass of a Laurent series (see Section \ref{Sect: 3}), it is more natural to consider the entries that appear on a given diagonal. This motivates the following definition.
\begin{definition}
    \label{diagonal} Let $\mathbf{S}=(s_i)_{0\le i < r}$ be a finite sequence over a field $\K$ of length $r$, and let $W_\K(\mathbf{S})=(W(\mathbf{S})[m,n])_{0\le m \le \left\lfloor\frac{r-1}{2}\right\rfloor, m\le n < r-m}$ be the number wall generated by $\mathbf{S}$. For an integer $k\ge0$, the $k^\text{th}$-\textbf{diagonal} is defined as \[D_k(\mathbf{S})=\left\{W_\K(\mathbf{S})[m,m+k]:0\le m\le \left\lceil\frac{r-k}{2}\right\rceil-1\right\}.\] When $\mathbf{S}$ is infinite, one defines $D_k(\mathbf{S})$ similarly as \[D_k(\mathbf{S})=\left\{W_\K(\mathbf{S})[m,m+k]:m\in\N\right\}.\]
    Finally, let $\widetilde W_{\mathbb{K}}(\mathbf{S})[m,k]:=W_{\mathbf{K}}(\mathbf{S})[m,m+k]$ be the entry on row $m$ and diagonal $k$, of the number wall generated by the sequence $\mathbf{S}$.
\end{definition}

\subsection{Identifying Sequences with Laurent Series}\label{sect: equiv_class}\hfill\\
\noindent  The zeroth partial quotient of $\Theta(t)$ (that is, the polynomial part of $\Theta(t)$) is irrelevant to any discussion on escape of mass of $\Theta(t)$ (see Corollary \ref{cor:EscWall}). This motivates an equivalence class on sequences over a field $\K$ defined by the following relationship: \[(s'_i)_{i\in\Z}\sim(s_i)_{i\in\Z} ~~~\text{ if and only if }~~~s'_i=s_i ~~~\text{for all}~~~ i\ge0.\] For convenience, let the sequence $(s_i)_{i\in\Z}$, where $s_i=0$ for all $i<0$, be the representative of each class. For ease, these representatives are denoted by $(s_i)_{i\ge0}$. \\

\noindent Then, one identifies the equivalence class of sequences over $\K$ represented by $\mathbf{S}=(s_i)_{i\ge0}$ with the Laurent series $\Theta(t)=t^{-1}\sum^\infty_{i=0}s_it^{-i}\in\K(\!(t^{-1})\!)$. Additionally, define $T_\Theta(n;m,v):=T_{\mathbf{S}}(n;m,v)$, $W_\K(\mathbf{S})[m,n]=W_\K(\Theta)[m,n]$ and $\widetilde W_\K(\mathbf{S})[m,k]:=\widetilde W_\K(\Theta)[m,k]$. 

 \section{Escape of Mass in Number Walls} \label{Sect: 3}
\noindent This section explains how number walls can be used to calculate the escape of mass of a quadratic irrational. This is achieved via the Lemma \ref{wind=CF}, which shows that for every $k\in\N$ and every sequence $\mathbf{S}=(s_i)_{i\ge0}\in\K^\infty$, the degrees of the partial quotients of $ t^{k-1}\cdot \sum_{i=0}^\infty s_i t^{-i}$ are displayed visually on the number wall of $\mathbf{S}$. Before Lemma \ref{wind=CF} is stated, some preliminary definitions are required.
\subsubsection{\textbf{The Norm and Fractional Part of Laurent Series}}\hfill\\
The following two classical functions on Laurent series are defined: \begin{definition}\label{def: norm+frac}
    Let $h\in\Z$ and $(s_i)_{i\ge -h}$ be a sequence over a field $\K$ with $s_{-h}\neq0$. Then, the \textbf{fractional part} of the Laurent series $\sum_{i=-h}^\infty s_it^{-i}$ is \[\left\langle\sum_{i=-h}^\infty s_it^{-i}\right\rangle:=\sum_{i=\max\{1,-h\}}^\infty s_it^{-i}.\] That is, the fractional part of $\Theta(t)\in\K(\!(t^{-1})\!)$ is obtained by removing all non-negative powers of $t$. \\
    
    \noindent The \textbf{absolute value} of a Laurent series is given by \begin{equation}\left|\sum_{i=-h}^\infty s_it^{-i}\right|:=q^{h},\label{eqn: norm}\end{equation} where $q=\vert \K\vert$ if $\K$ is finite, and $q=2$ otherwise.
\end{definition}
\subsubsection{\textbf{Statement of Lemma \ref{wind=CF}}}\hfill\\
\noindent For a sequence $\mathbf{S}$ over a field $\K$, recall the $k^\nth$ diagonal of $W_\K(\mathbf{S})$ (denoted $D_k(\mathbf{S})$) from Definition \ref{diagonal}.
\begin{lemma}\label{wind=CF}
    Let $\Theta(t)=t^{-1}\sum_{i=0}^\infty s_it^{-i}\in\K(\!(t^{-1})\!)$ be a Laurent series, let $\mathbf{S}=(s_i)_{i\ge0}$ be the sequence of its coefficients. For $k\in\N\cup\{0\}$, denote the continued fraction expansion of $\left\langle t^k\cdot\Theta(t)\right\rangle$ by
    \begin{equation*}
        \left\langle t^k\cdot\Theta(t)\right\rangle=\left[0;A_{1}^{[\Theta\cdot t^k]}(t),A_{2}^{[\Theta\cdot t^k]}(t)(\Theta),\dots\right].
    \end{equation*} 
    Then, $D_k(\mathbf{S})$ is comprised of successive zero portions of length $\deg\left(A_{i}^{[\Theta\cdot t^k]}(t)\right)-1$ separated by nonzero entries. Explicitly, $\widetilde W_{\mathbb{K}}(\mathbf{S})[m,k]\neq 0$ if and only if $\sum_{i=0}^j\deg\left(A_i^{[\Theta\cdot t^k]}(t)\right)=m$ for some $j\in\N$.
\end{lemma}
\noindent Lemma \ref{wind=CF} is depicted by the image below.
    \begin{figure}[H]
        \centering
        \includegraphics[width=0.9\linewidth]{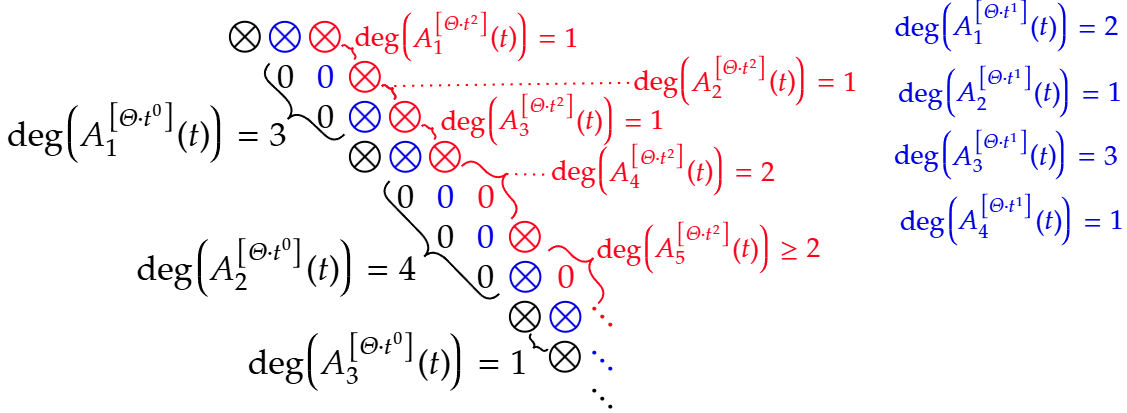}
        \caption{The profile of the first three diagonals of a number wall are depicted above, with the top row having index 0. The first, second and third diagonal are drawn in black, blue and red respectively. The degree of $A_i^{[\Theta\cdot t^k]}(t)$ is given by one plus the number of zeroes between the $i^\nth$ and the $(i+1)^\nth$ nonzero entry of diagonal $k$.}
    \end{figure}

\noindent Section \ref{sec: 4.1} is dedicated to proving Lemma \ref{wind=CF}, whereas Section \ref{sec: 4.2} uses Lemma \ref{wind=CF} to rephrase escape of mass (Definition \ref{def: QuadEscMass}) in terms of number walls.
\subsection{Proving Lemma \ref{wind=CF}}\label{sec: 4.1}\hfill\\
\noindent First, some further definitions are required.
\subsubsection{\textbf{Hankel Matrices}}\label{subsubsec:Hankel}\hfill\\
\noindent To bridge the gap between number walls and the partial quotients of a Laurent series $\Theta(t)\in\K(\!(t^{-1})\!)$, the following family of matrices are defined.
\begin{definition} \label{Hank}
Let $m$ and $n$ be natural numbers. A matrix $(s_{i,j})$ for $0\le i\le n, 0\le j \le m$ is called \textbf{Hankel} if all the entries on any anti-diagonal are equal to one another. Equivalently, $s_{i,j}=s_{i+1,j-1}$ for any $i,j\in\mathbb{N}$ for which this entry is defined.\\

 \noindent Given a doubly infinite sequence $\mathbf{S}:= (s_i)_{i\in\mathbb{Z}}$, natural numbers $m$ and $v$, and an integer $n$, define an $(m+1)\times (v+1)$ Hankel matrix as $H_\mathbf{S}(v):= (s_{i+j+n})_{0\le i \le m, 0\le j \le v}$; namely\\
\begin{align*}H_\mathbf{S}(n;m,v):=\begin{pmatrix}
s_n & s_{n+1} & s_{n+2}& \dots &s_{n+v}\\
s_{n+1} & \iddots &  &\iddots&\vdots\\
s_{n+2}&&\iddots&&\vdots\\
\vdots &&&& \vdots\\
\vdots&&\iddots&&\vdots\\
\vdots&\iddots&&&\vdots\\
s_{n+m} & \dots&\dots& \dots & s_{n+m+v}
\end{pmatrix}.&\end{align*} If $m=v$, this is shortened to $H_\mathbf{S}(n;m)$, respectively. Furthermore, the Laurent series $\Theta(t)=\sum^\infty_{i=1}s_it^{-1}\in\K(\!(t^{-1})\!)$ is identified with the sequence $\mathbf{S}=(s_i)_{i\ge0}$. Accordingly, define $H_\Theta(n;m,v)$\\$:=H_\mathbf{S}(n;m,v)$.
\end{definition}
\noindent It is clear that Hankel and Toeplitz matrices are related. Indeed, by elementary row operations one has that 
\begin{equation}
    \det(H_\textbf{S}(n;m))=(-1)^{\left\lceil\frac{m}{2}\right\rceil} \det(T_\textbf{S}(n+m;m)).\label{hanktoe}
\end{equation}

\subsubsection{\textbf{Number Walls and Approximations of Laurent Series}}\hfill\\
\noindent The connection between approximations to $t^{-1}\cdot\sum_{i=0}^\infty s_it^{-i}\in\K(\!(t^{-1})\!)$ and zero entries in $W_\K((s_i)_{i\ge0})$ is now established.

\begin{lemma}\label{lem: wind=approx}
    Let $q$ be as in Definition \ref{def: norm+frac}, let $\Theta(t)=t^{-1}\cdot \sum_{i=0}^\infty s_it^{-i}\in{\K}(\!(t^{-1})\!)$, let $N(t)\in\K[t]$ be a polynomial of degree $d_N$, and let $l$ and $k$ be natural numbers. Then, the following are equivalent: 
\begin{enumerate}
    \item \label{Cond:Beats_Dirichlet}The following inequality is satisfied: \begin{equation}|N(t)|\cdot \left|\left\langle N(t)\cdot t^k\cdot\Theta(t) \right\rangle\right|\le q^{-l}.\label{eqn: beats_dirichlet}\end{equation}
    \item For every $0\leq i\leq l-2$, one has $\widetilde W_\K(\mathbf{S})[d_N+i,k]=0$. \label{eqn:NumConsqZeros}
\end{enumerate}
\end{lemma}

\noindent To prove Lemma \ref{lem: wind=approx}, the following lemma from the literature is required. 
\begin{lemma}{\cite[Chapter 10]{G} or \cite[Lemma 2.3]{ANL}}\label{lem: famous} Let $\textbf{S}=(s_i)_{i\ge0}$ be a sequence, and let $H:= H_\textbf{S}(1;m-1) = (s_{i+j-1})_{1\le i,j\le m}$ be the $m \times m$ Hankel matrix with entries in a field $\mathbb{K}$. Assume that
the first $r$ columns of $H$ are linearly independent, but that the first $r + 1$ columns are linearly dependent (here, $1 \le r \le m-1$). Then, the principal minor of order r, that is, $\det(s_{i+j-1})_{1\le i,j\le r}$, does not vanish.
    \end{lemma}

\begin{proof}[Proof of Lemma \ref{lem: wind=approx}]
    \textit{Implication} (\ref{Cond:Beats_Dirichlet})$\Rightarrow:$(\ref{eqn:NumConsqZeros}) From the definition of the norm over function fields (Definition \ref{def: norm+frac}), equation (\ref{eqn: beats_dirichlet}) holds if and only if \begin{equation}\label{eqn: Lem_3.4_1}
        \left|\left\langle\Theta(t)\cdot t^k\cdot N(t)\right\rangle\right|\le q^{-l-d_N}.
    \end{equation}
    \noindent Let $n_i\in\K$ (for $0\le i \le d_N$) be such that \begin{equation}\sum_{i=0}^{d_N}n_it^i:=N(t).\label{eqn: Lem_3.4_3}\end{equation} Then, the coefficient of $t^{-i}$ in $\langle \Theta(t)\cdot t^k\cdot N(t)\rangle$ is given by \[\sum_{j=0}^{d_N}s_{i+k+j-1}n_j.\] Equation (\ref{eqn: Lem_3.4_1}) is equivalent to demanding that the coefficient of $t^{-i}$ in $\langle \Theta(t)\cdot t^k\cdot N(t)\rangle$ is equal to zero for all $1\le i\le l+d_N-1$. That is,\begin{equation}
        H_\Theta (k;d_N+l-2, d_N)\cdot \begin{pmatrix}
            n_0\\n_1\\\vdots\\n_{d_N}
        \end{pmatrix}=\mathbf{0}.\label{eqn: Lem_3.4_2}
    \end{equation} This implies that the $d_N+1$ columns of $H_\Theta (k;d_N+l-2, d_N)$ are linearly dependent. Therefore \begin{equation}\label{eqn:hank_det}\det(H_\Theta(k;d_N+i))=0,\end{equation} for every $0\le i \le l-2$. By equation (\ref{hanktoe}), one has that \begin{equation}\det(T_\Theta(k+d_N+i;d_N+i ))=0,\label{eqn: toe_det}\end{equation} for every $0\leq i\leq l-2$. By the definition of the number wall (Definition \ref{nw}), this implies Lemma \ref{lem: wind=approx}(\ref{eqn:NumConsqZeros}).\\

    \noindent \textit{Proof of implication} (\ref{eqn:NumConsqZeros})$\Rightarrow$(\ref{Cond:Beats_Dirichlet}): Assume that for $0\le i \le l-2$, equation (\ref{eqn: toe_det}) holds. By equation (\ref{hanktoe}), equation (\ref{eqn:hank_det}) also holds. Consider the matrices $H_\Theta(k;d_N+l-2)$ and $H_\Theta(k;d_N+l-3)$. Since $\det(H_\Theta(k;d_N+l-3))=0,$ the negation of Lemma \ref{lem: famous} with $r=l-2$ implies that at least one of the following two conditions is false: \begin{enumerate}
        \item The first $d_N+l-2$ columns of $H_\Theta(k;d_N+l-2)$ are linearly independent;
        \item The first $d_N+l-1$ columns of $H_\Theta(k;d_N+l-2)$ are linearly dependent. 
    \end{enumerate}  Since $\det(H_\Theta(k;d_N+l-2))=0$, condition (2) is true. Hence, condition (1) does not hold. That is, the columns of $H_\Theta(k;d_N+l-2,d_N+l-3)$ are linearly dependent. \\

    \noindent By repeating the process but with $r=l-3$ and using the fact that $\det(H_\Theta(k;d_N+l-4))=0$, one obtains that $H_\Theta(k;d_N+l-2,d_N+l-4)$ has linearly dependent columns. Repeating this step a total of $l-2$ times implies that equation (\ref{eqn: Lem_3.4_2}) holds for some vector $\mathbf{n}\in\K^{d_N+1}$. Defining $N(t)$ as in equation (\ref{eqn: Lem_3.4_3}) and using the equivalence between \eqref{eqn: Lem_3.4_1} and \eqref{eqn: Lem_3.4_2} completes the proof. 
\end{proof}

\subsubsection{\textbf{Classical Results on Continued Fractions over Function Fields}}\hfill\\
\noindent To obtain Lemma \ref{wind=CF} from Lemma \ref{lem: wind=approx}, further well-established results are required from the theory of continued fractions over function fields. Let $$\Theta(t)=\left[A_0^{[\Theta]}(t);A_1^{[\Theta]}(t),\dots\right]\in\K(\!(t^{-1})\!)$$ be a Laurent series. Then, the $i^\nth$ \textbf{convergent} of $\Theta(t)$ is the rational function defined by truncating the continued fraction expansion of $\Theta(t)$ after the first $i$ entries. That is, \begin{equation}\nonumber\frac{M_i(t)}{N_i(t)}:=\left[A_0^{[\Theta]}(t);A_1^{[\Theta]}(t),\dots, A_i^{[\Theta]}(t)\right]\in\K(t).\end{equation} Then, it is well known (\cite[Chapter 3.1]{BN}, for example) that \begin{equation}\deg(N_i(t))=\sum_{j=1}^{i} \deg\left(A_j^{[\Theta]}(t)\right).\label{eqn: CF_Fact_3}\end{equation}Furthermore, one has that\begin{equation}|N_i(t)|\cdot |\langle\Theta(t)\cdot N_i(t)\rangle|=\left|A^{[\Theta]}_{i+1}(t)\right|^{-1},\label{eqn: CF_fact_1}\end{equation} and hence \begin{equation}\label{eqn: CF_frac_2}\left|\Theta(t)-\frac{M_i(t)}{N_i(t)}\right|=\frac{1}{|N_i(t)|\cdot |N_{i+1}(t)|}\cdotp\end{equation}Finally, the convergents give the `best possible' approximations to $\Theta(t)$ with respect to the degree of the denominator. That is, for any polynomial $N(t)$ satisfying $\deg(N(t))\le\deg(N_i(t))$,\begin{equation}\label{eqn: best_possible}\left|\Theta(t)-\frac{M(t)}{N(t)}\right|\ge\left|\Theta(t)-\frac{M_i(t)}{N_i(t)}\right|.\end{equation}
\subsubsection{\textbf{Proof of Lemma \ref{wind=CF}}}\hfill\\
\begin{proof}[\unskip\nopunct]
    For $k\in \mathbb{N}$, let 
    $$\frac{M_{i,k}(t)}{N_{i,k}(t)}:=\left[A_0^{[\Theta(t)\cdot t^k]}(t);A_1^{[\Theta(t)\cdot t^k]}(t),\dots,A_i^{[\Theta(t)\cdot t^k]}(t)\right].$$
    By combining equation (\ref{eqn: CF_fact_1}) with Lemma \ref{lem: wind=approx}, one obtains that for every $i\in\N$ and every $0\le j \le \deg\left(A_i^{[\Theta\cdot t^k]}(t)\right)-2$, $\widetilde W_\K(\mathbf{S})[\deg(N_i(t))+j,k]=0$. Therefore, by equation (\ref{eqn: CF_Fact_3}), $D_k(\mathbf{S})$ is comprised of $\deg\left(A_i^{[\Theta\cdot t^k]}(t)\right)-1$ consecutive zero entries separated by elements of $\K$. These appear on rows with index $\deg(N_i(t))-1$. This is illustrated below.
    \begin{figure}[H]
        \centering
        \includegraphics[width=0.75\linewidth]{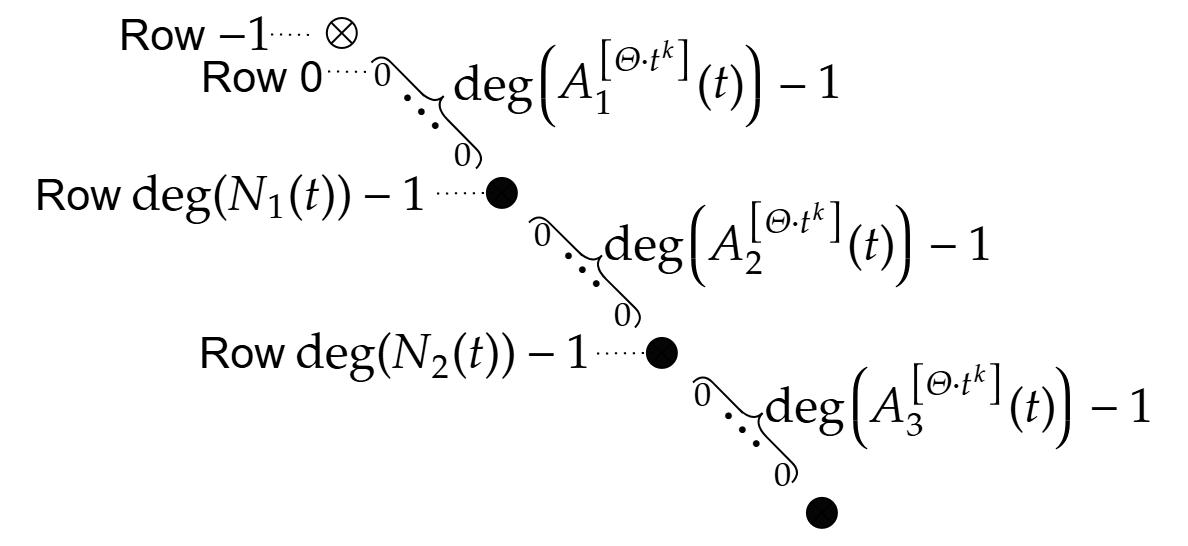}
        \caption{An illustration of diagonal $k$ of the number wall. The large dots represent a generic entry of $\K$. The entry on row $-1$ is nonzero from Definition \ref{nw}.}
    \end{figure}
    \noindent It remains to show that, for every $i$, the $\deg(N_i(t))-1^\nth$ entry of the $k^\nth$ diagonal is nonzero. This is achieved by contradiction. \\

    \noindent Indeed, assume that the entry on row $\deg(N_i(t))-1$ were equal to zero. Then, $D_k(\mathbf{S})$ contains a string of $\deg\left(A_{i}^{[\Theta\cdot t^k]}(t)\right)+\deg\left(A_{i+1}^{[\Theta\cdot t^k]}(t)\right)-1$ zero entries starting on row $\deg(N_{i-1}(t))$. By Lemma \ref{lem: wind=approx}, this is equivalent to the existence of some polynomial $N(t)$ with $\deg(N(t))=\deg(N_{i-1}(t))$ such that  \begin{equation}
        |N(t)|\cdot \left|\left\langle \Theta(t)\cdot t^k\cdot N(t)\right\rangle\right|\le q^{-\left(\deg\left(A_{i}^{[\Theta\cdot t^k]}(t)\right)+\deg\left(A_{i+1}^{[\Theta\cdot t^k]}(t)\right)\right)} .
    \end{equation}
    \noindent This implies that there exists a polynomial $M(t)\in\K[t]$ such that \begin{align}\left|\Theta(t)-\frac{M(t)}{N(t)}\right|\le&~ q^{-\deg\left(A_{i}^{[\Theta\cdot t^k]}(t)\right)-\deg\left(A_{i+1}^{[\Theta\cdot t^k]}(t)\right)-2\cdot \deg(N_{i-1}(t))}.\label{eqn: lem_3.1_proof}\end{align}\noindent By equation (\ref{eqn: CF_frac_2}), \[\deg(N_{i-1}(t))+\deg\left(A_i^{[t^k\Theta]}(t)\right)=\deg(N_i(t)).\] Hence,\begin{align}(\ref{eqn: lem_3.1_proof})&\substack{~~(\ref{eqn: CF_Fact_3})~~\\=}\frac{1}{|N_{i-1}(t)|\cdot |N_{i}(t)|\cdot \left|A^{[\Theta\cdot t^k]}_{i+1}(t)\right|}\nonumber\\&<\frac{1}{|N_{i-1}(t)|\cdot |N_{i}(t)|}\nonumber\\&=\left|\Theta(t)-\frac{M_{i-1}(t)}{N_{i-1}(t)}\right|,\nonumber\end{align} which contradicts equation \eqref{eqn: best_possible}.
\end{proof}
\subsection{Calculating Escape of Mass using Number Walls}\label{sec: 4.2}\hfill\\
\noindent Let $\Theta(t)=t^{-1}\cdot\sum_{i=0}^\infty s_it^{-i}\in\K(\!(t^{-1})\!)$ be a quadratic irrational, and let $\textbf{S}=(s_i)_{i\ge0}$ be the sequence of its coefficients. A simple corollary of Lemma \ref{wind=CF} is that $D_k(\mathbf{S})$ is eventually periodic for any $k\in\N$. This enables the following rephrasing of Definition \ref{def: QuadEscMass}: \begin{corollary}   
\label{cor:EscWall}
Let $\Theta(t)=t^{-1}\cdot \sum_{i=0}^\infty s_it^{-i}\in\K(\!(t^{-1})\!)$ be a quadratic irrational, and let $\textbf{S}=(s_i)_{i\ge0}$ be the sequence of its coefficients. For $k$ a non-negative integer, let $H_k(W_\K(\mathbf{S}))=(h_{k,i})_{0\le i \le \ell_k}$ be a finite sequence, where $h_{k,i}$ is the row index of the $i^\nth$ nonzero entry in the periodic part of $D_k(\mathbf{S})$. Then, for $c\in (0,1]$, $\Theta(t)$ exhibits $c$ escape of mass with respect to the sequence $\{t^i\}_{i\ge0}$ if and only if
    \begin{equation}
    \label{eqn:EscWall}
    \lim_{n\rightarrow\infty}\liminf_{k\rightarrow \infty}\frac{\sum_{i=1}^{\ell_{k}}\max\left\{h_{k,i}-h_{k,i-1}-n,0\right\}}{h_{k,\ell_k}-h_{k,0}}\ge c.
    \end{equation}
\end{corollary}
\begin{proof}
    By Lemma \ref{wind=CF}, $h_{k,i}-h_{k,i-1}=\deg\left(A_i^{[\Theta\cdot t^k]}(t)\right)$. To conclude the proof, note that the denominator of the equation in Definition \ref{def: QuadEscMass} becomes $\sum_{i=1}^{\ell_k}(h_{k,i}-h_{k,i-1})=h_{k,\ell_k}-h_{k,0}$.
\end{proof}

\noindent Therefore, a complete understanding of the profile (recall Definition \ref{profile_def}) of $W_\K(\textbf{S})$ enables the calculation of the escape of mass of $\Theta(t)$ with respect to the sequence $\{t^i\}_{i\ge0}$. 
\section{The $p$-Cantor Sequence and the Profile of its Number Wall}
\label{sec:pCantor}
\noindent The $p$-Cantor sequence is an automatic sequence, the profile of its number wall over $\F_p$ is a 2-dimensional automatic sequence \cite[Theorem 5.1]{AR}. Automatic sequences are a useful tool to construct quadratic irrationals due to Christol's Theorem \cite{Ch}. To define the $p$-Cantor sequence, automatic sequences are introduced. 
\subsection{Automatic Sequences}\hfill\\
\noindent An \textbf{alphabet} $\Sigma$ is a finite set of symbols. Given a finite alphabet $\Sigma$, a \textbf{finite word} over $\Sigma$ is a concatenation of finitely many letters from $\Sigma$. Denote the set of finite words over $\Sigma$ by $\Sigma^*$. 
\begin{definition}
    Given two words $w=a_1\dots a_\ell$ and $w'=a_1'\dots a_m'$ in $\Sigma^*$, define the concatenation operator by
    $$w*w'=a_1\dots a_\ell a_1'\dots a_m'\in \Sigma^*.$$
    Then, $\Sigma^*$ is a monoid with respect to the concatenation, where the identity element is the empty word $\varepsilon$.
\end{definition}
\noindent Denote the set of infinite words over $\Sigma$ by $\Sigma^{\infty}$, and $\Sigma^{\Omega}:=\Sigma^*\cup \Sigma^{\infty}$.
\subsubsection{One Dimensional Automatic Sequences}
\begin{definition}
    For a word $w\in \Sigma^{\Omega}$, denote the \textbf{length} of $w$ by $L(w)$. If $w\in \Sigma^{\infty}$, define $L(w)=\infty$.
\end{definition}
\begin{definition}
    Let $\Sigma,\Delta$ be finite alphabets. A \textbf{coding} is a function $\tau:\Sigma^{\Omega}\rightarrow \Delta^{\Omega}$, satisfying
    $$\tau(\sigma_0\sigma_1\dots)=\tau(\sigma_0)\tau(\sigma_1)\dots.$$
    Note that to define a coding $\tau$ it suffices to define $\tau(a)$ for every $a\in \Sigma$.\\ 

    \noindent A coding is called a \textbf{$d$ coding} if for every $w\in \Sigma^*$, one has $L(\tau(w))=dL(w)$. If $\Delta=\Sigma$, the $d$-coding $\tau$ is called a \textbf{uniform $d$-morphism}\footnote{These definitions are non-standard, but they are more convenient for the purposes of this paper.}. \\

    \noindent Let $\phi:\Sigma^{\Omega}\rightarrow \Sigma^{\Omega}$ be a uniform $d$-morphism. A letter $a\in \Sigma$ is called \textbf{$\phi$-prolongable} if the first letter of $\phi(a)$ is $a$. If $w=\sigma_0\sigma_1\dots \in \Sigma^{\infty}$ satisfies $\phi(w)=w$, then, $\sigma_0$ is $\phi$-prolongable. In this case, $w$ is a \textbf{fixed point} of $\phi$, denoted $\phi^{\infty}(\sigma_0)=w$. Furthermore, $w=\lim_{n\rightarrow \infty}\phi^n(\sigma_0)$, where $\phi^n(\sigma_0)$ is defined inductively as $\phi^n(\sigma_0)=\phi\left(\phi^{n-1}(\sigma_0)\right)$.
\end{definition}
\noindent For the purpose of this paper, the following definition of automaticity arising from Cobham's theorem \cite{C} is used.
\begin{definition}
    A sequence $\textbf{S}$ is $k$-\textbf{automatic} if there exist natural numbers $k,d\in\mathbb{N}$ and finite alphabets $\Sigma$ and $\Delta$, such that $\textbf{S}$ is the fixed point of some uniform $k$-morphism $\phi:\Sigma^\Omega\to\Sigma^\Omega$ under a $d$-coding $\tau:\Sigma^\Omega\to\Delta^\Omega$.
\end{definition}
\subsubsection{\textbf{Two Dimensional Automatic Sequences}}\hfill\\
\noindent Words have a natural generalization to higher dimensions. For the purposes of this work, only two-dimensional automatic sequences are required and thus, are defined below. \\

\noindent Let $\Sigma$ be a finite alphabet. A \textbf{two-dimensional finite word} $\mathbf{w}$ made up from $\Sigma$ is defined as a matrix $$\textbf{w}=(\sigma_{m,n})_{0\le m \le M, 0\le n \le N} ~~~~\text{ where }~~~~ M,N\ge0,$$
and $\sigma_{m,n}\in \Sigma$. In agreement with standard matrix terminology, $\sigma_{m,n}\in\Sigma$ denotes the entry of $\mathbf{w}$ in column $n$ and row $m$. \\

\noindent Let $\Sigma^*_2$ denote the set of finite two-dimensional words made up from $\Sigma$. Similarly, let $\Sigma^\infty_2$ be the set of two-dimensional infinite words, each of which is an infinite two-dimensional matrix $(\sigma_{m,n})_{m,n\in\N}$. Define $\Sigma_2^\Omega:=\Sigma_2^*\cup\Sigma_2^\infty$. \\

\noindent Given $k$ and, $l\in\N$, a \textbf{two-dimensional uniform $[k,l]$-coding }$\tau:\Sigma_2^\Omega\to\Sigma_2^\Omega$ is a function that satisfies the following two properties: \begin{enumerate}
    \item For any letter $\sigma\in\Sigma$, $\tau(\sigma)$ returns a $k\times l$ matrix in $\Sigma^*_2$.
    \item Let $M,N\in\N$ and let $\textbf{w}=(\sigma_{m,n})_{0\le m \le M,0\le n \le N}\in\Sigma_2^*$ be a two-dimensional word. Then, the function $\tau$ satisfies \[\tau\begin{pmatrix}
        \sigma_{0,0}&\sigma_{0,1}&\cdots &\sigma_{0,N}\\\vdots&&&\vdots\\\sigma_{M,0}&\sigma_{M,1}&\cdots&\sigma_{M,N}
    \end{pmatrix}=~\begin{pmatrix}
        \tau(\sigma_{0,0})&\tau(\sigma_{0,1})&\cdots &\tau(\sigma_{0,N})\\\vdots&&&\vdots\\\tau(\sigma_{M,0})&\tau(\sigma_{M,1})&\cdots&\tau(\sigma_{M,N})
    \end{pmatrix}\]
\end{enumerate} Similarly to the one-dimensional case, a $[k,l]$-coding is called a $[k,l]$-\textbf{morphism} if $\Delta=\Sigma$. \\

\noindent A letter $\sigma_0\in\Sigma$ is called $\phi$\textbf{-prolongable} if $\phi(\sigma_0)_{0,0}=\sigma_0$. Just as in the one-dimensional case, define $\phi^k(\sigma_0)$ iteratively as $\phi^k(\sigma_0)=\phi(\phi^{k-1}(\sigma_0))$ and $\phi^0(\sigma_0)=\sigma_0$. Then, the limiting sequence $\phi^\infty(\sigma_0)=\lim_{k\to\infty}\phi^k(\sigma_0)$ exists if and only if $\sigma_0$ is $\phi$-prolongable. \\

\noindent In the two dimensional case, automatic sequences are defined in a similar manner using a generalization of Cobham's theorem to higher dimensions \cite[Theorem 14.2.3]{AllSha}.
\begin{definition}
    A two dimensional sequence $\mathbf{S}$ is $[k,l]$-\textbf{automatic} if there exist natural numbers $k,l,k',l'$ such that $\mathbf{S}$ is the image (under a $[k',l']$-coding) of a fixed point of a two-dimensional uniform $[k,l]$-morphism.
\end{definition}
\subsection{The $p$-Cantor Sequence}\hfill\\
\subsubsection{\textbf{The} $p$\textbf{-Morphism that Defines the} $p$\textbf{-Cantor Sequence}}\hfill\\
\noindent Let $\Sigma$ be a finite alphabet and let $w\in\Sigma^*$. For a natural number $i<L(w)$,  let $w_i$ denote the $i^\nth$ letter of $w$. 
\begin{definition}\label{pcant}
    Let $p$ be an odd prime and let $$\Gamma_p=\{0,1,\dots,p-1\}$$ be the finite alphabet made up from the elements of $\F_p$ and associate $n\in\Gamma_p$ with $n\in\F_p$. Define a $p$-morphism $\phi_p:\Gamma_p^\Omega\to\Gamma_p^\Omega$ as \begin{equation}
        \phi_p(n)_i = \begin{cases}n\cdot \begin{pmatrix}p_2 \\\frac{i}{2}\end{pmatrix}\mod p &~\text{ if }~i\equiv0\mod2\\0&~\text{ otherwise}\end{cases}
    \end{equation} 
for every $n\in\F_p$ and for $0\le i <p$. Then, the \textbf{$\mathbf{p}$-Cantor sequence} $\mathbf{C}^{(p)}=\left(c^{(p)}_i\right)_{i\ge0}$ is defined as $\lim_{n\to\infty}\phi_p^n(1)$ over $\F_p$. Furthermore, to adapt $\mathbf{C}^{(p)}$ into a doubly infinite sequence, define $c^{(p)}_i=0$ for all $i<0$. Define $\Theta_p(t):=t^{-1}\sum_{i=0}^\infty c^{(p)}_it^{-i}\in\F_p(\!(t^{-1})\!)$ to be the Laurent series associated to $\mathbf{C}^{(p)}$.
\end{definition}
\begin{example}
    Let $\Gamma_7=\{0,1,2,3,4,5,6\}$. Then, the 7-Cantor sequence is given by the 7-morphism $\phi:\Gamma_7^\Omega\to\Gamma_7^\Omega$ defined by:\begin{align*}
         &\phi_7(0)=0000000& &\phi_7(A)=1030301& &\phi_7(2)=2060602& \\&\phi_7(3)=3020203& &\phi_7(4)=4050504& &\phi_7(5)=5010105& \\&&&\phi_7(6)=6040406.&
    \end{align*}
\end{example}

\subsubsection{\textbf{The Laurent Series} $\Theta_p(t)$ \textbf{is a Quadratic Irrational}}\hfill\\
\noindent Before diving deeper into the properties of $\mathbf{C}^{(p)}$, it is confirmed that the Laurent series $\Theta_p(t)$ is quadratic over $\F_p$. For brevity, define $p_2=\frac{p-1}{2}$.
\begin{lemma}\label{quad}
    For any odd prime $p$, the Laurent series $\Theta_p(t)\in\F_p(\!(t^{-1})\!)$ is a quadratic irrational. 
\end{lemma}
\begin{proof}
\noindent Instead of working with $\Theta_p(t)$ directly, the proof instead works with the Laurent series $\widetilde\Theta_p(t):=t\cdot \Theta_p(t)$. Clearly, $\widetilde\Theta_p(t)$ is a quadratic irrational if and only if $\Theta_p(t)$ is a quadratic irrational.\\

\noindent It suffices to prove that \begin{equation}\label{quad2}\widetilde\Theta_p(t)= \widetilde\Theta_p(t)^p \sum_{i=0}^{p_2} {p_2 \choose i}t^{-2i}.\end{equation} Indeed, rearranging this equation gives \begin{equation}\label{quad1}\widetilde\Theta(t)^{p-1}=\left(\sum_{i=0}^{p_2} {p_2 \choose i}  t^{-2i}\right)^{-1}.\end{equation} The denominator of the right hand side of equation (\ref{quad1}) is the binomial expansion of $$(1+t^{-2})^{p_2}.$$ Thus, \[\widetilde\Theta_p(t)=(1+t^{-2})^{-\frac{1}{2}},\] so that $\widetilde\Theta_p(t)$ is quadratic. All that remains is to prove (\ref{quad2}). To achieve this, one extends the definition of a $p$-morphism to apply to Laurent series in the following natural way: \\
    
    \noindent Let $\phi:\Gamma_p^\Omega\to\Gamma_p^\Omega$ be a $p$-morphism and let $\Theta(t)=\sum_{i=0}^\infty s_it^{-i}\in\F_p(\!(t^{-1})\!)$ be a Laurent series whose coefficients are given by a sequence $(s_i)_{i\in\N}$ over $\F_p$. Then, define \begin{equation}\phi(\Theta(t))=\sum_{i=0}^\infty t^{-pi}\sum_{j=0}^{p-1} (\phi(s_i)_j t^{-j}).\label{eqn: morph_to_LS}\end{equation}  Let $\phi_p$ be as in Definition \ref{pcant}. As $\mathbf{C}^{(p)}=\phi_p\left(\mathbf{C}^{(p)}\right)$ by definition, one has that $\widetilde\Theta_p(t)=\phi_p(\widetilde\Theta_p(t))$ and hence, \begin{align*}\widetilde\Theta_p(t):=\sum_{j=0}^\infty c^{(p)}_j t^{-j}=&\sum_{j=0}^{\infty}\sum_{i=0}^{p_2} {p_2\choose i}c_j^{(p)}t^{-(pj+2i)}\\
    =&\sum_{i=0}^{p_2} {p_2 \choose i} \sum_{j=0}^\infty c_j^{(p)} t^{-(pj+2i)}\\
    =&\sum_{i=0}^{p_2} {p_2 \choose i} \widetilde\Theta_p(t^p) t^{-2i}.\end{align*} Finally, as $\Theta(t^p)=\Theta(t)^p$ holds for any Laurent series $\Theta(t)\in\F_p(\!(t^{-1})\!)$, the proof is complete.
\end{proof}

\subsection{The Morphism that Generates The Right-Side of $W_p(\mathbf{C}^{(p)})$}\hfill\\
\noindent This subsection details the finite alphabet $\A$, the $[p,p]$-morphism $\Phi_p:\A_2^\Omega\to\A_2^\Omega$ and the $[1,1]$-coding $\Pi:\A_2^\Omega\to\{0,X\}_2^\Omega$ that generates $\left(\chi(W_p(\mathbf{C}^{(p)})[m,n])\right)_{m,n\ge0}$ as shown in \cite[Theorem 5.1]{AR}. The alphabet $\A$ is defined first.
\subsubsection{\textbf{The Finite Alphabet} $\mathcal{A}$}\hfill\\
\noindent Regardless of the choice of the prime $p$, the alphabet $\A$ consists of $12$ letters - $8$ of which are labeled by the cardinal and intercardinal directions: \begin{equation*}
        \mathcal{A}:=\{A,~ B,~ F,~ 0,~ E_N,~ E_E,~ E_S,~ E_W,~ C_{NE},~ C_{SE},~ C_{SW},~ C_{NW}\}.
    \end{equation*}
\noindent To make the roles of the letters of $\mathcal{A}$ more intuitive, the alphabet $\mathcal{A}$ is split into four sub-alphabets. 
\begin{align*}
    &\A_{units}:=\{A,B\},& &\A_{zeroes}=\{F,0\},& \\&\A_{edges}=\{E_N,E_E,E_S,E_W\},& &\A_{corners}=\{C_{NE},C_{SE},C_{SW},C_{NW}\}.&
\end{align*} 
\noindent The image of each letter under $\Phi_p$ and $\Pi$ is both explicitly given and illustrated. With the exception of $A,B$ and $F$, the images of all the above letters under $\Phi_p$ and $\Pi$ are very simple. For those three letters, it is recommended to look at the picture \textit{before} reading the definition.\\

\noindent For the remainder of this section, $0\le m,n\le p-1$ are natural numbers and let $\Phi_p(\cdot)_{m,n}$ is the item in the $m^\nth$ row and the $n^\nth$ column of the image of the input under $\Phi_p$.
\subsubsection{\textbf{The Image of} $\Phi_p$ \textbf{on} $\A_{units}$}\hfill\\
\noindent The $[p,p]$-morphism $\Phi_p$ acts on the letters $A$ and $B$ as,  \begin{align*}
    &\Phi_p(A)_{m,n}=\begin{cases}A &\text{ if } m\equiv 0 \equiv n \mod 2\\
    0 &\text{ if } m\equiv 0\not\equiv  n\mod 2\\
    F &\text{ if } m\not\equiv 0\equiv  n\mod 2\\
    B &\text{ if } m\equiv 1\equiv  n\mod 2\\\end{cases}& &\Phi_p(B)_{m,n}=\begin{cases}
    B &\text{ if } m\equiv 0 \equiv n \mod 2\\
    F &\text{ if } m\equiv 0\not\equiv  n\mod 2\\
    0 &\text{ if } m\not\equiv 0\equiv  n\mod 2\\
    A &\text{ if } m\equiv 1\equiv  n\mod 2\\\end{cases}&,
\end{align*}
\noindent that is, \begin{align*}
    &\Phi_p(A):=~\begin{matrix}
        A&0&A&0&A &\cdots& A&0&A\\
        F&B&F&B&F&\cdots &F&B&F\\
        A&0&A&0&A &\cdots& A&0&A\\
        F&B&F&B&F&\cdots &F&B&F\\
        A&0&A&0&A &\cdots& A&0&A\\
        \vdots &\vdots&\vdots&\vdots&\vdots&&\vdots&\vdots&\vdots\\
        A&0&A&0&A &\cdots& A&0&A\\
        F&B&F&B&F&\cdots &F&B&F\\
        A&0&A&0&A &\cdots& A&0&A\\
    \end{matrix}& &\Phi_p(B):=~\begin{matrix}
    B&F&B&F&B&\cdots&B&F&B\\
    0&A&0&A&0&\cdots &0&A&0\\
    B&F&B&F&B&\cdots&B&F&B\\
    0&A&0&A&0&\cdots &0&A&0\\
    B&F&B&F&B&\cdots&B&F&B\\
    \vdots &\vdots&\vdots&\vdots&\vdots&&\vdots&\vdots&\vdots\\
    B&F&B&F&B&\cdots&B&F&B\\
    0&A&0&A&0&\cdots &0&A&0\\
    B&F&B&F&B&\cdots&B&F&B\\
    \end{matrix}.
\end{align*}
\subsubsection{\textbf{The Image of} $\Phi_p$ \textbf{on} $\A_{zeroes}$}\hfill\\
\noindent Define $\Phi_p(F)$ and $\Phi_p(0)$ as\begin{align*}
    &\Phi_p(F)_{m,n}=\begin{cases}
        C_{NW} &\text{ if }~~~m=n=0,\\
        C_{NE} &\text{ if }~~~m=0~\text{ and }~n=p-1,\\
        C_{SW} &\text{ if }~~~m=p-1~\text{ and }~n=0,\\
        C_{SE} &\text{ if }~~~m=n=p-1,\\
        E_N &\text{ if }~~~m=0~\text{ and }~0\neq n \neq p-1,\\
        E_E &\text{ if }~~~0\neq m\neq p-1~\text{ and }~n=p-1,\\
        E_S &\text{ if }~~~m=p-1~\text{ and }~0\neq n \neq p-1,\\
        E_W &\text{ if }~~~0\neq m\neq p-1~\text{ and }~n=0,\\
        0 &\text{ otherwise},
    \end{cases}& &\Phi_p(0)_{m,n}=0 \text{   for all }m,n,&
\end{align*}\noindent  that is,\begin{align*}
    &\Phi_p(F):=~\begin{matrix}
        C_{NW}&E_N&E_N&\cdots&E_N&E_N&C_{NE}\\
        E_W&0&0&\cdots&0&0&E_E\\
        E_W&0&0&\cdots&0&0&E_E\\
        \vdots&\vdots&\vdots&&\vdots&\vdots&\vdots&\\
        E_W&0&0&\cdots&0&0&E_E\\
        E_W&0&0&\cdots&0&0&E_E\\
        C_{SW}&E_S&E_S&\cdots&E_S&E_S&C_{SE}
    \end{matrix}& &\Phi_p(0):=~\begin{matrix}
        0&0&0&\cdots&0&0&0\\
        0&0&0&\cdots&0&0&0\\
        0&0&0&\cdots&0&0&0\\
        \vdots&\vdots&\vdots&&\vdots&\vdots&\vdots&\\
        0&0&0&\cdots&0&0&0\\
        0&0&0&\cdots&0&0&0\\
        0&0&0&\cdots&0&0&0\\
    \end{matrix}&.
\end{align*}
\begin{remark}
    The letter $F$ was chosen to stand for the word ``Frame", as the image of $F$ under $\Phi_p$ is a $p-2 \times p-2$ square of zeroes \textit{framed} by the letters labelled with the corresponding cardinal and intercardinal directions.
\end{remark}  
\subsubsection{\textbf{The Image of} $\Phi_p$\textbf{ on} $\A_{edges}$}\hfill\\
\noindent The $[p,p]$-morphism $\Phi_p$ acts on the elements of $\A_{edges}$ as,  \begin{align*}
    &\Phi_p(E_N)_{m,n}=\begin{cases}E_N &\text{ if } m= 0 \\
    0 &\text{ otherwise } 
\end{cases}& &\Phi_p(E_E)_{m,n}=\begin{cases}
    E_E &\text{ if } n=p-1\\
    0 &\text{ otherwise } \end{cases}&\\&\Phi_p(E_S)_{m,n}=\begin{cases}E_S &\text{ if } m= p-1 \\
    0 &\text{ otherwise } 
\end{cases}& &\Phi_p(E_W)_{m,n}=\begin{cases}
    E_W &\text{ if } n=0\\
    0 &\text{ otherwise } \end{cases}&,
\end{align*}
\noindent that is,
\begin{align*}
    &\begin{matrix}
        &E_N&E_N&E_N&\cdots&E_N&E_N&E_N\\
        &0&0&0&\cdots&0&0&0\\
        &0&0&0&\cdots&0&0&0\\
       \Phi_p(E_N):=&\vdots&\vdots&\vdots&&\vdots&\vdots&\vdots\\
        &0&0&0&\cdots&0&0&0\\
        &0&0&0&\cdots&0&0&0\\
        &0&0&0&\cdots&0&0&0\\
        \\
        &0&0&0&\cdots&0&0&0\\
        &0&0&0&\cdots&0&0&0\\
        &0&0&0&\cdots&0&0&0\\
        \Phi_p(E_S):=&\vdots&\vdots&\vdots&&\vdots&\vdots&\vdots\\
        &0&0&0&\cdots&0&0&0\\
        &0&0&0&\cdots&0&0&0\\
        &E_S&E_S&E_S&\cdots&E_S&E_S&E_S
    \end{matrix}&  &\begin{matrix}
        &0&0&0&\cdots&0&0&E_E\\
        &0&0&0&\cdots&0&0&E_E\\
        &0&0&0&\cdots&0&0&0_E\\
        \Phi_p(E_E):=&\vdots&\vdots&\vdots&&\vdots&\vdots&\vdots\\
        &0&0&0&\cdots&0&0&E_E\\
        &0&0&0&\cdots&0&0&E_E\\
        &0&0&0&\cdots&0&0&E_E\\
        \\
        &E_W&0&0&\cdots&0&0&0\\
        &E_W&0&0&\cdots&0&0&0\\
        &E_W&0&0&\cdots&0&0&0\\
        \Phi_p(E_W):=&\vdots&\vdots&\vdots&&\vdots&\vdots&\vdots\\
        &E_0&0&0&\cdots&0&0&0\\
        &E_W&0&0&\cdots&0&0&0\\
        &E_W&0&0&\cdots&0&0&0
    \end{matrix}&.
\end{align*}
\subsubsection{\textbf{The Image of} $\Phi_p$ \textbf{on} $\A_{corners}$}\hfill\\
\noindent The $[p,p]$-morphism $\Phi_p$ acts on the elements of $\A_{corners}$ as,  \begin{align*}
    &\Phi_p(C_{NE})_{m,n}=\begin{cases}C_{NE} &\text{ if } m= 0 ~\text{ and }~n=p-1 \\
    E_N &\text{ if } m=0 ~\text{ and }~n\neq p-1\\
    E_E &\text{ if }m\neq0~\text{ and }~n=p-1\\
    0 &\text{ otherwise } 
\end{cases}& &\Phi_p(C_{SE})_{m,n}=\begin{cases}C_{SE} &\text{ if } m=n= p-1 \\
    E_S &\text{ if } m=p-1 \neq n\\
    E_E &\text{ if }m\neq p-1=n~\\
    0 &\text{ otherwise } 
\end{cases}&\\&\Phi_p(C_{SW})_{m,n}=\begin{cases}C_{SW} &\text{ if } m= p-1~\text{ and }~n=0 \\
    E_S &\text{ if } m=p-1 ~\text{ and } ~n\neq0\\
    E_E &\text{ if }m\neq p-1~\text{ and }~n=0\\
    0 &\text{ otherwise } 
\end{cases}& &\Phi_p(C_{NW})_{m,n}=\begin{cases}C_{NW} &\text{ if } m=n= 0 \\
    E_N &\text{ if } m=0 \neq n\\
    E_W &\text{ if }m\neq 0=n~\\
    0 &\text{ otherwise } 
\end{cases}&,
\end{align*}
\noindent that is,
\begin{align*}
    & \Phi_p(C_{NW}):=\hspace{-0.2cm}\begin{matrix}
        C_{NW}&E_N&E_N&\cdots&E_N&E_N&E_N\\
        E_W&0&0&\cdots&0&0&0\\
        E_W&0&0&\cdots&0&0&0\\
       \vdots&\vdots&\vdots&&\vdots&\vdots&\vdots\\
        E_W&0&0&\cdots&0&0&0\\
        E_W&0&0&\cdots&0&0&0\\
        E_W&0&0&\cdots&0&0&0
    \end{matrix}\Phi_p(C_{NE}):=\hspace{-0.2cm}\begin{matrix}
        E_N&E_N&E_N&\cdots&E_N&E_N&C_{NE}\\
        0&0&0&\cdots&0&0&E_E\\
        0&0&0&\cdots&0&0&E_E\\
       \vdots&\vdots&\vdots&&\vdots&\vdots&\vdots\\
        0&0&0&\cdots&0&0&E_E\\
        0&0&0&\cdots&0&0&E_E\\
        0&0&0&\cdots&0&0&E_E
    \end{matrix}&,
\end{align*}
\begin{align*}
    &\Phi_p(C_{SW}):=\begin{matrix}
        E_W&0&0&\cdots&0&0&0\\
        E_W&0&0&\cdots&0&0&0\\
        E_W&0&0&\cdots&0&0&0\\
        \vdots&\vdots&\vdots&&\vdots&\vdots&\vdots\\
        E_W&0&0&\cdots&0&0&0\\
        E_W&0&0&\cdots&0&0&0\\
        C_{SW}&E_S&E_S&\cdots&E_S&E_S&E_S
    \end{matrix}~\Phi_p(C_{SE}):=\begin{matrix}
        0&0&0&\cdots&0&0&E_E\\
        0&0&0&\cdots&0&0&E_E\\
        0&0&0&\cdots&0&0&E_E\\
        \vdots&\vdots&\vdots&&\vdots&\vdots&\vdots\\
        0&0&0&\cdots&0&0&E_E\\
        0&0&0&\cdots&0&0&E_E\\
        E_S&E_S&E_S&\cdots&E_S&E_S&C_{SE}
    \end{matrix}&.
\end{align*}

\subsubsection{\textbf{The} $[1,1]$\textbf{-coding }$\Pi$ \textbf{on} $\A$}\hfill\\
\noindent Let $L\in\A$. Then, the $[1,1]$-coding $\Pi:\mathcal{A}^\Omega_2\to \{0,X\}_2^\Omega$ is defined by,  \begin{align*}
    &\Pi(L)_{m,n}=\begin{cases}
    0 &\text{ if } L = 0\\
    X &\text{ otherwise } \end{cases}& 
\end{align*}
The following Theorem from \cite{AR} lies at the heart of every proof in this paper.
\begin{theorem}{\cite[Theorem 5.1]{AR}}\label{thm: morphism}
    Let $\A$, $\Phi_p:\mathcal{A}^\Omega_2\to\A_2^\Omega$ and $\Pi:\A_2^\Omega\to\{0,X\}_2^\Omega$ be as described above. Then, \[\chi\left(\left(W_p\left(\mathbf{C}^{(p)}\right)[m,n]\right)_{m,n\in\N}\right)=\Pi(\Phi^\infty(A)).\]
\end{theorem}
\section{Properties of the Number Wall of the $p$-Cantor Sequence}
\noindent
The goal of this section is  carefully analyze the morphism $\phi_p$ and the coding $\Pi$, so to gain a precise understanding of $\Phi_p^\infty(A)$. 
\subsection{The Morphism and Group Actions}\hfill\\
\noindent In this section, a group action on $\mathcal{A}$ is defined. To compute escape of mass, the symmetries of this group action are exploited. Let $G=\langle \rho,\eta,\iota|\rho^4,\eta^2,\iota^2\rangle$ be a product of a cyclic group of order $4$ and two cyclic groups of order $2$. 
\begin{enumerate}
\item Then, $\rho$ acts on $\mathcal{A}$ in the following manner:
\begin{align}
    \,\rho(X)&=X,~\text{for }X\in \mathcal{A}_{units}\cup\mathcal{A}_{zeroes},\nonumber\end{align}\vspace{-0.7cm}\begin{align}
    \rho(E_N)&=E_E\,,&\rho(E_E)&=E_S\,,&\rho(E_S)&=E_W\,,&\rho(E_W)&=E_N,\\
    \rho(C_{NW})&=C_{NE}\,,&\rho(C_{NE})&=C_{SE}\,,&\rho(C_{SE})&=C_{SW}\,,&\rho(C_{SW})&=C_{NW}.\nonumber 
\end{align}

\item The $\iota$ action on $\mathcal{A}$ is defined as:
\begin{equation}
\begin{aligned}
    \iota(A)&=B\,,&\iota(B)&=A,\\
    \iota(X)&=X\,,&\text{for }&X\in \mathcal{A}\setminus \mathcal{A}_{units}.
\end{aligned}
\end{equation}
\item The $\eta$ action is defined as:
\begin{equation}
\begin{aligned}
    \eta(0)&=F\,&\eta(F)&=0\\
    \eta(X)=&X,& \text{ for }&X\in \mathcal{A}\setminus\mathcal{A}_{zeroes}.
\end{aligned}
\end{equation}
\end{enumerate}
\subsubsection{\textbf{Reducing the Alphabet}}\hfill\\
Due to these symmetries, $\mathcal{A}$ is described through the action of $G$ on a subset of $\mathcal{A}$, which is denoted by $\mathbf{T}$. Let $\mathbf{T}=\{0,A,E_E,C_{NE},F\}$. To avoid visual clutter, the subscripts are removed to obtain $\mathbf{T}=\{0,A,E,C,F\}$. The motivation behind these symbols is the following
\begin{equation}
\begin{aligned}
    0\leftrightarrow\text{zero}\,,&&A\leftrightarrow \text{first}\,,&&E\leftrightarrow \text{edge}\,,&&C\leftrightarrow \text{corner}\,,&&
    F\leftrightarrow \text{frame}.
\end{aligned}
\end{equation}
The $\rho$ action lifts to blocks of size $p^k\times p^k$ with entries in $\mathbf{T}$ by
\begin{equation}
    \rho\begin{pmatrix}
        a_{0,0}&\dots&a_{0,p^k-1}\\
        \vdots&\ddots&\vdots\\
        a_{p^k-1,0}&\dots&a_{p^k-1,p^k-1}
    \end{pmatrix}=\begin{cases}
        \begin{pmatrix}
        a_{0,0}&\dots&a_{0,p^k-1}\\
        \vdots&\ddots&\vdots\\
        a_{p^k-1,0}&\dots&a_{p^k-1,p^k-1}
    \end{pmatrix}&a_{i,j}\in \{A,B\}\text{ for some }i,j\\
    \begin{pmatrix}
        \rho(a_{p^k-1,0})&\dots&\rho(a_{0,0})\\
        \vdots&\ddots&\vdots\\
        \rho(a_{p^k-1,p^k-1})&\dots&\rho(a_{0,p^k-1})
    \end{pmatrix}&\text{else}\end{cases}.
    \end{equation}
The actions of $\eta$ and $\iota$ also lift to an action on $p^k\times p^k$ blocks, which are described shortly.
\subsubsection{\textbf{The Morphism} $\Phi_p$ \textbf{on} $\mathbf{T}$}\hfill\\
\noindent Using the alphabet $\mathbf{T}$, the morphism $\Phi_p:\mathbf{T}\rightarrow \mathbf{T}^{p\times p}$ is given by
\begin{equation}
\label{eqn:A,iota(A)}
    A\mapsto \begin{matrix}
        A&0&A&\dots&A&0&A\\
        F&\iota(A)&F&\dots&F&\iota(A)&F\\
        A&0&A&\dots&A&0&A\\
        \vdots&\vdots&\vdots&\ddots&\vdots&\vdots&\vdots\\
        A&0&A&\dots&A&0&A\\
        F&\iota(A)&F&\dots&F&\iota(A)&F\\
        A&0&A&\dots&A&0&A
    \end{matrix},\text{     } \iota(A)\mapsto \begin{matrix}
        \iota(A)&F&\iota(A)&\dots&\iota(A)&F&\iota(A)\\
        0&A&0&\dots&A&0&A\\
        \iota(A)&F&\iota(A)&\dots&\iota(A)&F&\iota(A)\\
        \vdots&\vdots&\vdots&\ddots&\vdots&\vdots&\vdots\\
        \iota(A)&F&\iota(A)&\dots&\iota(A)&F&\iota(A)\\
        0&A&0&\dots&A&0&A\\
        \iota(A)&F&\iota(A)&\dots&\iota(A)&F&\iota(A)
    \end{matrix},
\end{equation}
\begin{equation}
0\mapsto \begin{matrix}
    0&0&\dots&\dots&0\\
    0&0&\dots&\dots&0\\
    \vdots&\vdots&\ddots&&\vdots\\
    \vdots&\vdots&&\ddots&\vdots\\ 
    0&0&\dots&\dots&0
\end{matrix}, ~~~~~~~~~F\mapsto \begin{matrix}
        \rho^{-1}(C)&\rho^{-1}(E)&\dots&\rho^{-1}(E)&C\\
        \rho^2(E)&0&\dots&0&E\\
        \vdots&\vdots&\ddots&\vdots&\vdots\\
        \rho(E)&0&\dots&0&E\\
        \rho^2(C)&\rho(E)&\dots&\rho(E)&\rho(C)
    \end{matrix};
\end{equation}
\begin{equation}
    E\mapsto \begin{matrix}
        0&0&\dots&0&E\\
        0&0&\dots&0&E\\
        \vdots&\vdots&\ddots&\vdots&\vdots\\
        0&0&\dots&0&E\\
        0&0&\dots&0&E
    \end{matrix},~~~~~~~~~ C\mapsto \begin{matrix}
        \rho^{-1}(E)&\rho^{-1}(E)&\dots&\rho^{-1}(E)&C\\
        0&0&\dots&0&E\\
        0&0&\dots&0&E\\
        \vdots&\vdots&\ddots&\vdots&\vdots\\
        0&0&\dots&0&E
    \end{matrix}.
\end{equation}
Note that $\Phi_p$ is $\rho$ equivariant on the letters $a\in \mathbf{T}$, since $\Phi_p(\rho^i(E))=\rho^i(\Phi_p(E))$, $\Phi_p(\rho^i(C))=\rho^i(\Phi_p(C))$, and $\rho(\Phi_p(X))=\Phi_p(X)=\Phi_p(\rho(X))$ for $X\in \mathcal{A}_{units}\cup \mathcal{A}_{zeroes}$. 
\newline\newline
\noindent The remainder of this subsection establishes properties of the morphism $\Phi_p$ which are crucial in computing escape of mass in Section \ref{sec:EscComp}. 
\subsubsection{$k$\textbf{-Pseudo Windows}}
\begin{definition}
    For $k\in \mathbb{N}$, denote the set of two dimensional words with entries in $\mathcal{A}$ of side length $p^k\times p^k$ by $\mathcal{A}^{p^k\times p^k}$. A \textbf{large} $k$\textbf{-pseudo-window}\footnote{The word ``window'' comes from the fact that zero regions in a number wall always come in squares, which are called `windows'. Since this is a square of zeroes but it is not within a number wall, the word `pseudo' has been added.}  of size $p^k$ is a square zero block in $\mathcal{A}^{p^k\times p^k}$ and is denoted by $Z_k=\mathbf{0}_{p^k\times p^k}$. A \textbf{small }$k$\textbf{-pseudo-window} is a word $W_k\in \mathcal{A}^{p^k\times p^k}$ of the form
    \begin{equation}
        W_k=\begin{matrix}
            \rho^{3}(C)&\rho^{3}(E)&\rho^{3}(E)&\dots&\rho^{3}(E)&\rho^{3}(E)&C\\
            \rho^2(E)&0&0&\dots&0&0&E\\
            \rho^2(E)&0&0&\dots&0&0&E\\
            \vdots&\vdots&\vdots&\ddots&\vdots&\vdots&\vdots\\
            \rho^2(E)&0&0&\dots&0&0&E\\
            \rho^2(E)&0&0&\dots&0&0&E\\
            \rho^2(C)&\rho(E)&\rho(E)&\dots&\rho(E)&\rho(E)&\rho(C)
        \end{matrix},
    \end{equation}
    Define $W_0=F$ and $Z_0=0$. 
\end{definition}
\noindent The following lemma is an immediate consequence of the definition of $\Phi_p$.
\begin{lemma}
\label{lem:WindTrans}
For every $k\in \mathbb{N}$, one has
    \begin{enumerate}
        \item $\Phi_p(W_k)=W_{k+1}$;
        \item $\Phi_p(Z_k)=Z_{k+1}$'
    \end{enumerate}
\end{lemma}
\subsubsection{\textbf{The Repeated Image of }$\Phi$ \textbf{on the Letters of} $\mathbf{T}$}\hfill\\
The next theorem complements Lemma \ref{lem:WindTrans} by establishing the image of the remaining letters of $\mathbf{T}$ under $\Phi_p^k$.
\begin{theorem}
\label{thm:PhiTrans}
For every $k$, 
\begin{enumerate}
\item \label{eqn:Phi(A)Trans}
    $$\Phi_p^k(A)=~\begin{matrix}
        \Phi_p^{k-1}(A)&Z_{k-1}&\Phi_p^{k-1}(A)&\dots&\Phi_p^{k-1}(A)&Z_{k-1}&\Phi_p^{k-1}(A)\\
        W_{k-1}&\Phi_p^{k-1}(\iota(A))&W_{k-1}&\dots&W_{k-1}&\Phi_p^{k-1}(\iota(A))&W_{k-1}\\
        \Phi_p^{k-1}(A)&Z_{k-1}&\Phi_p^{k-1}(A)&\dots&\Phi_p^{k-1}(A)&Z_{k-1}&\Phi_p^{k-1}(A)\\
        \vdots&\vdots&\vdots&&\vdots&\vdots&\vdots\\
        \Phi_p^{k-1}(A)&Z_{k-1}&\Phi_p^{k-1}(A)&\dots&\Phi_p^{k-1}(A)&Z_{k-1}&\Phi_p^{k-1}(A)\\
        W_{k-1}&\Phi_p^{k-1}(\iota(A))&W_{k-1}&\dots&W_{k-1}&\Phi_p^{k-1}(\iota(A))&W_{k-1}\\
        \Phi_p^{k-1}(A)&Z_{k-1}&\Phi_p^{k-1}(A)&\dots&\Phi_p^{k-1}(A)&Z_{k-1}&\Phi_p^{k-1}(A)
    \end{matrix};$$
    \item \label{eqn:Phi(1bar)Trans}
    $$\Phi_p^k(\iota(A))=\begin{matrix}
        \Phi_p^{k-1}(\iota(A))&W_{k-1}&\Phi_p^{k-1}(\iota(A))&\dots&\Phi_p^{k-1}(\iota(A))&W_{k-1}&\Phi_p^{k-1}(\iota(A))\\
        Z_{k-1}&\Phi_p^{k-1}(A)&Z_{k-1}&\dots&Z_{k-1}&\Phi_p^{k-1}(A)&Z_{k-1}\\
        \Phi_p^{k-1}(\iota(A))&W_{k-1}&\Phi_p^{k-1}(\iota(A))&\dots&\Phi_p^{k-1}(\iota(A))&W_{k-1}&\Phi_p^{k-1}(\iota(A))\\
        \vdots&\vdots&\vdots&&\vdots&\vdots&\vdots\\
        \Phi_p^{k-1}(\iota(A))&W_{k-1}&\Phi_p^{k-1}(\iota(A))&\dots&\Phi_p^{k-1}(\iota(A))&W_{k-1}&\Phi_p^{k-1}(\iota(A))\\
        Z_{k-1}&\Phi_p^{k-1}(A)&Z_{k-1}&\dots&Z_{k-1}&\Phi_p^{k-1}(A)&Z_{k-1}\\
        \Phi_p^{k-1}(\iota(A))&W_{k-1}&\Phi_p^{k-1}(\iota(A))&\dots&\Phi_p^{k-1}(\iota(A))&W_{k-1}&\Phi_p^{k-1}(\iota(A))
    \end{matrix}.$$
\end{enumerate}
\end{theorem}
\begin{proof}[Proof of \eqref{eqn:Phi(A)Trans}]
    The proof proceeds by induction on $k$. For $k=1$, this holds due to the definition of $\Phi_p(A)$. Assume that the claim holds up to $k$. By the definition of $\Phi_p$, the induction hypothesis, and Lemma \ref{lem:WindTrans}, one has
    \begin{equation*}
        \Phi_p^{k+1}(A)=\Phi_p\begin{pmatrix}
        \Phi_p^{k-1}(A)&Z_{k-1}&\Phi_p^{k-1}(A)&\dots&\Phi_p^{k-1}(A)&Z_{k-1}&\Phi_p^{k-1}(A)\\
        W_{k-1}&\Phi_p^{k-1}(\iota(A))&W_{k-1}&\dots&W_{k-1}&\Phi_p^{k-1}(\iota(A))&W_{k-1}\\
        \Phi_p^{k-1}(A)&Z_{k-1}&\Phi_p^{k-1}(A)&\dots&\Phi_p^{k-1}(A)&Z_{k-1}&\Phi_p^{k-1}(A)\\
        \vdots&\vdots&\vdots&&\vdots&\vdots&\vdots\\
        \Phi_p^{k-1}(A)&Z_{k-1}&\Phi_p^{k-1}(A)&\dots&\Phi_p^{k-1}(A)&Z_{k-1}&\Phi_p^{k-1}(A)\\
        W_{k-1}&\Phi_p^{k-1}(\iota(A))&W_{k-1}&\dots&W_{k-1}&\Phi_p^{k-1}(\iota(A))&W_{k-1}\\
        \Phi_p^{k-1}(A)&Z_{k-1}&\Phi_p^{k-1}(A)&\dots&\Phi_p^{k-1}(A)&Z_{k-1}&\Phi_p^{k-1}(A)\\
    \end{pmatrix}\\
    \end{equation*}
    \begin{equation*}
        =\begin{matrix}\Phi_p^{k}(A)&Z_{k}&\Phi_p^{k}(A)&\dots&\Phi_p^{k}(A)&Z_{k}&\Phi_p^{k}(A)\\
        W_{k}&\Phi_p^{k}(\iota(A))&W_{k}&\dots&W_{k}&\Phi_p^{k}(\iota(A))&W_{k}\\
        \Phi_p^{k}(A)&Z_{k}&\Phi_p^{k}(A)&\dots&\Phi_p^{k}(A)&Z_{k}&\Phi_p^{k}(A)\\
        \vdots&\vdots&\vdots&&\vdots&\vdots&\vdots\\
        \Phi_p^{k}(A)&Z_{k}&\Phi_p^{k}(A)&\dots&\Phi_p^{k}(A)&Z_{k}&\Phi_p^{k}(A)\\
        W_{k}&\Phi_p^{k}(\iota(A))&W_{k}&\dots&W_{k}&\Phi_p^{k}(\iota(A))&W_{k}\\
        \Phi_p^{k}(A)&Z_{k}&\Phi_p^{k}(A)&\dots&\Phi_p^{k}(A)&Z_{k}&\Phi_p^{k}(A)\\
    \end{matrix}.
    \end{equation*}
    The proof of \eqref{eqn:Phi(1bar)Trans} is very similar and thus is omitted. 
\end{proof}
\noindent The elements of $\mathcal{A}^{p^{k-1}\times p^{k-1}}$ in $\Phi_p^k(A)$ denoted by $\Phi_p^{k-1}(A)$, $Z_{k-1}$, $W_{k-1}$, and $\Phi_p^{k-1}(\iota(A))$ are called $p^{k-1}\times p^{k-1}$ \textbf{words}.
\subsubsection{\textbf{The Actions} $\iota$ \textbf{and} $\eta$ \textbf{on} $\mathcal{A}^{p^k\times p^k}$}\hfill\\
\begin{definition}
    Define the action of $\eta$ on $p^k\times p^k$ words $X\in \mathcal{A}^{p^k\times p^k}$ by
    \begin{equation}
        \eta(X)=\begin{cases}
            Z_k&X=W_k\\
            W_k&X=Z_k\\
            X&\text{else}
        \end{cases}.
    \end{equation} 
    Define the $\iota$ action recursively. For 
    $$X=\begin{pmatrix}
        X_{0,0}&\dots&X_{0,p-1}\\
        \vdots&\ddots&\vdots\\
        X_{p-1,0}&\dots&X_{p-1,p-1}
    \end{pmatrix}\in \mathcal{A}^{p^{k+1}\times p^{k+1}},$$
    where $X_{i,j}\in \mathcal{A}^{p^k\times p^k}$ for all $0\leq i,j\leq p-1$, define
    \begin{equation}
    \label{eqn:iotaDef}
        \iota\begin{pmatrix}
            X_{0,0}&\dots&X_{0,p-1}\\
            \vdots&\ddots&\vdots\\
            X_{p-1,0}&\dots&X_{p-1,p-1}
        \end{pmatrix}=\begin{cases}
            \Phi_p^{k}(\iota(A))&X=\Phi_p^k(A)\\
            \Phi_p^{k}(A)&X=\Phi_p^k(\iota(A))\\
            \begin{pmatrix}
                \iota(X_{0,0})&\dots&\iota(X_{0,p-1})\\
                \vdots&\ddots&\vdots\\
                \iota(X_{p-1,0})&\dots&\iota(X_{p-1,p-1})
            \end{pmatrix}&\text{else}
        \end{cases}.
    \end{equation}
\end{definition}
\noindent With this definition, $\rho,\eta,\iota$ are $\Phi_p$-equivariant.
\begin{lemma}
\label{lem:Equivar}
    For every $g\in G=\langle \rho,\eta,\iota\rangle$, for every $k\in \mathbb{N}$ and for every $X\in \mathcal{A}^{p^k\times p^k}$, one has $\Phi_p(gX)=g\Phi_p(X)$. 
\end{lemma}
\begin{proof}
    It suffices to verify this for every $X\in \mathcal{A}$. By definition, $\Phi_p(0),\Phi_p(F),\Phi_p(A),\Phi_p(\iota(A))$ are $\rho$ invariant, as are $0$ and $F$, and it is clear that $\Phi_p$ is $\rho$ equivariant on $E$ and $C$.
\newline\newline
    \noindent By the definition of $\eta$, for every word $X\in \mathcal{A}^{p^k\times p^k} \setminus\{\Phi_p^k(F),\Phi_p^k(0)\}$, one has that $\eta(\Phi_p(X))=\Phi_p(X)=\Phi_p(\eta(X))$. Moreover,
    $$\eta(\Phi_p(W_k))=\eta(W_{k+1})=Z_{k+1}=\Phi_p(Z_k)=\Phi_p(\eta(W_k)),$$
    and $$\eta(\Phi_p(Z_k))=\eta(Z_{k+1})=W_{k+1}=\Phi_p(W_k)=\Phi_p(\eta(Z_k)).$$
    
    \noindent Induction is now used to prove that for every $k\geq 0$, $\Phi_p$ is $\iota$ equivariant on $\mathcal{A}^{p^k\times p^k}$. The case $k=0$ holds trivially. Assume that $\Phi_p$ is $\iota$ equivariant on $\mathcal{A}^{p^k\times p^k}$. Let $X\in \mathcal{A}^{p^{k+1}\times p^{k+1}}$. 
    \newline\newline
    \noindent If $X=\Phi_p^k(A)$, then, by the definition \eqref{eqn:iotaDef},
    $$\iota(\Phi_p(X))=\iota(\Phi_p^{k+1}(A))=\Phi_p^{k+1}(\iota(A))=\Phi_p(\Phi_p^{k}(\iota(A)))=\Phi_p(\iota(\Phi_p^k(A)))=\Phi_p(\iota(X)).$$ 
    Similarly, if $X=\Phi_p^k(\iota(A))$, then, $$\iota(\Phi_p(X))=\iota(\Phi_p^{k+1}(\iota(A)))=\Phi_p^{k+1}(A)=\Phi_p(\Phi_p^k(A))=\Phi_p(\iota(X)).$$
    Otherwise, let $X_{i,j}\in \mathcal{A}^{p^k\times p^k}$ for $0\leq i,j\leq p^k-1$, be a word such that 
    $$X=\begin{pmatrix}
        X_{0,0}&\dots&X_{0,p-1}\\
        \vdots&\ddots&\vdots\\
        X_{p-1,0}&\dots&X_{p-1,p-1}
    \end{pmatrix}.$$ 
    Then, by the induction hypothesis,
    $$\Phi_p(\iota(X))=\Phi_p\begin{pmatrix}
        \iota(X_{0,0})&\dots&\iota(X_{0,p-1})\\
        \vdots &\ddots&\vdots\\
        \iota(X_{p-1,0})&\dots&\iota(X_{p-1,p-1})
    \end{pmatrix}=\begin{pmatrix}
        \Phi_p(\iota(X_{0,0}))&\dots&\Phi_p(\iota(X_{0,p-1}))\\
        \vdots&\ddots&\vdots\\
        \Phi_p(\iota(X_{p-1,0}))&\dots&\Phi_p(\iota(X_{p-1,p-1}))
    \end{pmatrix}$$
    $$=\begin{pmatrix}
        \iota(\Phi_p(X_{0,0}))&\dots&\iota(\Phi_p(X_{0,p-1}))\\
        \vdots&\ddots&\vdots\\
        \iota(\Phi_p(X_{p-1,0}))&\dots&\iota(\Phi_p(X_{p-1,p-1}))
    \end{pmatrix}=\iota(\Phi_p(X)).$$
\end{proof}
\subsubsection{\textbf{The Outermost Edges of Words in} $\mathcal{A}^{p^k\times p^k}$}\hfill\\
Next, some claims about the outermost rows and columns of $\Phi_p^k(A)$ are proved. 
\begin{lemma}
\label{lem:row1Phi^k(A)}
For every $k\in \mathbb{N}$,
\begin{enumerate}
    \item \label{lem:Phi^k(A)FirstRow} the first row of $\Phi_p^k(A)$ contains only the letters $0$ and $A$
    \item \label{lem:Phi^k(iota(A))FirstCol} the first column of $\Phi_p^k(\iota(A))$ contains only the letters $0$ and $\iota(A)$.
\end{enumerate}
\end{lemma}
\begin{proof}
    Claim \eqref{lem:Phi^k(A)FirstRow} is proved explicitly and claim \eqref{lem:Phi^k(iota(A))FirstCol} follows from \eqref{lem:Phi^k(A)FirstRow} and Lemma \ref{lem:Equivar} with $g=\rho$. When $k=1$, this is clear by the definition of $\Phi_p(A)$. Assume that the claim holds up to $k$. By Theorem \ref{thm:PhiTrans}, the first row of $\Phi_p^{k+1}(A)$ is composed of copies of the first row of $\Phi_p^k(A)$ and copies of the first row of $Z_k$. By the definition of $Z_k$ and the induction hypothesis, both of these contain only zeroes and $A$'s. Hence, the claim follows. 
\end{proof}
\begin{lemma}
\label{lem:LastColFirstRow}
    For every $k\in \mathbb{N}$,
    \begin{enumerate}
        \item \label{Phi^k(A)LastCol} the first and last columns of $\Phi_p^k(A)$ contains only the letters $F$, $A$, $E$, $C$, and $\rho(C)$;
        \item \label{Phi^k(bar(A))FirstRow} the first and last rows of $\Phi_p^k(\iota(A))$ contains only the letters $F$, $\iota(A)$, $\rho(E)$, $\rho(C)$, and $\rho^2(C)$.
    \end{enumerate}
\end{lemma}
\begin{proof}
    Part \eqref{Phi^k(A)LastCol} is proved, and \eqref{Phi^k(bar(A))FirstRow} follows from \eqref{Phi^k(A)LastCol} and Lemma \ref{lem:Equivar}. By the symmetry of Theorem \ref{thm:PhiTrans}, it suffices to prove the claim for the last column of $\Phi_p^k(A)$. When $k=1$, the claim is true due to the definition of $\Phi_p(A)$. Assume that the claim holds for $k$. By Theorem \ref{thm:PhiTrans}, the last column of $\Phi_p^{k+1}(A)$ consists of copies of the last columns of $\Phi_p^k(A)$ and of $W_k$. By the induction hypothesis and the definition of $W_k$, the last column of $\Phi_p^{k+1}(A)$ consists only of $F,A,E,C,\rho(C)$.
\end{proof}
\subsubsection{\textbf{The Maximal Size of Pseudo Windows}}\hfill\\
\noindent As a consequence, all (large) pseudo windows in $\Phi_p^{k}(A)$ are of size at most $p^{k-1}$.
\begin{lemma}
\label{lem:Phi^k(A)WindSize}
    The largest zero regions in $\Phi_p^{k}(A)$ (respectively $\Phi_p^k(\iota(A))$) are given by large $(k-1)$-psuedo-windows. 
\end{lemma}
\begin{proof}
The proof is completed for $\Phi_p^k(A)$. The method is identical for $\Phi_p^k(\iota(A))$ and is therefore omitted. \\

\noindent A simple induction shows that all zero regions in $\Phi_p^k(A)$ are either large or small $k'$-psuedo-windows for some $0\le k'\le k$. If $\Phi^k_p(A)$ did contain a large or small $k$-pseudo-window, it must intersect at least two distinct $p^{k-1}\times p^{k-1}$ subwords of $\Phi_p^k(A)$. By Lemma \ref{lem:row1Phi^k(A)}, the first and last row of $\Phi_p^{k-1}(\iota(A))$ contain no zeroes. Hence, by Theorem \ref{thm:PhiTrans} and Lemma \ref{lem:row1Phi^k(A)}, a large or small $k$-pseudo-window intersecting $Z_{k-1}$ cannot intersect $\Phi_p^{k-1}(\iota(A))$ and $\Phi_p^{k-1}(A)$. On the other hand, by the definition of $W_{k-1}$, a large or small $k$-pseudo-window intersecting $W_{k-1}$ must be contained in $W_{k-1}$. Hence, by Lemma \ref{lem:row1Phi^k(A)}, all zero windows in $\Phi_p^k(A)$ must be contained in some $p^{k-1}\times p^{k-1}$ subword of $\Phi_p^k(A)$, which is a contradiction.
\end{proof}
\subsubsection{\textbf{Diagonals of} $\Phi_p^{\infty}(A)$}\hfill\\
Let $j\in \mathbb{N}$. Similarly to in Section \ref{subsub:DiagsNW}, define the $j$-th diagonal of $\Phi_p^{\infty}(A)$ by 
$$\mathcal{D}_j(\Phi_p^{\infty}(A))=\left\{\Phi^\infty_p(A)[n,n+j]:n\in \mathbb{N}\cup \{0\}\right\},$$ 
where $\Phi^\infty_p(A)[m,n]$ is the entry in row $m$ and column $n$ of $\Phi^\infty_p(A)$.\\

\noindent For $X\in \mathcal{A}^{p^k\times p^k}$ a subword of $\Phi_p^\infty(A)$ and $j\in \mathbb{N}$, denote the intersection of the $j$-th diagonal of $\Phi_p^{\infty}(A)$ with $X$ by $\mathcal{D}_j(\Phi_p^{\infty}(A))\sqcap X$.
\begin{definition}
    Let $j\in \mathbb{N}$. Then, $\mathcal{D}_j(\Phi_p^{\infty}(A))$ is \textbf{eventually periodic with period $\ell$} if there exists $N\in \mathbb{N}$, such that for every $n\geq N$,
    $$\Phi^{\infty}(A)[n+\ell,n+\ell+j]=\Phi^{\infty}(A)[n,n+j].$$
   Let $\mathcal{P}(j)$ denote the period length of $\mathcal{D}_j(\Phi_p^{\infty}(A))$. If $\mathcal{D}_j(\Phi^{\infty}(A))$ is not eventually periodic, define $\mathcal{P}(j)=\infty$. 
\end{definition}
\noindent The following corollary follows immediately from Theorem \ref{thm:PhiTrans} and plays a crucial role in proving that all the diagonals of $\Phi^{\infty}(A)$ are periodic.
\begin{corollary}
\label{thm:IntZ_kW_k}
    Let $0\leq j\leq p^k-1$. Then, the first $2p^k$ entries of $\mathcal{D}_j(\Phi^{\infty}(A))$ are split into the four following consecutive blocks:
    \begin{enumerate}
        \item \label{int:Phi^k(A)} $\Phi_p^k(A)[0,j],\dots,\Phi_p^k(A)[p^k-1-j,p^k-1]$;
        \item \label{eqn:Z_kInt} $Z_k[p^k-j,0],\dots,Z_k[p^k-1,j-1]$;
        \item \label{int:Phi^k(iota(A))} $\Phi_p^k(\iota(A))[0,j],\dots,\Phi_p^k(\iota(A))[p^k-1-j,p^k-1]$;
        \item \label{eqn:W_kInt} $W_k[p^k-j,0],\dots,W_k[p^k-1,j-1]$.
    \end{enumerate}
\end{corollary}
\begin{remark}
    Note that when $j=0$, \eqref{eqn:Z_kInt} and \eqref{eqn:W_kInt} are not defined, since indeed the $0^{\nth}$ diagonal does not intersect $Z_k$ and $W_k$.
\end{remark}
\begin{remark}
\label{rem:Z_k/W_kInts}
Let $0\leq j\leq p^k-1$.
    \begin{enumerate}
            \item \label{D_jcapZ_k} By Corollary \ref{thm:IntZ_kW_k}(\ref{eqn:Z_kInt}), $\mathcal{D}_j(\Phi^{\infty}(A))\sqcap Z_k=
            \underbrace{0,\dots,0}_{j\text{ times}}$. 
            \item \label{D_jcapW_k} By Corollary \ref{thm:IntZ_kW_k}(\ref{eqn:W_kInt}), 
            $$\mathcal{D}_j(\Phi^{\infty}(A))\sqcap W_k=\{\rho^2(E),\underbrace{0,\dots,0}_{j-2\text{ times}},\rho(E)\}.$$
    \end{enumerate}
\end{remark}
\subsubsection{\textbf{The Diagonals of} $\Phi_p^{\infty}(A)$ \textbf{are Periodic}}\hfill\\
\noindent Through Corollary \ref{thm:IntZ_kW_k} and Remark \ref{rem:Z_k/W_kInts}, one obtains the following lower bound on $\mathcal{P}(j)$, which is illustrated below.
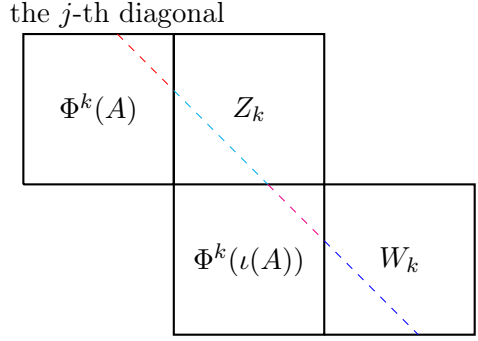
\begin{figure}[H]

\begin{tikzpicture}
    \draw[thick](-2,0)--(0,0)--(0,2)--(-2,2)--(-2,0);
    \node at (-1,1) {$\Phi^k(A)$};
    \draw[thick](0,0)--(2,0)--(2,2)--(0,2)--(0,0);
    \node at (1,1) {$Z_k$};
    \draw[thick](0,0)--(0,-2)--(2,-2)--(2,0);
    \node at (1,-1) {$\Phi^k(\iota(A))$};
    \draw[thick](2,-2)--(4,-2)--(4,0)--(2,0);
    \node at (3,-1) {$W_k$};
    \draw[dashed,red](-0.75,2)--(0,1.25);
    \draw[dashed,cyan](0,1.25)--(1.25,0);
    \draw[dashed,magenta](1.25,0)--(2,-0.75);
    \draw[dashed,blue](2,-0.75)--(3.25,-2);
    \node at (-0.75,2.25) {the $j$-th diagonal};
\end{tikzpicture}
\caption{The $j$-th diagonal is periodic, with its period being comprised of its journey through
four components, where the red part is its intersection with $\Phi^k(A)$, the cyan part is its intersection with $Z_k$, the pink part is its intersection with $\Phi^k(\iota(A))$, and the blue part is its intersection with $W_k$.}
\label{fig:j-thDiag}

\end{figure}
\begin{corollary}
\label{cor:DiagPerUppBnd}
    Let $p^k-1\leq j\leq p^{k+1}$. Then, $\mathcal{P}(j)\geq 2p^{k+1}$.
\end{corollary}
\begin{proof}   
     By Lemma \ref{lem:Phi^k(A)WindSize}, the largest zero portion in $\Phi_p^{k+1}(A)$ is a large $k$-pseudo-window, which has side lengths $p^k$. Since $j\geq p^k-1$, then, by Remark \ref{rem:Z_k/W_kInts}(\ref{D_jcapZ_k}), $\mathcal{D}_j(\Phi_p^{\infty}(A))\sqcap Z_{k+1}$ does not appear in $\mathcal{D}_j(\Phi_p^{\infty}(A))\sqcap \Phi_p^{k+1}(A)$, that is the cyan part of the diagonal in Figure \ref{fig:j-thDiag} did not appear in the red part of that diagonal. Hence, by Corollary \ref{thm:IntZ_kW_k}, $\mathcal{P}(j)\geq p^{k+1}$.\\
    
     \noindent It is now shown that $\mathcal{D}_j(\Phi_p^{\infty}(A))\sqcap W_{k+1}$ does not appear in 
     $$\left(\mathcal{D}_j(\Phi_p^{\infty}(A))\sqcap \Phi_p^{k+1}(A)\right)\bigsqcup\left(\mathcal{D}_j (\Phi_p^{\infty}(A))\sqcap Z_{k+1}\right)\bigsqcup\left(\mathcal{D}_j(\Phi_p^{\infty}(A))\sqcap \Phi_p^{k+1}(\iota(A))\right),$$ 
     that is, the blue part of the diagonal in Figure \ref{fig:j-thDiag} does not appear in the red, cyan, or pink parts of that diagonal. By Remark \ref{rem:Z_k/W_kInts}(\ref{D_jcapW_k}), the blue part of the diagonal in Figure \ref{fig:j-thDiag} is given by
     $$\mathcal{D}_j(\Phi_p^{\infty}(A))\sqcap W_{k+1}=\rho^2(E),\overbrace{0,\dots,0}_{j-2 \text{ times}},\rho(E).$$
     Firstly, if $j\geq p^k+3$, then, $j-2\geq p^k+1$. Therefore, by Lemma \ref{lem:Phi^k(A)WindSize}, $\mathcal{D}_j(\Phi_p^{\infty}(A))\sqcap W_{k+1}$ is not contained in $\mathcal{D}_j(\Phi_p^{\infty}(A))\sqcap \Phi_p^{k+1}(A)$ or in $\mathcal{D}_j(\Phi_p^{\infty}(A))\sqcap \Phi_p^{k+1}(\iota(A))$. Thus, the only other place the $j-2$ zeroes in $\mathcal{D}_j(\Phi_p^{\infty}(A))\sqcap W_{k+1}$ could appear is in $\mathcal{D}_j(\Phi_p^{\infty}(A))\sqcap Z_{k+1}$ (the cyan part of the diagonal in Figure \ref{fig:j-thDiag}). By Remark \ref{rem:Z_k/W_kInts}(\ref{D_jcapZ_k}), $\mathcal{D}_j(\Phi_p^{\infty}(A))\sqcap Z_{k+1}=\underbrace{0,\dots,0}_{j\text{ times}}$, and thus, by Lemma \ref{lem:LastColFirstRow}, $\mathcal{D}_j(\Phi_p^{\infty}(A))\sqcap W_{k+1}$ does not intersect $Z_{k+1}$. 
    \newline\newline
    \noindent Thus, the period of $\mathcal{D}_j(\Phi_p^{\infty}(A))$ consists of its intersections with $\Phi_p^{k+1}(A)$, $Z_{k+1}$, $\Phi_p^{k+1}(\iota(A))$, and $W_{k+1}$. Therefore, the period of the $j$-th diagonal is at least $2p^{k+1}$ when $j\geq p^k+3$.
    \newline\newline
    \noindent If $j=p^k+i$, where $i=1,2$, then, by Corollary \ref{thm:IntZ_kW_k}, $\mathcal{D}_j(\Phi_p^{\infty}(A))\sqcap Z_{k+1}$ (the cyan part in Figure \ref{fig:j-thDiag}) consists of exactly $j$ consecutive entries. Hence, if $D_j(\Phi_p^{\infty}(A))\sqcap W_{k+1}=\rho^2(E),\underbrace{0,\dots,0}_{j-2 \text{ times}},\rho(E)$ (the blue part of the diagonal in Figure \ref{fig:j-thDiag}) appeared in 
    $$\left(\mathcal{D}_j(\Phi_p^{\infty}(A))\sqcap \Phi_p^{k+1}(A)\right)\bigsqcup \left(\mathcal{D}_j(\Phi_p^{\infty}(A))\sqcap Z_{k+1}\right)\bigsqcup \left(\mathcal{D}_j(\Phi_p^{\infty}(A))\sqcap \Phi_p^{k+1}(\iota(A))\right),$$ then by Lemmas \ref{lem:row1Phi^k(A)} and \ref{lem:LastColFirstRow}, this configuration has to be contained in $\Phi_p^{k+1}(A)$ (the red part of the diagonal in Figure \ref{fig:j-thDiag}) or in $\Phi_p^{k+1}(\iota(A))$ (the pink part of the diagonal in Figure \ref{fig:j-thDiag}). 
\newline\newline
    \noindent By Lemma \ref{lem:Phi^k(A)WindSize} and the definition of $W_k$, the only zero windows of $\Phi_p^{k+1}(A)$ (respectively $\Phi_p^{k+1}(\iota(A))$) of size at least $p^k-1\leq j-2$ are contained in $Z_k$. On the other hand, by Lemma \ref{lem:LastColFirstRow}, a diagonal which enters $Z_k$ first goes through the last column of $\Phi_p^k(A)$ (respectively, the last row of $\Phi_p^{k-1}(\iota(A))$), does not consist of $\rho^2(E)$. Thus, the configuration $\rho^2(E),\underbrace{0,\dots,0}_{j-2\text{ times}},\rho(E)$ does not appear in $\Phi_p^{k+1}(A)$ (the red part of the diagonal in Figure \ref{fig:j-thDiag}) or in $\Phi_p^{k+1}(\iota(A))$ (the pink part in Figure \ref{fig:j-thDiag}). As a consequence, the period of the $j$-th diagonal consists of its intersections with $\Phi_p^{k+1}(A)$, $Z_{k+1}$, $\Phi_p^{k+1}(\iota(A))$, and $W_{k+1}$. Hence, the period of the $j$-th diagonal of $\Phi_p^{\infty}(A)$ is at least $2p^{k+1}$.
\end{proof}
\begin{lemma}
\label{lem:DiagPerLowBnd}
    Let $k\in \mathbb{N}$. For any $p^k+1\leq j\leq p^{k+1}$, $\mathcal{D}_j(\Phi_p^{\infty}(A))$ is periodic with $\mathcal{P}(j)\leq 2p^{k+1}$.
\end{lemma}
\begin{proof}
    The proof proceeds inductively on $k$. For $k=0$, this is a consequence of the definition of $\Phi_p(A)$, the definition of $\Phi_p(\iota(A))$, Theorem \ref{thm:PhiTrans}, and Corollary \ref{thm:IntZ_kW_k}.
    \newline\newline
    \noindent For the purposes of this proof, for a multiset $X=\{x_1,\dots,x_m\}\in \mathcal{A}^m$, define $$\Phi_p(X)=\{\Phi_p(x_1),\dots,\Phi_p(x_m)\},$$ so that $\Phi_p$ also acts on diagonals. Let $1\leq i\leq 2p^{k+1}-1$. Then, $\Phi_p$ maps $\Phi^\infty_p(A)[i,j]$ to the square
    $$\left\{\Phi_p^\infty(A)[\ell,k]:(p-1)i+1\leq \ell\leq i, (p-1)j+1\leq k\leq pj\right\}.$$
    Hence, $\Phi_p^\infty(A)[i,i+j]$ arises as part of the image under $\Phi_p$ of the $\Phi_p^\infty(A)\left[\left\lfloor\frac{i}{p}\right\rfloor,\left\lfloor\frac{i+j}{p}\right\rfloor\right]$-th entry. Moreover if $j=\sum_{\ell=0}^kj_{\ell}p^{\ell}$, then,
    $$\Phi_p^\infty(A)\left[\left\lfloor\frac{i}{p}\right\rfloor,\left\lfloor\frac{i+\sum_{\ell=0}^kj_{\ell}p^{\ell}}{p}\right\rfloor\right]=\Phi_p^\infty(A)\left[\left\lfloor\frac{i}{p}\right\rfloor,\sum_{\ell=1}^kj_{\ell}p^{\ell-1}+\left\lfloor\frac{i+j_0}{p}\right\rfloor\right].$$
     Thus, $\Phi_p^\infty(A)[i,i+j]$ arises from the  $\sum_{\ell=1}^kj_{\ell}p^{\ell-1}+\left\lfloor\frac{i+j_0}{p}\right\rfloor-\left\lfloor\frac{i}{p}\right\rfloor$-th diagonal. 
     \newline\newline
     \noindent Let $i=\sum_{\ell=0}^{k+1}i_{\ell}p^{\ell}$ be the base $p$ expansion of $i$. Then,
     \begin{equation}\begin{split}
         \left\lfloor\frac{i+j_0}{p}\right\rfloor-\left\lfloor\frac{i}{p}\right\rfloor=\left\lfloor\frac{\sum_{\ell=0}^{k+1}i_{\ell}p^{\ell}+j_0}{p}\right\rfloor-\left\lfloor\frac{\sum_{\ell=0}^{k+1}i_{\ell}p^{\ell-1}}{p}\right\rfloor\\
         =\left\lfloor\frac{i_0+j_0}{p}\right\rfloor-\left\lfloor\frac{i_0}{p}\right\rfloor=\left\lfloor\frac{i_0+j_0}{p}\right\rfloor=\begin{cases}
             0&i_0<p-j_0,\text{ or }j_0=0,\\
             1&p-j_0\leq i_0\leq p-1
         \end{cases}.
         \end{split}
     \end{equation}
     \noindent Let $p^k+1\leq j\leq p^{k+1}$. The proof now splits into two cases. The first case is $j_0=0$. In this case, one has $\mathcal{D}_j(\Phi_p^{\infty}(A))=\Phi_p\left(\mathcal{D}_{\frac{j}{p}}(\Phi_p^{\infty}(A))\right)$. Hence, by the induction hypothesis for $\mathcal{D}_{\frac{j}{p}}(\Phi^{\infty}(A))$ and the fact that $\Phi_p$ is a function from $\mathcal{A}$ to $\mathcal{A}^{p\times p}$, the $j$-th diagonal is periodic with $\mathcal{P}(j)\leq 2p^{k+1}$.
     \newline\newline
     \noindent If $j_0\neq 0$, for every $i=1,\dots,2p^{k+1}-1$, then, $\Phi_p^{\infty}(A)[i,i+j]$ arises from 
     $$\Phi_p\left\{\Phi_p^{\infty}(A)\left[\left\lfloor\frac{i}{p}\right\rfloor,\left\lfloor\frac{j}{p}\right\rfloor+\left\lfloor\frac{i}{p}\right\rfloor\right]:i_0\leq p-j_0\right\}\bigsqcup\Phi_p\left\{\Phi_p^{\infty}(A)\left[\left\lfloor\frac{i}{p}\right\rfloor,\left\lfloor\frac{j}{p}\right\rfloor+1+\left\lfloor\frac{i}{p}\right\rfloor\right]:p-j_0\leq i_0\leq p-1\right\}.$$ 
     This means that $\mathcal{D}_j(\Phi^{\infty}(A))$ arises as parts of the images of $\mathcal{D}_{\left\lfloor\frac{j}{p}\right\rfloor}(\Phi_p^{\infty}(A))$ and $\mathcal{D}_{\left\lfloor\frac{j}{p}\right\rfloor+1}(\Phi_p^{\infty}(A))$ under $\Phi_p$.
     \newline\newline
      \noindent By the induction hypothesis, $\mathcal{D}_{\left\lfloor\frac{j}{p}\right\rfloor}(\Phi_p^{\infty}(A))$ and the $\mathcal{D}_{\left\lfloor \frac{j}{p}\right\rfloor+1}(\Phi_p^{\infty}(A))$ are periodic with period length at most $2p^{k}$. Thus, their images under $\Phi_p$ are periodic with period length at most $2p^{k+1}$, so that the $\mathcal{D}_j(\Phi^{\infty}(A))$ is periodic with period length at most $2p^{k+1}$. 
\end{proof}
\begin{corollary}
\label{cor:PerSize}
    Let $k\in \mathbb{N}$ and let $p^k\leq j\leq p^{k+1}-1$. Then, $\mathcal{P}(j)=2p^{k+1}$.
\end{corollary}
\section{Escape of Mass Computations}\label{Sect: 7,0}
\noindent Theorem \ref{thm: morphism} and the results of the previous section are now used to prove Theorems \ref{cor:minEsc}, \ref{thm:MaxEsc1}, \ref{thm:GenEsc=1}, and \ref{thm:TMFullEsc}.
\subsubsection{\textbf{Escape of Mass of} $\mathcal{D}_j(\Phi^{\infty}(A))$}\hfill\\
\noindent Recall that by Corollary \ref{cor:EscWall}, escape of mass can be computed by using the number wall through equation \eqref{eqn:EscWall}. Furthermore, by Theorem \ref{thm: morphism}, $\Pi(\Phi_p^{\infty}(A))$ generates $\chi\left(W_p\left(\textbf{C}^{(p)}\right)\right)$. Hence, the escape of mass of the $j$-th diagonal of $W_p(\mathbf{C}^{(p)})$ can be computed through the similar notion of \textit{escape of mass of $\mathcal{D}_j(\Phi_p^{\infty}(A))$}.
\newline\newline 
\noindent For $j,N\in \mathbb{N}$, let $e_{j,N}$ be the proportion of mass greater than $N$ which escapes in the $j$-th diagonal in $\Phi_p^{\infty}(A)$. By Corollary \ref{cor:PerSize}, one can express $e_{j,N}$ explicitly by
\begin{equation}
    \label{eqn:escDef}
    e_{j,N}=\frac{\sum_{h_{i+1,j}\leq 2p^{k+1}}\max\{h_{i+1,j}-h_{i,j}-N,0\}}{2p^{k+1}},\end{equation}
where $h_{1,j}<h_{2,j}<\dots $ are the indices at which $\Phi^{\infty}(A)[h_{i,j},h_{i,j}+j]\neq 0$. For $j\leq p^{k+1}$, define $\phi_{k+1,j,N}$ to be the amount of mass of size at least $N$ which escapes in the $j$-th diagonal of $\Phi^{k+1}(A)$, that is
$$\phi_{k+1,j,N}=\sum_{h_{i+1,j}\leq p^{k+1}-1-j}\max\{h_{i+1,j}-h_{i,j}-N,0\}.$$
Note that the periodic nature of the diagonals implies that for every $j\in \mathbb{N}$, there exists $N$ for which $\phi_{k+1,j,N}=0$. Hence, one can take $N$ large enough so that $\phi_{k+1,0,N}=0$. \\

\noindent The following lemma is a direct consequence of Lemmas \ref{lem:Equivar}, \ref{lem:row1Phi^k(A)}, and \ref{lem:LastColFirstRow}. 
\begin{lemma}
\label{lem:EscIsIotaInv}
    Let $k\in \mathbb{N}$, let $p^{k}+1\leq j\leq p^{k+1}$, and let $\tilde{h}_{1,j}<\tilde{h}_{2,j}<\dots$ be the indices for which $\Phi^{k+1}(\iota(A))[\tilde{h}_{i,j},\tilde{h}_{i,j}+j]\neq 0$. Then, 
    \begin{equation}
    \label{eqn:e_j,N=e(iota(j))}
        \phi_{k+1,j,N}=\sum_{\tilde{h}_{i+1,j}\leq p^{k+1}-1-j}\max\{\tilde{h}_{i+1,j}-\tilde{h}_{i,j}-N,0\}+O_N(1).
    \end{equation}
\end{lemma}
\begin{proof}
    The proof proceeds by induction. For $k=-1$, Corollary \ref{cor:PerSize} implies that the first diagonal of $\iota(A)$ has the following path:
    \begin{enumerate}
        \item $\Phi(\iota(A))[0,0],\dots,\Phi(\iota(A))[p-1,p-1]$;
        \item $\Phi(A)[0,0],\dots,\Phi(A)[p-1,p-1]$.
    \end{enumerate}
    Note that by \eqref{eqn:A,iota(A)}, $\Phi(\iota(A))[i,i],\Phi(A)[j,j]\in \left\{A,\iota(A)\right\}$ for every $i,j\leq p-1$. Thus, the right hand side of \eqref{eqn:e_j,N=e(iota(j))} is equal to zero for every $N\in \mathbb{N}$ and $\Phi_{0,1,N}=0$ for every $N\in \mathbb{N}$. 

    \noindent Assume that \eqref{eqn:A,iota(A)} holds up to $k-1$, and \eqref{eqn:A,iota(A)} is now proved for $k$. Let $p^{k}+1\leq j=j_kp^k+\sum_{i=0}^{k-1}j_ip^i\leq p^{k+1}$ and let $\operatorname{suf}(j)=j-j_kp^k$. By Theorem \ref{thm:PhiTrans} Corollary \ref{cor:PerSize}\eqref{eqn:Phi(1bar)Trans}, the $j^\nth$ diagonal of $\Phi^{k+1}(\iota(A))$ consists of 
    \begin{enumerate}
        \item \label{iota(A)1} $\Phi^{k+1}(\iota(A))[0,j],\dots,\Phi^{k+1}(\iota(A))[p^{k+1}-1-j,p^{k+1}-1]$;
        \item \label{iota(A)2} $W_{k+1}[p^{k+1}-j,0],\dots,W_{k+1}[p^{k+1}-1,j-1]$;
        \item \label{iota(A)3}$\Phi^{k+1}(A)[0,j],\dots,\Phi^{k+1}(A)[p^{k+1}-1-j,p^{k+1}-1]$;
        \item \label{iota(A)4}$Z_{k+1}[p^{k+1}-j,0],\dots,Z_{k+1}[p^{k+1}-1,j-1]$.
    \end{enumerate}
    On the other hand, the $j^\nth$ diagonal of $\Phi^{k+1}(A)$ consists of \eqref{iota(A)3},\eqref{iota(A)4},\eqref{iota(A)1},\eqref{iota(A)2}. By Lemma \ref{lem:row1Phi^k(A)} and Lemma \ref{lem:LastColFirstRow}, reordering the steps \eqref{iota(A)1},\eqref{iota(A)2}, \eqref{iota(A)3},\eqref{iota(A)4} does not affect $\phi_{k+1,j,N}$. Therefore, by the induction hypothesis, on has \eqref{eqn:e_j,N=e(iota(j))}. 
\end{proof}
\begin{lemma}
\label{lem:e_jFormula}
Let $k\in \mathbb{N}$ and let $p^k+1
\leq j=\sum_{\ell=0}^kj_{\ell}p^{\ell}\leq p^{k+1}$. Then, for every sufficiently large $k$, 
\begin{equation}
  \label{eqn:e_jCalcn=1}
    e_{j,N}=\frac{j+\phi_{k+1,j,N}+O_N(1)}{p^{k+1}}\geq \frac{j_k}{p}+O_N(p^{-k}).
\end{equation}
\end{lemma}
\begin{proof}
    By Corollary \ref{cor:PerSize}, the period of the $j$-th diagonal is $2p^{k+1}$. Hence by Lemma \ref{lem:EscIsIotaInv}, for any large enough $k$, $e_{j,N}$ can be computed via Corollary \ref{thm:IntZ_kW_k} to obtain
    \begin{equation}
        e_{j,N}=\frac{2\phi_{k+1,j,N}}{2p^{k+1}}+\frac{2j+O_N(1)}{2p^{k+1}}=\frac{j+\phi_{k+1,j,N}+O_N(1)}{p^{k+1}}\geq \frac{j_k}{p}+O_N(p^{-k}).
    \end{equation}
\end{proof}
\begin{remark}
\label{rem:MorphVsCode}
    Due to the definition of $\Pi$ and Theorem \ref{thm: morphism}, the escape of mass of the $j^{\nth}$ diagonal of the number wall of $\Theta_p(t)$, denoted $e_{j,N}(\Theta_p(t))$ and defined in equation \eqref{eqn:EscDefWall}, is equal to $e_{j,N}$ for every pair $(j,N)$, such that $j$ is large enough with respect to $N$.
\end{remark}
\subsection{Computing $\phi_{k+1,j,N}$}
\begin{figure}[H]

\begin{tikzpicture}
    \node at (3.4,3.5) {$\Phi_p^{k+1}(A)$};
    \draw[thick](0,-4)--(0,3)--(7,3)--(7,-4)--(0,-4);
    \draw[blue](2,3)--(2,2)--(3,2)--(3,3);
    \node at (2.45,2.45) {$\Phi_p^k(A)$};
    \draw[blue](1,2)--(2,2)--(2,3)--(1,3)--(1,2);
    \node at (1.45,2.45) {$Z_k$};
    \draw[blue](3,2)--(3,1)--(4,1)--(4,2)--(3,2);
    \node at (3.45,1.45) {\begin{tiny}$\Phi_p^k(\iota(A))$\end{tiny}};
    \draw[blue](2,1)--(3,1)--(3,2)--(2,2)--(2,1);
    \node at (2.45,1.45) {$W_k$};
    \draw[blue](7,-2)--(6,-2)--(6,-1)--(7,-1);
    \node at (6.45,-1.45) {$\Phi_p^k(A)$};
    \draw[blue](6,-2)--(6,-3)--(7,-3);
    \node at (6.45,-2.45) {$W_k$};
    \draw[blue](6,-2)--(5,-2)--(5,-1)--(6,-1);
    \node at (5.45,-1.45) {$Z_k$};
    \draw[thick,dotted,magenta](1.45,3)--(7,-2.55);
    \draw[violet,thick,dotted](2.3,3)--(7,-1.7);
\end{tikzpicture}
\caption{The path of the $j$-th diagonal in $\Phi_p^{k+1}(A)$. The purple diagonal corresponds to even $j_k$, whereas the pink diagonal corresponds to odd $j_k$.}
\label{pic:DiagParity}
\end{figure}
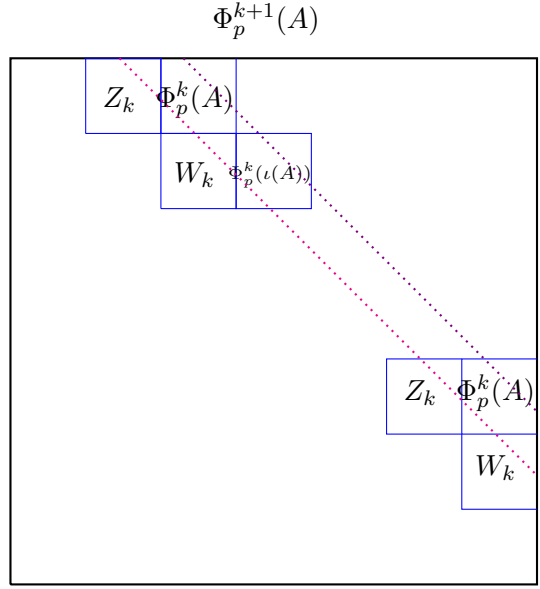

\hfill\\
\noindent Equation (\ref{eqn:e_jCalcn=1}) shows that to bound $e_{j,N}$, it suffices to bound $\phi_{k+1,j,N}$. Assume that $k$ is sufficiently large with respect to $N$, so that $\frac{O_N(1)}{p^{k+1}}$ is very small, and write $\phi_{k+1,j}$ instead of $\phi_{k+1,j,N}$ to ease on notations.
\begin{definition}\label{def: 7.4}
For $j=\sum_{i=0}^{k}j_ip^i$, where $0\leq j_i\leq p-1$ for $0\leq i\leq k$ and $j_k\neq0$, define the \textbf{suffix} of $j$ by
$$\operatorname{suf}(j):=\sum_{i=0}^{k-1}j_ip^i$$ and 
$$\alpha\left(\sum_{i=0}^kj_ip^i\right)=\sum_{i=0}^{k}\alpha_ip^i:=p^{k+1}-j.$$
Then, $\alpha_i=p-j_i$ for every $i=0,\dots,k$. Note that for every $j$, the orbit of $j$ under the operations $\operatorname{suf}$ and $\alpha$ is finite. 
\end{definition}
\begin{remark}
\label{rem:DiagExit}
    The motivation for defining $\alpha(j)$ is the fact that, for $j$ as in Definition \ref{def: 7.4}, $\mathcal{D}_j(\Phi_p^{\infty}(A))\sqcap \Phi_p^{k}(A)$ contains exactly $\alpha(j)$ entries. Moreover, $j_k$ is even if and only if $\alpha_k$ is odd. Hence, by Theorem \ref{thm:PhiTrans},
    \begin{enumerate}
        \item if $j_k$ is even, then, the $j$-th diagonal enters and leaves $\Phi_p^{k+1}(A)$ at the word $\Phi_p^{k}(A)\in \mathcal{A}^{p^{k}\times p^{k}}$;
        \item if $j_k$ is odd, then, the $j$-th diagonal enters $\Phi_p^{k+1}(A)$ at the $p^k\times p^k$ word $Z_k$ and leaves $\Phi_p^{k+1}(A)$ at the word $W_k\in \mathcal{A}^{p^k\times p^k}$.
    \end{enumerate}
    This is illustrated in Figure \ref{pic:DiagParity}
\end{remark}
\noindent The goal is now to write $\phi_{k+1,j}$ as a function of $\phi_{k+1,j'}$, where $j'$ is in the (finite) orbit of $j$ under the operators $\alpha$ and $\operatorname{suf}$. In particular, the following theorem is established.
\begin{theorem}
\label{thm:phi_k,jCalc}
Up to an additive constant of size at most $p$,
\begin{equation}
\label{eqn:phi_k+1,jCalc}
    \phi_{k+1,j}=\begin{cases}
        \phi_{k,\operatorname{suf}(j)}&j_k=p-1,\operatorname{suf}(j)\neq 0\\
        2\cdot (p^k-\operatorname{suf}(j))+\phi_{k,p^k-\operatorname{suf}(j)}&j_k=p-2,\operatorname{suf}(j)\neq 0\\
        (p-j_k-3)\cdot(p^k-\operatorname{suf}(j))+(p-j_k-2)\cdot\phi_{k,\operatorname{suf}(j)}&0\leq j_k\leq p-3 \text{ is even}, \operatorname{suf}(j)\neq 0\\
        (p-j_k)\cdot(p^k-\operatorname{suf}(j))+(p-j_k-1)\cdot\phi_{k,p^k-\operatorname{suf}(j)}&j_k\leq p-4\text{ is odd},\operatorname{suf}(j)\neq 0\\
        (p-j_k)\cdot p^k&j=j_kp^k, j_k\text{ is odd}\\
        0&j=j_kp^k,j_k\text{ is even.}
    \end{cases}
\end{equation}
If the expression in the right hand side of \eqref{eqn:phi_k+1,jCalc} is negative, then, $\phi_{k+1,j}=0$.
\end{theorem}
\begin{proof}
In order to prove Theorem \ref{thm:phi_k,jCalc}, it suffices to go over four different cases. 
\begin{enumerate}
\item[] \textbf{Case 1: $j_k=p-1$.}
\newline
If $j_k=p-1$, then by Theorem \ref{thm:PhiTrans}, the $j$-th diagonal intersects $\Phi_p^{k+1}(A)$ only at $$\Phi_p^k(A)[0,\operatorname{suf}(j)],\dots,\Phi_p^k(A)[p^k-\operatorname{suf}(j)-1,p^k-1].$$ Hence, $\phi_{k+1,j}=\phi_{k,\operatorname{suf}(j)}$.
\begin{figure}[H]
\begin{tikzpicture}
    \draw[thick](0,-2)--(0,2)--(4,2)--(4,-2)--(0,-2);
    \node at (2,-2.4) {$\Phi_p^{k+1}(A)$};
    \draw[blue](2.75,2)--(2.75,0.75)--(4.0,0.75);
    \node at (3.35,1.35) {$\Phi_p^k(A)$};
    \draw[thick,red,dotted](3.2,2)--(4,1.2);
\end{tikzpicture}
\caption{the path of the $j$-th diagonal when $j_k=p-1$}
\end{figure}
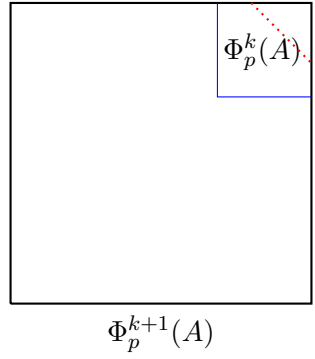

\item[] \textbf{Case 2: $j_k=p-2$.}

\noindent If $j_k=p-2$, then by Theorem \ref{thm:PhiTrans}, the $j$-th diagonal intersects $\Phi_p^{k+1}(A)$ at the following entries:
\begin{enumerate}
    \item $Z_k[0,\operatorname{suf}(j)],\dots,Z_k[p^k-\operatorname{suf}(j)-1,p^k-1]$;
    \item $\Phi_p^k(A)[p^k-\operatorname{suf}(j),0],\dots,\Phi_p^k(A)[p^k-1,\operatorname{suf}(j)-1]$;
    \item $W_k[0,\operatorname{suf}(j)],\dots,W_k[p^k-\operatorname{suf}(j)-1,p^k-1]$.
\end{enumerate}
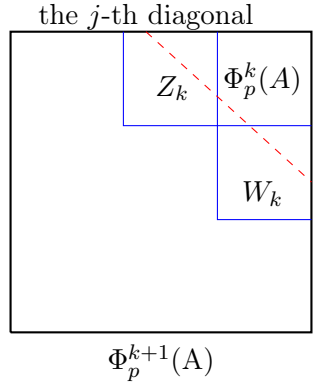
\begin{figure}[H]
\begin{tikzpicture}
    \draw[thick](0,-2)--(0,2)--(4,2)--(4,-2)--(0,-2);
    \node at (2,-2.4) {$\Phi_p^{k+1}$(A)};
    \draw[blue](2.75,2)--(2.75,0.75)--(4,0.75);
    \node at (3.35,1.35) {$\Phi_p^k(A)$};
    \draw[blue](1.5,2)--(1.5,0.75)--(2.75,0.75);
    \node at (2.15,1.25) {$Z_k$};
    \draw[blue](2.75,0.75)--(2.75,-0.5)--(4,-0.5);
    \node at (3.35,-0.15) {$W_k$};
    \draw[red,dashed] (1.8,2)--(4,0);
    \node at (1.8,2.2) {the $j$-th diagonal};
\end{tikzpicture}
\caption{the path of the $j$-th diagonal when $j_k=p-2$}
\end{figure}
Thus, $\phi_{k+1,j}=2(p^k-\operatorname{suf}(j))-2+\phi_{k,p^k-\operatorname{suf}(j)}$.\\ 

\item[] \textbf{Case 3: $2\leq j_k\leq p-3$ is even.}

\noindent If $j_k\leq p-3$ is even, the $j$-th diagonal begins by going over the following route $\left\lfloor\frac{p^{k+1}-j}{2p^k}\right\rfloor$ times:
\begin{enumerate}
    \hypertarget{route:Phi^k(A)}{\item} $\Phi_p^k(A)[0,\operatorname{suf}(j)],\dots,\Phi_p^k(A)[p^k-\operatorname{suf}(j)-1,p^k-1]$;
    \hypertarget{route:Z_k}{\item} $Z_k[p^k-\operatorname{suf}(j),0],\dots,Z_k[p^k-1,\operatorname{suf}(j)-1]$;
    \hypertarget{route:Phi^k(bar(A))}{\item} $\Phi_p^k(\iota(A))[0,\operatorname{suf}(j)],\dots,\Phi_p^k(\iota(A))[p^k-\operatorname{suf}(j)-1,p^k-1]$;
    \hypertarget{route:W_k}{\item}  $W_k[p^k-\operatorname{suf}(j),0],\dots,W_k[p^k-1,\operatorname{suf}(j)-1]$.
\end{enumerate}
Furthermore,
\begin{align*}\left\lfloor\frac{p^{k+1}-j}{2p^k}\right\rfloor&=\left\lfloor\frac{p^{k+1}-\sum_{i=0}^kj_ip^i}{2p^k}\right\rfloor\\&=\left\lfloor\frac{p-j_k}{2}-\frac{\sum_{i=0}^{k-1}j_ip^i}{2p^k}\right\rfloor\\&=\protect{\begin{cases}
    \frac{p-j_k-1}{2},&\operatorname{suf}(j)=0\\
    \frac{p-j_k-3}{2},&\text{else}
\end{cases}}.\end{align*}
By Remark \ref{rem:DiagExit}, after going through the route
\hyperlink{route:Phi^k(A)}{(a)}, \hyperlink{route:Z_k}{(b)}, \hyperlink{route:Phi^k(bar(A))}{(c)}, \hyperlink{route:W_k}{(d)} $\left\lfloor\frac{p^{k+1}-j}{2p^k}\right\rfloor$ times, the $j^\nth$ diagonal enters $\Phi_p^k(A)[0,\operatorname{suf}(j)],\dots,\Phi_p^k(A)[p^k-\operatorname{suf}(j)-1,p^k-1]$.
\begin{figure}[H]
\begin{tikzpicture}
    \node at (-0.5,2.5) {$\Phi_p^{k+1}(A)$};
    \draw[thick](-5,2)--(-5,-5)--(2,-5)--(2,2)--(-5,2);
    \draw[blue](-3,1)--(-2,1)--(-2,2)--(-3,2)--(-3,1);
    \node at (-2.45,1.45) {$\Phi_p^k(A)$};
    \draw[blue](-2,1)--(-1,1)--(-1,2)--(-2,2)--(-2,1);
    \node at (-1.5,1.5){$Z_k$};
    \draw[blue](-2,1)--(-2,0)--(-1,0)--(-1,1);
    \node at (-1.45,0.45) {\begin{tiny}$\Phi_p^k(\iota(A))$\end{tiny}};
    \draw[blue](-1,1)--(0,1)--(0,0)--(-1,0)--(-1,1);
    \node at (-0.5,0.5) {$W_k$};
    \draw[blue](2,-3)--(1,-3)--(1,-2)--(2,-2)--(2,-3);
    \node at (1.45,-2.45) {$\Phi_p^k(A)$};
    \draw[violet,thick,dotted](-3,2)--(2,-3);
    \draw[purple,thick,dotted](-2.5,2)--(2,-2.5);
\end{tikzpicture}
\caption{The purple line denotes Case \protect\hyperlink{case:j=j_kp^kEven}{3.1} and the burgundy line is Case \protect\hyperlink{case:j_kEven}{3.2}, both defined below.}
\end{figure}
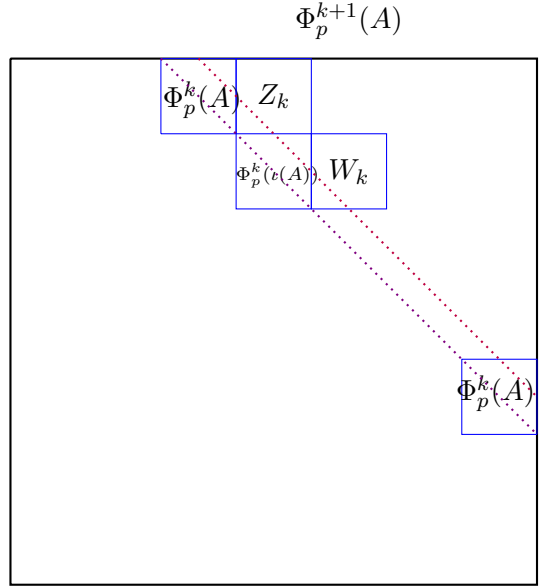
\begin{enumerate}
\item\hypertarget{case:j=j_kp^kEven}{\textbf{Case 3.1: $2\leq j_k\leq p-3$ is even and $\operatorname{suf}(j_k)=0$.}} 
\noindent In this case, $j=j_kp^k$. Thus, the $j_kp^k$-th diagonal does not intersect $Z_k$ and $W_k$ at all. That is, it enters $\Phi_p^k(A)$ as soon as it leaves $\Phi_p^k(\iota(A))$ (and vice versa). Hence, $\phi_{k+1,j_kp^k}=2\phi_{k,0}=0$. 

\item\hypertarget{case:j_kEven} {\textbf{Case 3.2: $2\leq j_k\leq p-3$ is even and $\operatorname{suf}(j_k)\neq 0$.}}
\noindent By Remark \ref{rem:DiagExit} and the steps \hyperlink{route:Phi^k(A)}{(a)}, \hyperlink{route:Z_k}{(b)}, \hyperlink{route:Phi^k(bar(A))}{(c)}, \hyperlink{route:W_k}{(d)}, 
\begin{equation}
\begin{split}
    \phi_{k+1,j}=\frac{p-j_k-3}{2}\left(2\cdot (p^k-\operatorname{suf}(j))+2\phi_{k,\operatorname{suf}(j)}\right)+\phi_{k,\operatorname{suf}(j)}\\
    =(p-j_k-3)\left(p^k-\operatorname{suf}(j)+\phi_{k,\operatorname{suf}(j)}\right)+\phi_{k,\operatorname{suf}(j)}\\
    =(p-j_k-3)\cdot (p^k-\operatorname{suf}(j))+(p-j_k-2)\phi_{k,\operatorname{suf}(j)}.
\end{split}
\end{equation}
\end{enumerate}
\item[] \textbf{Case 4: $j_k\leq p-4$ is odd.}

\noindent If $j_k\leq p-4$ is odd, then, the $j$-th diagonal intersects $\Phi_p^{k+1}(A)$ in the following route which is repeated $\left\lfloor\frac{p^{k+1}-j}{2p^k}\right\rfloor$ times:
\begin{enumerate}
    \item \label{route:Z_kodd} $Z_k[0,\operatorname{suf}(j)],\dots,Z_k[p^k-\operatorname{suf}(j)-1,p^k-1]$;
    \item \label{route:Phi^k(iot(A))Odd} $\Phi_p^k(\iota(A))[p^k-\operatorname{suf}(j),0],\dots,\Phi_p^k(A)[p^k-1,\operatorname{suf}(j)-1]$;
    \item \label{route:W_kodd} $W_k[0,\operatorname{suf}(j)],\dots,W_k[p^k-\operatorname{suf}(j)-1,p^k-1]$;
    \item \label{route:Phi^k(A)Odd} $\Phi_p^k(A)[p^k-\operatorname{suf}(j),0],\dots,\Phi_p^k(A)[p^k-1,\operatorname{suf}(j)-1]$.
\end{enumerate}
Since $\frac{p^{k+1}-j}{2p^k}=\frac{p-j_k}{2}-\frac{\operatorname{suf}(j)}{2p^k}$, one has that
$$\left\lfloor\frac{p^{k+1}-j}{2p^k}\right\rfloor=\begin{cases}
    \frac{p-j_k}{2}&\operatorname{suf}(j)=0\\
    \frac{p-j_k-2}{2}&\text{else}
\end{cases}.$$ 
By Remark \ref{rem:DiagExit}, after repeating \hyperlink{route:Z_kodd}{(a)}, \hyperlink{route:Phi^k(iot(A))Odd}{(b)}, \hyperlink{route:W_kodd}{(c)}, \hyperlink{route:Phi^k(A)Odd}{(d)} $\left\lfloor\frac{p^{k+1}-j}{2p^k}\right\rfloor$ times, then, the $j$-th diagonal goes through steps \hyperlink{route:Z_kodd}{(a)}, \hyperlink{route:Phi^k(iot(A))Odd}{(b)}, and \hyperlink{route:W_kodd}{(c)}. 
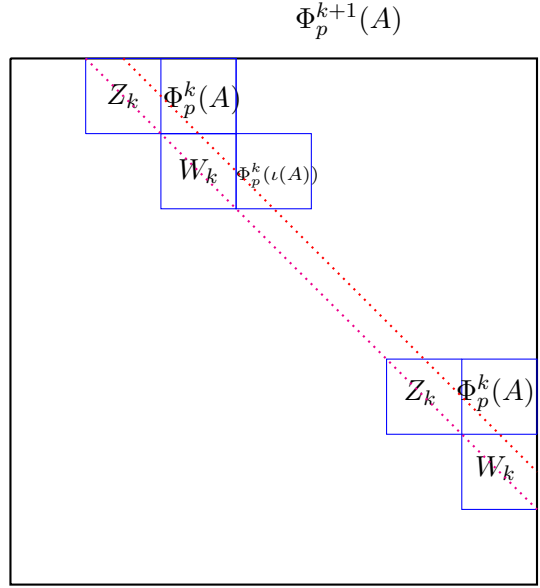
\begin{figure}[H]
\begin{tikzpicture}
    \node at (-0.5,2.5) {$\Phi_p^{k+1}(A)$};
    \draw[thick](-5,2)--(-5,-5)--(2,-5)--(2,2)--(-5,2);
    \draw[blue](-3,1)--(-2,1)--(-2,2)--(-3,2)--(-3,1);
    \node at (-2.45,1.45) {$\Phi_p^k(A)$};
    \draw[blue](-4,1)--(-2,1)--(-2,2)--(-4,2)--(-4,1);
    \node at (-3.5,1.5){$Z_k$};
    \draw[blue](-2,1)--(-2,0)--(-1,0)--(-1,1)--(-2,1);
    \node at (-1.45,0.45) {\begin{tiny}$\Phi_p^k(\iota(A))$\end{tiny}};
    \draw[blue](-3,1)--(-2,1)--(-2,0)--(-3,0)--(-3,1);
    \node at (-2.5,0.5) {$W_k$};
    \draw[blue](2,-3)--(1,-3)--(1,-2)--(2,-2)--(2,-3);
    \node at (1.45,-2.45) {$\Phi_p^k(A)$};
    \draw[blue](1,-3)--(1,-4)--(2,-4);
    \node at (1.45,-3.45) {$W_k$};
    \draw[blue](1,-3)--(0,-3)--(0,-2)--(1,-2);
    \node at (0.45,-2.45) {$Z_k$};
    \draw[magenta,thick,dotted](-4,2)--(2,-4);
    \draw[red,thick,dotted](-3.5,2)--(2,-3.5);
\end{tikzpicture}
\caption{The pink line is Case \eqref{case:j_k=j_kp^kOdd} and the red line is Case \eqref{case:j_kOddOther}.}
\end{figure}
\begin{enumerate}\item\label{case:j_k=j_kp^kOdd}\textbf{Case 4.1: $j_k\leq p-4$ is odd and $\operatorname{suf}(j)=0$.} 

\noindent In this case, $j=j_kp^k$, so that the $j$-th diagonal does not intersect $\Phi_p^k(A)$ or $\Phi_p^k(\iota(A))$. Hence,
\begin{equation}\nonumber
\begin{split}
    \phi_{k+1,j}=\frac{p-j_k}{2}\cdot 2p^k=(p-j_k)p^k
\end{split}
\end{equation}
\item\label{case:j_kOddOther}\textbf{Case 4.2: $j_k\leq p-4$ is odd and $\operatorname{suf}(j)\neq 0$.}
In this case, 
\begin{equation}\nonumber    
\begin{split}
    \phi_{k+1,j}=\frac{p-j_k-2}{2}\left(2(p^k-\operatorname{suf}(j))+2\phi_{k,p^k-\operatorname{suf}(j)}\right)+2(p^k-\operatorname{suf}(j))+\phi_{k,p^k-\operatorname{suf}(j)}\\
    =(p-j_k)(p^k-\operatorname{suf}(j))+(p-j_k-1)\phi_{k,p^k-\operatorname{suf}(j)}.
\end{split}
\end{equation}
\end{enumerate}
\end{enumerate}
\end{proof}
\subsection{Escape of Mass and and Maximal Escape of Mass}
\label{sec:EscComp}\hfill\\
The goal of this section is to deduce Theorems \ref{cor:minEsc} and \ref{thm:MaxEsc1} from Theorem \ref{thm: morphism}.
\subsubsection{\textbf{Maximal Escape of Mass}}\hfill\\
The proof of Theorem \ref{thm:MaxEsc1} is short, and so it is completed first. Indeed, it suffices to find a sequence of diagonals which exhibits full escape of mass. 
\begin{lemma}
\label{lem:fullEsc}
    Let $k\in \mathbb{N}$ and let $1\leq j_k\leq p-2$ be an odd integer. Then, the sequence $\{j_kp^k\}_{k\in \mathbb{N}}$ exhibits full escape of mass, that is $$\lim_{N\rightarrow \infty}\lim_{k\rightarrow \infty}e_{j_kp^k,N}(\Theta)=1.$$
\end{lemma}
\begin{proof}
    By Lemma \ref{lem:e_jFormula} and Theorem \ref{thm:phi_k,jCalc}, for every $k$ sufficiently large with respect to $N$, 
    \begin{equation}\nonumber
        e_{j_kp^k,N}=\frac{(p-j_k)p^k+j_kp^k-1+O_N(1)}{p^{k+1}}=\frac{p^{k+1}+O_N(1)}{p^{k+1}}.
    \end{equation}
    Hence for every $N>1$, $e_{j_kp^k,N}\rightarrow 1$ as $k\rightarrow \infty$. Thus, by Theorem \ref{thm: morphism} and Remark \ref{rem:MorphVsCode},
    $$\lim_{N\rightarrow \infty}\lim_{k\rightarrow \infty}e_{j_kp^k,N}(\Theta)=\lim_{N\rightarrow \infty}\lim_{k\rightarrow \infty}e_{j_kp^k,N}=1.$$
\end{proof}
\subsubsection{\textbf{Bounding the Escape of Mass}}
\begin{lemma}
\label{lem:MinEsc}
    For every $k\in\mathbb{N}$ and $p^k\leq j<p^{k+1}$,  one has $e_{j,N}\geq \frac{2}{p}+O_N(p^{-k})$. 
\end{lemma}
\begin{proof}
By Lemma \ref{lem:e_jFormula}, if $j_k\geq 2$, then, 
\begin{equation}\nonumber
    e_{j,N}=\frac{j+\phi_{k+1,j}+O_N(1)}{p^{k+1}}\geq \frac{j_kp^k}{p^{k+1}}+O_N(p^{-(k+1)})\geq \frac{2}{p}+O_N(p^{-(k+1)}).
\end{equation}
Thus, by Remark \ref{rem:MorphVsCode}, for every sequence $j^{(n)}=\sum_{i=0}^{k_n}j_i^{(n)}p^i$ with $j^{(n)}_{k_n}\geq 2$ and $k_n\rightarrow \infty$, 
$$\lim_{N\rightarrow \infty}\liminf_{n\rightarrow \infty}e_{j^{(n)},N}(\Theta)=\lim_{N\rightarrow \infty}\liminf_{n\rightarrow \infty}e_{j^{(n)},N}\geq \frac{2}{p}.$$
Hence, by Lemma \ref{lem:e_jFormula}, Lemma \ref{lem:fullEsc}, and Remark \ref{rem:MorphVsCode}, it suffices to prove that $e_{j,N}\geq \frac{2}{p}+O_N(p^{-k})$ when $j_k=1$ and $\operatorname{suf}(j)\neq 0$. In this case, Theorem \ref{thm:phi_k,jCalc},
\begin{equation}\nonumber
\begin{split}
    \phi_{k+1,j,N}=(p-1)(p^k-\operatorname{suf}(j))+(p-2)\phi_{k,\alpha(\operatorname{suf}(j))}+O(1)\\
    =p^{k+1}-p^k-(p-1)\operatorname{suf}(j)+(p-2)\phi_{k,\alpha(\operatorname{suf}(j))}+O(1).
\end{split}
\end{equation}
Therefore, by Lemma \ref{lem:e_jFormula} and Remark \ref{rem:MorphVsCode},
\begin{align}
    e_{j,N}(\Theta)&=e_{j,N}=\frac{j+\phi_{k+1,j}+O_N(1)}{p^{k+1}}\nonumber\\
    &=\frac{p^k+\operatorname{suf}(j)+p^{k+1}-p^k-(p-1)\operatorname{suf}(j)+(p-2)\phi_{k,p^k-\operatorname{suf}(j)}+O_N(1)}{p^{k+1}}  \nonumber\\
    &=1-\frac{(p-2)\operatorname{suf}(j)-(p-2)\phi_{k,p^{k}-\operatorname{suf}(j)}}{p^{k+1}}+O_N\left(\frac{1}{p^k}\right).\label{eqn:e_jCalcCasej_k=1}
\end{align}
Since $\operatorname{suf}(j)\leq p^k$, then, $\operatorname{suf}(j)-\phi_{k,p^k-\operatorname{suf}(j)}\leq p^k$. Thus, the right hand side of (\ref{eqn:e_jCalcCasej_k=1}) is at greater than or equal to
$$1-\frac{p-2}{p^{k+1}}p^k-O_N(p^{-k})=\frac{2}{p}+O_N(p^{-k}).$$
As a consequence, for any sequence $j^{(n)}=\sum_{i=0}^{k_n}j^{(n)}_ip^k$ with $j^{(n)}_{k_n}=1$, one has
$$\lim_{N\rightarrow \infty}\liminf_{n\rightarrow \infty}e_{j^{(n)},N}\geq \lim_{N\rightarrow \infty}\lim_{n\rightarrow \infty}\left(\frac{2}{p}+O_N(p^{-k_n})\right)=\frac{2}{p}.$$
\end{proof}
\noindent To conclude the proof of Theorem \ref{cor:minEsc}, it suffices to find a sequence which exhibits exactly $\frac{2}{p}$ escape of mass. 
\begin{lemma}
\label{lem:minEsc}
    The $\{2p^k\}_{k\in \mathbb{N}}$ diagonals exhibit $\frac{2}{p}$ escape of mass. 
\end{lemma}
\begin{proof}
    By Lemma \ref{lem:e_jFormula} and Remark \ref{rem:MorphVsCode}, for every $N>1$, and for every sufficiently large $k$, 
    \begin{equation}
        e_{2p^k,N}=\frac{2}{p}+\frac{\phi_{k+1,2p^k}}{p^{k+1}}-\frac{O_N(1)}{p^{k+1}}.\nonumber
    \end{equation}
    By Theorem \ref{thm:phi_k,jCalc}, $\phi_{k+1,2p^k,N}=0$. Therefore, 
    \begin{equation}\nonumber
        e_{2p^k,N}=\frac{2}{p}+O_N\left(\frac{1}{p^{k+1}}\right)\rightarrow \frac{2}{p},
    \end{equation}
    as $k\rightarrow \infty$. Hence, by Lemma \ref{lem:MinEsc}, the $\{2p^k\}_{k\in \mathbb{N}}$ diagonals exhibit exactly $\frac{2}{p}$ escape of mass. 
\end{proof} \newpage
\subsection{Generic Escape of Mass}
\subsubsection{\textbf{A General Method to Establish Generic Escape of Mass}}\hfill\\
\noindent The following terminology is introduced to provide a method of proving generic escape of mass claims. This methodology is sequence agnostic, but in this paper it is applied to the $p$-Cantor and Thue-Morse sequences, the number wall of the latter being constructed in \cite{AN}.
\begin{definition}
\label{def:TileEquiv}
    Let $\boldsymbol{X}$ be a finite alphabet and let $\zeta:\boldsymbol{X_2}^\Omega\rightarrow \boldsymbol{X}_2^{\Omega}$ be a $[p,p]$-morphism. For $a\in \boldsymbol{X}$ and a square block $B\in \boldsymbol{X}^{p^k\times p^k}$, let $N_a(B)$ be the number of appearances of $a$ in the $B$. 
\end{definition}
\noindent This function is central to the following Fubini like claim which is used to prove that the $p$-Cantor and Thue Morse sequences exhibit full generic escape of mass.
\begin{theorem}
\label{thm:Til-->Esc}
    Let $p$ be a prime power, let $m\in \mathbb{N}$, and let $\Theta(t)\in \mathbb{F}_p(\!(t^{-1})\!)$ be a quadratic irrational. Let $\boldsymbol{X}$ be an alphabet containing the symbols $0$ and $A$, let $\zeta:\boldsymbol{X}^\Omega_2\rightarrow \boldsymbol{X}_2^{\Omega}$ be a $[p,p]$-morphism such that $A$ is $\zeta$-prolongable, and let $\xi:\boldsymbol{X}_2^\Omega\rightarrow \{0,1\}_2^{\Omega}$ be an $[m,m]$-coding. Assume that $\xi(0)=\boldsymbol{0}_{m\times m}$, $\zeta(0)=\boldsymbol{0}_{p\times p}$, and that the profile of the number wall of $\Theta(t)$ is equal to $\xi(\zeta^{\infty}(A))$. If
    \begin{equation}
    \label{eqn:N_0/p^2k}
        \lim_{k\rightarrow \infty}\frac{N_0(\zeta^k(A))}{p^{2k}}=1,
    \end{equation}
    then, for every $\varepsilon>0$, 
    $$\lim_{n\rightarrow \infty}d(k\in \mathbb{N}:e_{k,n}(\Theta(t))>1-\varepsilon)=1.$$
\end{theorem}
\begin{remark}
    For the $p$-Cantor sequence, $m=1$ and $\xi=\Pi$.
\end{remark}
\begin{proof}[Proof of Theorem \ref{thm:Til-->Esc}]
    For $k\in \mathbb{N}\cup \{\infty\}$, define 
    $$B_k=\bigcup_{\begin{matrix}
    (i,j)\in \mathbb{Z}^2:\\0\leq i,j\leq p^{k-1},\\ \zeta^k(A)[i,j]=0\end{matrix}}\left[ip^{-k},(i+1)p^{-k}\right]\times \left[jp^{-k},(j+1)p^{-k}\right].$$
    Then, 
    \begin{equation}
    \label{eqn:lambda_2(B_k)}
    \lambda_2(B_k)=\frac{N_0(\zeta^k(A))}{p^{2k}},\end{equation}
    where $\lambda_2$ is the Lebesgue measure on $[0,1]\times [0,1]$. Since $\zeta^{\infty}(A)$ exists, then, \begin{equation}B_{\infty}:=\lim_{k\rightarrow \infty}B_k\label{eqn: B_infty}\end{equation} is well defined. Note that $B_k$ are simple measurable sets approximating $B_{\infty}$. Hence, by \eqref{eqn:N_0/p^2k} and \eqref{eqn:lambda_2(B_k)}, $\lambda_2(B_{\infty})=\lim_{k\rightarrow \infty}\lambda_2(B_k)=1$.\\

\noindent First, to account for the limit on $n\rightarrow \infty$ in the definition of escape of mass, a variation of $\xi$ is defined. For $k\in \mathbb{N}\cup \{\infty\}$ and $n\in \mathbb{N}\cup \{0\}$, define $\xi_n$ in the following manner: For every diagonal $0\leq j\leq p^k-1$, and every $2\leq i\leq p^k-j-1$, define 
    $$\xi_n(\zeta^k(A))[\ell+h_{j,i},h_{j,i}+\ell+j]=\begin{cases}
        1&\ell=1,\dots,\min\{n,h_{j,i+1}-h_{j,i}\}\\
        \xi(\zeta^k(A))[h_{j,i}+\ell,h_{j,i}+\ell+j]&\text{else}
    \end{cases}.$$
\noindent That is, $\xi_n$ is identical to $\xi$, but where every string of at least $n$ zeroes on a diagonal has had the first $n$ zeroes replaced with ones.\\

\noindent Additionally, define
    $$B_{k,n}=\bigcup_{\begin{matrix}
        (i,j)\in \mathbb{Z}^2\\
        0\leq i,j\leq mp^k-1\\
        \xi_n(\zeta^k(A))[i,j]=0
    \end{matrix}}\left[ip^{-k}m^{-1},(i+1)p^{-k}m^{-1}\right]\times \left[jp^{-k}m^{-1},(j+1)p^{-k}m^{-1}\right].$$

\noindent Since $\lim_{k\rightarrow \infty}\zeta^k(A)$ exists, then, the limit $B_{\infty,n}:=\lim_{k\rightarrow \infty}B_{k,n}$ is well defined. Moreover, $B_{k,n}$ are simple measurable sets approximating $B_{\infty,n}$. Therefore, for every $n\in \mathbb{N}\cup \{0\}$ and for every $k\in \mathbb{N}\cup \{\infty\}$,
$$\lim_{k\rightarrow \infty}\frac{N_{0}(\xi_n(\zeta^k(A)))}{p^{2k}m^2}=\lim_{k\rightarrow \infty}\lambda(B_{k,n})=\lambda(B_{\infty,n}).$$ 
Note that the $j^\nth$ diagonal of $\xi_n(\zeta^k(A))$ differs from the $j^\nth$ diagonal of $\xi(\zeta^k(A))$ in at most $n\ell_k$ entries, where $\ell_k$ is the length of the period of the $k$-th diagonal. By \cite[Theorem 10]{KPS} and the analysis of \cite{P,PS}, there exists $C_{\Theta(t)}>0$ such that $\ell_k\leq C_{\Theta(t)}k$. Hence, 
\begin{equation}\nonumber
\begin{split}
\label{eqn:propXi(zeta^k(A))vsN}
    \frac{N_{0}(\xi(\zeta^k(A)))}{p^{2k}{m^2}}\geq \frac{N_{0}(\xi_n(\zeta^k(A)))}{p^{2k}m^2}\geq \frac{N_{0}(\xi(\zeta^k(A)))-nC_{\Theta(t)}kp^k}{p^{2k}m^2}\\
    =\frac{N_{0}(\xi(\zeta^k(A)))}{p^{2k}m^2}-O_{m,\Theta(t),n}(kp^{-k}).
\end{split}
\end{equation}
Consequently,
\begin{equation}
\label{eqn:lambda(B_infty,0)Ineq}
\lambda_2(B_{\infty,0})=\lim_{k\rightarrow \infty}\frac{N_0(\xi(\zeta^k(A)))}{p^{2k}m^2}=\lim_{k\rightarrow \infty}\frac{N_{0}(\xi_n(\zeta^k(A)))}{p^{2k}m^2}=\lambda_2(B_{\infty,n}),\end{equation}
for every $n\in \mathbb{N}$. Since $\xi(0)=\mathbf{0}_{m\times m}$, then, 
\begin{equation}\nonumber
    \frac{N_{0}(\xi(\zeta^k(A)))}{p^{2k}m^2}\geq \frac{N_{0}(\xi^k(A))m^2}{p^{2k}m^2}=\frac{N_0(\xi^k(A))}{p^{2k}}.
\end{equation}
By passing to a limit, 
\begin{equation}
\label{eqn:lambda(B_infty,0)}
    \lambda_2(B_{\infty,0})\geq \lim_{k\rightarrow \infty}\frac{N_0(\zeta^k(A))}{p^{2k}}=\lambda_2(B_\infty)=1,
\end{equation}
so that $\lambda_2(B_{\infty,0})=1$. Thus, by \eqref{eqn:lambda(B_infty,0)Ineq} and \eqref{eqn:lambda(B_infty,0)}, for every $n\in \mathbb{N}$, $\lambda_2(B_{\infty,n})=1$. \\

\noindent Let $\lambda_1$ be the one dimensional Lebesgue measure in $\mathbb{R}^2$ and for $a\in [0,1]$, define
$$\mathcal{D}_a=\left\{(x,1-x-a):a\leq x\leq 1-a\right\}.$$
By Fubini's Theorem, for every $n\in \mathbb{N}$ and for $\lambda_1$ almost every $a\in [0,1]$, one has $\lambda_1(B_{\infty,n}\cap \mathcal{D}_a)=1$.\\
\noindent The fact that $\lambda_1(B_{\infty,n}\cap \mathcal{D}_a)=1$ is now translated to generic escape of mass. Let $$s_{j,1}^{(k)}(n)<s_{j,2}^{(k)}(n)<\dots<s_{j,b_{k,n,j}}^{(k)}(n)$$ be the indices satisfying $\xi_N(\zeta^k(A))\left[s_{j,i}^{(k)}(n),j+s_{j,i}^{(k)}(n)\right]=1$. 
By the definition of escape of mass and Lebesgue measure, for any sequence $j_k$ with $j_k\leq p^k$ and $\lim_{k\rightarrow \infty}j_k$, 
    \begin{equation}
    \label{eqn:e_jFunctofEsc}
    \lim_{k\rightarrow \infty}e_{j_k,n}=\lim_{k\rightarrow \infty}\frac{\sum_{i=1}^{b_{k,n,j_k}}\left((s_{j_k,i}^{(k)}(n)-s_{j_k,i-1}^{(k)}(n)\right)}{m^2p^k}=\lim_{k\rightarrow \infty}\lambda_1(B_{k,n}\cap \mathcal{D}_{k,j_k}),\end{equation}
    where $e_{j,n}$ is defined in \eqref{eqn:escDef} and 
    $$\mathcal{D}_{k,j}=\left\{\left(\frac{i}{p^k},\frac{i+j}{p^k}\right):0\leq i\leq p^k-j-1\right\}.$$
    Thus, by \eqref{eqn:e_jFunctofEsc} for $\lambda_1$ almost every $a\in [0,1]$, and for every increasing sequence $j_k$ with $\frac{j_k}{p^k}\rightarrow a$, 
    \begin{equation}
        \lim_{k\rightarrow \infty}e_{j_k,n}=\lim_{k\rightarrow \infty}\lambda_1(B_{k,n}\cap \mathcal{D}_{k,j_k})=\lambda_1(B_{\infty,n}\cap \mathcal{D}_a)=1.
    \end{equation}
Since every $j\in \mathbb{N}$ is part of some sequence $j_k\leq p^k$ with $\lim_{k\rightarrow \infty}\frac{j_k}{p^k}=a$, then, for every $n\in \mathbb{N}$ and for every $\varepsilon>0$, one has
$$\lim_{k\rightarrow \infty}\frac{\left|\{j\leq k:e_{j,n}>1-\varepsilon\}\right|}{K}=\lambda_1\left\{a\in [0,1]:\lambda_1(B_{\infty,n}\cap \mathcal{D}_a)=1\right\}=1.$$

\end{proof}
\noindent Therefore, to prove Theorems \ref{thm:GenEsc=1} and \ref{thm:TMFullEsc}, it suffices to prove that $\Phi_p$ (respectively the morphism from \cite{AN}) satisfies the conditions of Theorem \ref{thm:Til-->Esc}. 
\begin{proof}[Proof of Theorem \ref{thm:GenEsc=1}]
By \cite[Theorem 1.6]{AR}, the set $B_\infty$ (see equation (\ref{eqn: B_infty}) defined by $\Phi_p$ has Hausdorff dimension less than 1. Therefore, $\lim_{k\rightarrow \infty}\frac{N_0(\Phi_p^k(A))}{p^{2k}}=1$, and hence, by Theorem \ref{thm:Til-->Esc}, Theorem \ref{thm:GenEsc=1} follows. 
\end{proof}
\subsubsection{\textbf{Proof of Theorem} \ref{thm:TMFullEsc}}\hfill\\
\noindent Theorem \ref{thm:TMFullEsc} is now established using a more general method that can be applied to any sequence whose number wall is fully understood. In particular, one requires the morphism from \cite{AN} that generates the number wall of the Thue-Morse sequence over $\mathbb{F}_2$.
\begin{definition}
    Let $\mathbf{X}=\{x_0,x_1,\dots,x_n\}$ be a finite alphabet. With the symbols of Definition \ref{def:TileEquiv}, define the \textbf{transition matrix} of $\zeta$ to be the matrix $M\in \boldsymbol{X}^{n\times n}$ by
    $$M_{i,j}=\frac{N_{x_{i-1}}(\zeta(x_{j-1}))}{p^{2}}=\mathbb{P}(x_{i-1}\text{ in }\zeta(x_{j-1})),$$ 
    where $\mathbb{P}(x_{i-1}\text{ in }\zeta(x_{j-1}))$ is the proportion of appearances of $x_{i-1}$ in $\zeta(x_{j-1})$. For $j\in \{0,1,\dots,n-1\}$ and $k\geq 1$, denote $\mathbf{N}(\zeta^k(x_j))=(N_{x_0}(\zeta^k(x_j)),\dots,N_{x_{n-1}}(\zeta^k(x_j)))^T$, and let $\mathbb{P}(\zeta^k(x_j))=\frac{1}{p^{2k}}N(\zeta^k(x_j))$ be the probability distribution vector of $\zeta(x_j)$. Note that $M^T$ is a \textit{stochastic matrix}, that is the sum of the elements in every row is $1$.
\end{definition}
\begin{remark}
    The motivation behind the definition of $M$ is the following fact: for every $a,b\in \boldsymbol{X}$, 
    \begin{equation}
    \label{eqn:Nb(zeta^k(a))}
        N_b(\zeta^k(a))=\sum_{x\in \boldsymbol{X}}N_x(\zeta^{k-1}(a))N_b(\zeta(a)).
    \end{equation}
    Hence, for every $x\in \boldsymbol{X}$, 
    \begin{equation}\nonumber
        \mathbb{P}(\zeta^k(x))=\frac{1}{p^{2k}}\mathbf{N}(\zeta^k(x))=\frac{1}{p^2}M^{k-1}\mathbf{N}\left(\zeta(x)\right)=M^{k-1}\mathbb{P}(\zeta(x)).
    \end{equation}
    Thus, if $A$ is the initial symbol, and $x_0=0$, then $\mathbb{P}(0 \text{ in }\zeta^k(A))$ is equal to the first coordinate of $M^{k-1}\mathbb{P}(\zeta(A))$.
\end{remark}
\noindent The methods to prove Theorem \ref{thm:GenEsc=1} and the morphism from \cite{AN} are used to prove Theorem \ref{thm:TMFullEsc}. 
\begin{definition}
    The Thue Morse sequence $\{\tau_n\}_{n\in \mathbb{N}}$ is defined by the morphism $\mu:\{0,1\}\rightarrow \{0,1\}^2$ given by 
    $$\mu(0)=01\, \text{ and } \mu(1)=10.$$ 
    Let $\tau=\sum_{n=1}^{\infty}\tau_nt^{-n}\in\F_2(\!(t^{-1})\!)$. It is known that $\tau$ is a quadratic irrational \cite{APWW}.
\end{definition}
\noindent In \cite{AN}, the first named author, Erez Nesharim, and Uri Shapira found a morphism and coding which generates the rotated number wall (which is called the diagonally aligned number wall) of $\tau$. Explicitly, the diagonally aligned number wall of a sequence $\Theta=\{\Theta_n\}_{n\in \mathbb{N}}$ is defined by
$$\mathfrak{F}_{m,n}(\Theta)=\begin{cases}
    0&n>m+1\\
    1&n=m+1\\
    0&m=n\mod 2\\
    \det H_{\Theta}\left(\frac{m-n}{2},\frac{m+n}{2}\right)&m>n\text{ and }m\neq n\mod 2
\end{cases},$$

\noindent where $H_{\Theta}(m,n)$ was defined in Definition \ref{Hank}. Essentially, $\mathfrak{F}_{m,n}(\Theta)$ is a rotation by $\ang{45}$ of $W(\Theta)$. Hence, to evaluate escape of mass in the diagonally aligned number wall, one computes escape of mass along the columns, instead of along the diagonals (see \cite{AN} for details). Explicitly, the following lemma is an immediate consequence of Corollary \ref{cor:EscWall} and the definition of a diagonally aligned number wall.
\begin{lemma}
    Let $\Theta(t)=\sum_{i=1}^{\infty}\theta_it^{-i}$ be a quadratic irrational and let $t^k\Theta(t)=[0;\overline{a_{k,1},\dots,a_{k,\ell_k}}]$. Let $d_{k,j}$ be the $j^\nth$ non-zero coordinate in the $k^\nth$ diagonal of $\mathfrak{F}(\Theta)$. Then, $\Theta(t)$ exhibits $c$-escape of mass if and only if
        $$\lim_{n\rightarrow \infty}\liminf_{k\rightarrow \infty}\frac{\sum_{j=1,\dots,\ell_k-1}\max\{d_{k,j+1}-d_{k,j}-n,0\}}{\sum_{j=1}^{\ell_k-1}(d_{k,j+1}-d_{k,j})}\geq c.$$
\end{lemma}

\noindent  In \cite{AN}, the authors proved that $\tau$ exhibits $\frac{2}{3}$ escape of mass and moreover, found a sequence $\{t^{n_k}\tau\}$ which exhibits full escape of mass. To do so, they used the following $[2,2]$-morphism to generate $\mathfrak{F}(\tau)$.\\

\noindent Let $\Sigma=\{\mathbf{a}_i,\mathbf{b}_i,\mathbf{c}_i:0\leq i\leq 3\}\cup\{\mathbf{d_a},\mathbf{d_b},\mathbf{o}\}$ and let 
$$G=\langle \rho,\eta,\iota|\rho^4,\eta^2,\iota^2,\rho\eta\rho\eta,\rho\iota\rho\iota,\eta\iota\eta\iota\rangle=D_4\times\mathbb{Z}/2\mathbb{Z}=\langle \rho,\eta|\rho^2,\eta^2,\rho\eta\rho\eta\rangle\times \langle \iota|\iota^2\rangle$$
be the product of the dihedral group of order $8$ by a group of order two. For an integer $i$ and an integer $n$, denote $[i]_n=i\mod n$. Then $D_4$ acts on $\Sigma$ in the following manner:
\begin{equation}\label{eq:G-action}
	\begin{aligned}
		\rho(\mathbf{o})&=\mathbf{o}\, ,&\eta(\mathbf{o})&=\mathbf{o},\\
        \rho(\mathbf{a}_i)&=\mathbf{a}_{[i+1]_4}\,,&\text{for }i&=0,1,2,3\\
        \eta(\mathbf{a}_0)&=\mathbf{a}_0\,,&\eta(\mathbf{a}_2)&=\mathbf{a}_2\,,&\eta(\mathbf{a}_1)&=\mathbf{a}_3\,,&\eta(\mathbf{a}_3)&=\mathbf{a}_1,\\
        \rho(\mathbf{b}_i)&=\mathbf{b}_{[i+1]_4}&\text{for }i&=0,1,2,3\\
        \eta(\mathbf{b}_0)&=\mathbf{b}_0\,,&\eta(\mathbf{b}_2)&=\mathbf{b}_2\,,&\eta(\mathbf{b}_1)&=\mathbf{b}_3\\
        \rho(\mathbf{c}_i)&=\mathbf{c}_{[i+1]_4}\,,&\text{ for }i&=0,1,2,3,\\
        \eta(\mathbf{c}_0)&=\mathbf{c}_0\,,&\eta(\mathbf{c}_1)&=\mathbf{c}_1\,,&\eta(\mathbf{c}_1)&=\mathbf{c}_3.
	\end{aligned}
\end{equation}
Moreover, $\iota$ is the operation which switches between $\mathbf{a}$ and $\mathbf{b}$. Explicitly, $\mathbb{Z}/2\mathbb{Z}$ acts on $\Sigma$ in the following manner:
\begin{equation}
    \begin{aligned}
    \iota(\mathbf{x})=\begin{cases}
        \mathbf{d}_{\mathbf{b}}&\mathbf{x}=\mathbf{d}_\mathbf{a}\\
        \mathbf{d}_{\mathbf{a}}&\mathbf{x}=\mathbf{d}_{\mathbf{b}}\\
        \mathbf{x}&\text{else}
    \end{cases}.
    \end{aligned}
\end{equation}
Note that $\Sigma$ is a union of exactly four $G$ orbits whose representatives are $\{\mathbf{o},\mathbf{a}_0,\mathbf{c}_0,\mathbf{d_a}\}$.

\noindent The $[2,2]$-morphism from \cite{AN} is the unique $\rho$, $\iota$, and $\eta$ equivariant $[2,2]$-morphism defined by $\sigma:\Sigma\rightarrow \Sigma^{2\times 2}$ by
\begin{equation}
\begin{aligned}
\label{eqn:sigmaDef}
    \sigma(\mathbf{o})=\begin{matrix}
        \mathbf{o}&\mathbf{o}\\
        \mathbf{o}&\mathbf{o}
    \end{matrix},\quad&&\sigma(\mathbf{a}_0)=\begin{matrix}
        \mathbf{a}_0&\mathbf{c}_0\\
        \mathbf{d_a}&\mathbf{b}_0
    \end{matrix},\\
    \sigma(\mathbf{d_a})=\begin{matrix}
        \mathbf{b}_1&\mathbf{a}_2\\
        \mathbf{a}_0&\mathbf{b}_3
    \end{matrix},\quad&&
    \sigma(\mathbf{c}_0)=\begin{matrix}
        \mathbf{o}&\mathbf{o}\\
        \mathbf{c}_0&\mathbf{o}
    \end{matrix},
\end{aligned}
\end{equation}
where for every $M=\begin{matrix}
    a&b\\
    c&d
\end{matrix}\in \Sigma^{2\times 2}$, $\rho,\eta,$ and $\iota$ act on $M$ by
\begin{equation}\label{eq:G-of-matrices}
	\begin{aligned}
		\rho (M)&:=\begin{array}{cc}
			\iota(\rho(c)) &\iota(\rho(a)) \\
			\iota(\rho(d)) &\iota(\rho(b))
		\end{array},\quad &\eta (M)&:=\begin{array}{cc}
		\iota(\eta(d))&\iota(\eta(b))\\
		\iota(\eta(c))&\iota(\eta(a))
		\end{array},\quad &\iota (M)&:=\begin{array}{cc}
		\iota(a)&\iota(b)\\
		\iota(c)&\iota(d)
		\end{array}.
	\end{aligned}
\end{equation}
The $\rho$ equivariant $[4,4]$-coding $\kappa:\Sigma\rightarrow \{0,1\}^{4\times 4}$ is defined by
\begin{equation}
\begin{aligned}
    \kappa(\mathbf{o})&=\begin{matrix}
        0&0&0&0\\
        0&0&0&0\\
        0&0&0&0\\
        0&0&0&0
    \end{matrix}\,,&\kappa(\mathbf{c}_0)&=\begin{matrix}
        0&0&0&0\\
        0&0&0&0\\
        0&0&0&0\\
        1&0&0&0
    \end{matrix}\\
    \kappa(\mathbf{d_a})&=\begin{matrix}
        0&1&0&0\\
        0&0&1&0\\
        0&1&0&1\\
        1&0&1&0
    \end{matrix}\,,&\kappa(\mathbf{d_b})&=\begin{matrix}
        0&0&0&1\\
        1&0&1&0\\
        0&1&0&1\\
        0&0&1&0
    \end{matrix}\\
    \kappa(\mathbf{a}_0)&=\begin{matrix}
        0&1&0&0\\
        0&0&1&0\\
        0&1&0&1\\
        1&0&1&0
    \end{matrix}\,,&\kappa(\mathbf{b}_0)&=\begin{matrix}
        0&1&0&0\\
        1&0&1&0\\
        0&0&0&1\\
        0&0&0&0
    \end{matrix},
\end{aligned}
\end{equation}
where $\rho$ is the action which rotates by $\ang{90}$ clockwise. Let $\sigma^{\infty}(\mathbf{a}_0)=\lim_{n\rightarrow \infty}\sigma^n(\mathbf{a}_0)$. In this case, the symbol $\mathbf{o}$ plays the role of the symbol $0$ from Theorem \ref{thm:Til-->Esc}. By \cite[Theorem 1.13]{AN}, $\kappa(\sigma^{\infty}(\mathbf{d_a}))$ is the diagonally aligned number wall of $\tau$. 

\begin{proof}[Proof of Theorem \ref{thm:TMFullEsc}]
Firstly, one can obtain the conclusion of Theorem \ref{thm:Til-->Esc} holds for diagonally aligned number walls, since these are essentially number walls rotated by $45$ degrees and as Lebesgue measure is rotation invariant. Note that for every $a\in \Sigma$, the number of zeroes in $\kappa(a)$ is $G$ invariant, that is the number of zeroes in $\kappa(ga)$ is equal to the number of zeroes in $\kappa(a)$ for every $g\in G$.

\noindent Hence, $a$ can be identified with $\iota(a)$, $\rho(a)$, and $\eta(a)$ without changing the number of zeroes in $\kappa(\sigma^k(\mathbf{d_a}))$. Since this identification does not effect escape of mass, it suffices to work with a simplified morphism $\sigma'$ defined on the alphabet $\{\mathbf{o},\mathbf{c}_0,\mathbf{a}_0,\mathbf{d}_{\mathbf{a}}\}$ by \eqref{eqn:sigmaDef}. The transition matrix of $\sigma$ with respect to this alphabet is
\begin{equation}
    M=\frac{1}{4}\begin{pmatrix}
        4&3&0&0\\
        0&1&1&0\\
        0&0&2&4\\
        0&0&1&0
    \end{pmatrix}.
\end{equation}
Note that the eigenvalues of $M$ are $\left\{1,\frac{1}{4},\frac{1-\sqrt{5}}{4},\frac{1+\sqrt{5}}{4}\right\}$ with respective eigenvectors
\begin{equation}
    \begin{pmatrix}
        1\\
        0\\
        0\\
        0
    \end{pmatrix},\begin{pmatrix}
        -1\\
        1\\
        0\\
        0
    \end{pmatrix},\begin{pmatrix}
        -3+\frac{6\sqrt{5}}{5}\\
        1+\frac{\sqrt{5}}{5}\\
        1+\frac{\sqrt{5}}{5}\\
        1
    \end{pmatrix},\begin{pmatrix}
        -3-\frac{6\sqrt{5}}{5}\\
        1-\frac{\sqrt{5}}{5}\\
        1-\frac{\sqrt{5}}{5}\\
        1
    \end{pmatrix}.
\end{equation}
Note that by decomposing $\mathbb{P}(\sigma(\mathbf{d}_{\mathbf{a}}))$ as a linear combination of eigenvectors, one has
\begin{equation}
    \frac{1}{4}\begin{pmatrix}
        0\\
        1\\
        2\\
        1
    \end{pmatrix}=\begin{pmatrix}
        1\\
        0\\
        0\\
        0
    \end{pmatrix}-\frac{1}{20}\begin{pmatrix}
        -1\\
        1\\
        0\\
        0
    \end{pmatrix}+\left(\frac{1}{8}-\frac{\sqrt{5}}{40}\right)\begin{pmatrix}
        -3+\frac{6\sqrt{5}}{5}\\
        1+\frac{\sqrt{5}}{5}\\
        1+\frac{\sqrt{5}}{5}\\
        1
    \end{pmatrix}+\left(\frac{1}{8}+\frac{\sqrt{5}}{40}\right)\begin{pmatrix}
        -3-\frac{6\sqrt{5}}{5}\\
        1-\frac{\sqrt{5}}{5}\\
        1-\frac{\sqrt{5}}{5}\\
        1
    \end{pmatrix}.
\end{equation}
Hence,
\begin{equation}
\begin{split}
    \mathbb{P}\left(\mathbf{o}\text{ in }\sigma^k(\mathbf{d_a})\right)=1+\frac{1}{20\cdot 4^k}-\left(\frac{1-\sqrt{5}}{4}\right)^k\left(\frac{21}{40}-\frac{9\sqrt{5}}{40}\right)-\left(\frac{1+\sqrt{5}}{4}\right)^k\left(\frac{1}{8}+\frac{\sqrt{5}}{40}\right)\\=1-O(4^{-k})\rightarrow 1
\end{split}
\end{equation}
as $k\rightarrow \infty$. Thus, by Theorem \ref{thm:Til-->Esc}, for every $\varepsilon>0$, 
\begin{equation}
    \lim_{N\rightarrow \infty}d(k\in \mathbb{N}:e_{k,N}(\tau)>1-\varepsilon)=1.
\end{equation}
\end{proof}

\bibliography{Ref}
\bibliographystyle{amsalpha}
\end{document}